\documentclass[11pt]{article}

\usepackage[text={153mm,235mm},centering]{geometry}
\usepackage{enumerate}
\usepackage[mathscr]{eucal}
\usepackage[T1]{fontenc}
\usepackage[utf8]{inputenc}
\usepackage{lmodern}
\usepackage{amsmath,amssymb,amsthm,mathtools}
\usepackage{booktabs}
\usepackage{mathrsfs}
\usepackage{graphicx}
\usepackage[all]{xy}
\usepackage{enumitem}
\usepackage{hyperref}
\numberwithin{equation}{section}
\usepackage{tikz-cd}
\usepackage{stmaryrd}
\usepackage{microtype}
\usepackage{authblk}

\allowdisplaybreaks
\hypersetup{
  colorlinks=true,
  linkcolor=blue,
  citecolor=blue,
  urlcolor=blue
}

\newtheorem{theorem}{Theorem}[section]
\newtheorem{proposition}[theorem]{Proposition}
\newtheorem{corollary}[theorem]{Corollary}
\newtheorem{lemma}[theorem]{Lemma}

\theoremstyle{definition}
\newtheorem{definition}[theorem]{Definition}
\newtheorem{remark}[theorem]{Remark}
\newtheorem{example}[theorem]{Example}

\newtheorem*{ack}{Acknowledgements}

\newcommand{\thick}{\mathrm{thick}}

\newcommand{\id}{\mathrm{id}}

\title{An explicit Dennis trace and \\its compatibility with Koszul duality
\footnotetext{Email: huang.maozhou@outlook.com}}

\author{Maozhou Huang}

\affil{Department of Mathematics, Chongqing University of Technology, \\Chongqing City, 401320 P.R.\ China}

\date{}

\begin{document}

\maketitle

\begin{abstract}
In this paper, we construct explicit Dennis trace maps from algebraic $K$-theory to Hochschild homology for dg module and curved dg comodule categories, and prove their compatibility with derived Koszul duality. 
For this, we construct a natural chain map $\tau$ from the chains on the newly defined thickened nerve to the total complex of the dg cyclic nerve, with explicit formulas incorporating morphisms and coherent higher homotopies, which generalizes McCarthy's approach.  
This $\tau$ construction works both for modules over dg algebras and for cdg comodules over counital conilpotent curved dg coalgebras. 
Combining the $\tau$ with the dg-enriched Waldhausen $S_\bullet$-construction and a Hochschild delooping theorem yields the Dennis trace maps. 
For thick subcategories corresponding under derived Koszul duality, we prove that the induced $K$-theory equivalences and Hochschild homology isomorphisms fit into natural commutative diagrams with these trace maps.



\end{abstract}

\section{Introduction}\label{sect:intro}

The Dennis trace map relates algebraic $K$-theory to Hochschild homology. This paper constructs an explicit dg version of McCarthy's construction, in which morphisms and their higher homotopies enter the trace formula.

Let $k$ be a field of characteristic zero and let $\mathcal P$ be the category of finitely generated projective modules over an ordinary unital $k$-algebra $B$. McCarthy \cite[Definition~4.4.1]{McCarthy} defines the Dennis trace map by realizing and looping the map
\[\mathrm{obj}\,S_\bullet\mathcal P\longrightarrow\mathrm{CN}_0(S_\bullet\mathcal P)\longrightarrow\mathrm{CN}_\bullet(S_\bullet\mathcal P),\quad M\longmapsto\id_M,\]
where $S_\bullet\mathcal P$ denotes Waldhausen's $S_\bullet$-construction and $\mathrm{CN}_\bullet$ denotes the cyclic nerve. The assignment $X\mapsto\id_X$ also appears in Tabuada's comparison of the dg-categorical construction with the classical Dennis trace \cite[Theorem~2.8 and Section~7]{Tabuada}.

In the above expression, the information carried by the Dennis trace on morphisms remains in the background. This information is made more explicit by McCarthy's alternative description \cite[Section~4.5]{McCarthy}. Let $\mathrm i\mathcal P$ be the subcategory with the same objects as $\mathcal P$ and only isomorphisms. The cyclic structure on its nerve gives the following map of cyclic sets:
\begin{equation}
\begin{aligned}
\tau:N_\bullet \mathrm i\mathcal P&\longrightarrow\mathrm{CN}_\bullet(\mathcal P),\\
\bigl(X_1\xleftarrow{g_1}\cdots X_n\xleftarrow{g_n}X_0\bigr)
&\longmapsto(g_1\cdots g_n)^{-1}\otimes g_1\otimes\cdots\otimes g_n.
\end{aligned}
\label{formula:originalDennistrace}
\end{equation}
Applying the $S_\bullet$-construction, realizing, and looping recovers the Dennis trace $\mathrm{tr}$. 
Moreover, this gives a diagram commutative up to homotopy,
\begin{equation}
\vcenter{\xymatrix@=1.5em{|N_\bullet\mathrm{i}\mathcal P|\ar[r]\ar[d]_{\tau} & \Omega|N_\bullet \mathrm i S_\bullet \mathcal P|\ar[d]^{\mathrm{tr} = \Omega\tau_{S}}\\
|\mathrm{CN}_\bullet\mathcal P|\ar[r]^{\simeq\hspace{.5cm}} & \Omega|\mathrm{CN}_\bullet S_\bullet \mathcal P|.}}
\label{formula:McCarthyDennistracediagram}
\end{equation}

Trace maps for dg categories already fit into the theory of additive invariants \cite{Tabuada}. In the spectral setting, Blumberg-Gepner-Tabuada characterize the topological Dennis trace by a universal property \cite[Corollary~1.13]{BGT}, and Campbell-Lind-Malkiewich-Ponto-Zakharevich give a construction of the topological Dennis trace for spectral Waldhausen categories \cite{CLMPZ}. Our aim here is to give a dg refinement of the morphism formula~(\ref{formula:originalDennistrace}), together with its higher homotopies, and to obtain a natural dg analogue of diagram~(\ref{formula:McCarthyDennistracediagram}).

When $B$ is replaced by a dg algebra $A$, its $K$-theory is modeled by the full subcategory $\mathcal M$ of cofibrant-fibrant dg $A$-modules whose images in the derived category belong to the thick full subcategory $\thick_A(A)$ generated by $A$; see, e.g., \cite[Section~5.2]{Keller} and \cite[Definition~4.4]{BM1}. 
The subcategory $\mathrm w\mathcal M$ of weak equivalences plays the role of $\mathrm i\mathcal P$. 

However, weak equivalences admit homotopy inverses, but these inverses need not be strict. 
McCarthy's formula relies on strict inverse and composition identities to intertwine the simplicial and Hochschild boundaries. 
Simply substituting homotopy inverses therefore does not generally give a chain map: the required identities hold
only up to internal boundaries. 
A chain-level extension must incorporate the homotopies witnessing these identities, together with their higher coherence data.

We retain this information by replacing the ordinary nerve with the thickened nerve $\mathbb N_\bullet w\underline{\mathcal M}$ of Section~\ref{sec:thickenednerve}. 
Its simplices record weak equivalences in both directions, homotopies witnessing the inverse and composition identities, and coherent higher homotopies. 
These homotopies are part of the simplex data and are used directly in the trace formula.
Applied degreewise to the dg $S_\bullet$-construction, this replacement preserves the Waldhausen $K$-theory space by Proposition~\ref{prop:comparison}.

Our main new construction is an explicit natural chain map
\[\tau_{\underline{\mathcal M}}: \widetilde{k}[\mathbb N_\bullet w\underline{\mathcal M}] \longrightarrow \underline{\mathrm C}_\bullet(\underline{\mathcal M}),\] from reduced unnormalized simplicial chains to the unnormalized direct-sum total complex of the dg cyclic nerve. 
The construction extends McCarthy's formula by explicit correction terms built from the coherent homotopy data.
More precisely, for a thickened $n$-simplex $F$, let $T_0(F)$ denote the na\"ive tensor formed from its degree-zero edges. 
Although $T_0(F)$ is internally closed, its Hochschild boundary need not equal $T_0(\partial F)$.
The degree-$(-1)$ homotopies recorded by $F$ supply a correction $T_1(F)$ satisfying
\[ d_{\mathrm{int}}T_1(F) = T_0(\partial F)-\bar b\,T_0(F), \]
where $\bar b$ is the Hochschild differential with the conventions of Section~\ref{sec:tau}. 
The higher homotopies supply further corrections $T_2(F),\ldots,T_n(F)$, each correcting the remaining defect in the next bidegree. 
Lemma~\ref{lemma:tauisachainmap} proves that these corrections assemble into the finite sum
\[\tau_{\underline{\mathcal M}}(F)=\sum_{r=0}^{n}T_r(F), \quad D\tau_{\underline{\mathcal M}}(F) = \tau_{\underline{\mathcal M}}(\partial F).\]
Thus the homotopies recorded by the thickened nerve explicitly cancel the failure of the na\"ive formula to define a chain map.

When the dg category is concentrated in degree zero and its weak equivalences are isomorphisms, all higher correction terms vanish and the construction reduces to McCarthy's map~(\ref{formula:originalDennistrace}). 
More generally, the leading term induces a morphism of cyclic $k$-modules after taking internal cohomology in degree zero, as shown in Corollary~\ref{Cor:cyclicmap}.

Together with the delooping comparison, this chain map yields the following dg version of McCarthy’s diagram.
The bottom is the delooping theorem for cyclic nerves. Here we use the dg-categorical framework and the derived Morita theory described by Keller \cite{Keller}.
\begin{theorem}[Theorem~\ref{thm:Dennistracediagram}]
Let $C$ be a counital conilpotent cdg $k$-coalgebra, $A$ a unital dg $k$-algebra, $M$ a left cdg $C$-comodule, and $N$ a left dg $A$-module.
For the dg Waldhausen categories $\underline{\mathcal M}$ corresponding to the thick full subcategories $\mathrm{thick}_C(M)$ and $\mathrm{thick}_A(N)$, explained later in the introduction; 
let $\mathcal H(\underline{\mathcal M})$ and $\mathcal H_S(\underline{\mathcal M})$ be the homology spaces associated, respectively, to $\underline{\mathrm C}_\bullet(\underline{\mathcal M})$ and the total complex of $\underline{\mathrm C}_\bullet(\underline{S}_\bullet\underline{\mathcal M})$, see Section~\text{\rm \ref{sec:explicitDennistraceconstruction}}. 
There is a diagram
\begin{equation}
\xymatrix@=1.5em{
|\Bbb{N}_\bullet\underline{\mathrm w\mathcal M}|\ar[rr]^{\eta_{\underline{\mathcal M}}\hspace{.5cm}}\ar[d]_{\tau^t_{\underline{\mathcal M}}}
&&\Omega|\Bbb{N}_\bullet\underline{\mathrm wS}_\bullet\underline{\mathcal M}|\ar[d]^{\Omega\tau^t_{S_\bullet\underline{\mathcal M}}}\\
\mathcal H(\underline{\mathcal M})\ar[rr]^{\beta_{\underline{\mathcal M}}\hspace{.5cm}}_{\simeq
\hspace{.5cm}}&&\Omega\mathcal H_S(\underline{\mathcal M}),}\nonumber
\end{equation}
commutative up to based homotopy and natural in the dg exact functors under consideration. The bottom map is a weak homotopy equivalence.
\end{theorem}
The upper-right corner is a model of the $K$-theory $K(\mathcal M)$ of $\mathcal M$. 
We define the Dennis trace to be the right vertical map, retaining $\Omega\mathcal H_S(\underline{\mathcal M})$ as its target, as in Definition~\ref{def:dgDennistracemap}. 
By the delooping of Lemma~\ref{lemma:delooping},
\[\pi_n\Omega\mathcal H_S(\underline{\mathcal M})\cong\underline{\mathrm{HH}}_n(\underline{\mathcal M})\text{ for all $n\geq 0$}.\]
Thus taking homotopy groups gives $\mathrm{tr}_{\underline{\mathcal M}}:K_n(\mathcal M)\to\underline{\mathrm{HH}}_n(\underline{\mathcal M})$.

A second motivation of the current paper is to include $K$-theory for cdg coalgebras and its relation to the $K$-theory of dg algebras. 
Positselski's coderived categories and model structures provide the setting for treating cdg comodules \cite[Sections~4, 6, and~8]{Posit}. 
For a cdg $C$-comodule $M$ and a dg $A$-module $N$, we use the dg Waldhausen models $\underline{\mathcal M}_C(M)$ and $\underline{\mathcal M}_A(N)$ of the thick subcategories $\thick_C(M)$ and $\thick_A(N)$ from Sections~\ref{section:algebraicKfordgmodulesandcdgcomodules} and \text{\rm \ref{sec:dgK-theory}}. 
On the coalgebra side, the thick subcategory is taken in the coderived category. This gives $K$-theory and a Dennis trace for these categories of cdg comodules, together with their counterparts for dg modules.
Related work on coalgebraic $K$-theory includes Klanderman-P\'eroux’s constructions of Hattori-Stallings (co)traces involving coHochschild homology \cite{KP}, and Gerhardt-Péroux-Sor\'e’s comparisons between the $K$-theories of algebras and coalgebras compatible with these traces \cite{GPS}.

The choice $N=A$ recovers the usual $K$-theory of $A$. When $C$ is itself a cdg $C$-comodule, one can similarly use $\thick_C(C)$; in particular, this is available when $C$ is dg. 
For a general cdg coalgebra, the construction is attached to a chosen cdg comodule $M$. 

Derived Koszul duality supplies that comparison through the twisting-tensor functors. 
The role of thick subcategories in $K$-theoretic Koszul duality is already visible in the work of Blumberg-Mandell \cite{BM1}. 
In our setting, Propositions~\ref{prop:derivedKoszulequivalence} and~\ref{prop:KKoszulequivalence} relate the module and comodule models, and the naturality of the explicit trace extends these equivalences to the following compatibility statement, proved in Theorem~\ref{thm:compatibility}.

\begin{theorem}
Under the assumptions of Proposition~\text{\rm{\ref{prop:derivedKoszulequivalence}}}, let
\[
(\underline{\mathcal P},\underline{\mathcal Q})=
\begin{cases}
\bigl(\underline{\mathcal M}_C(k),\underline{\mathcal M}_A(A)\bigr),
&\text{without further conditions};\\[5pt]
\bigl(\underline{\mathcal M}_C(C),\underline{\mathcal M}_A(k)\bigr),
&\left(\begin{split}&\small \text{if $C$ is dg and $A$ has a dg augmentation}\\[-1ex] 
&\small \text{$v:A\to k$ satisfying $v\circ \tau = 0$}\end{split}\right).
\end{cases}
\]
The functor $L=A\otimes^\tau-:\underline{\mathcal P}\to\underline{\mathcal Q}$ induces a square of $K$-theory and Hochschild homology spaces, commutative up to based homotopy, whose horizontal maps are weak equivalences and whose vertical maps are the Dennis trace maps. 
In particular, for every $n\geq0$, there is a commutative diagram
\begin{equation}
\vcenter{\xymatrix@=1.5em{K_\bullet(\mathcal P)\ar[rr]^{K_\bullet(L)}_{\cong}\ar[d]_{\mathrm{tr}_{\underline{\mathcal P}}} && K_\bullet(\mathcal Q)\ar[d]^{\mathrm{tr}_{\underline{\mathcal Q}}}\\
\mathrm{HH}_\bullet(\underline{\mathcal P})\ar[rr]^{\mathrm{HH}_\bullet(L)}_{\cong} && \mathrm{HH}_\bullet(\underline{\mathcal Q}).}}\label{formula:KKoszulhomotopygroup}\end{equation}
\end{theorem}

McCarthy’s work suggests several further questions for this explicit dg model. 
These include comparison with Segal’s group-completion construction \cite{Segal} and the construction of a Jones–Goodwillie Chern map with equally explicit higher homotopies. 
The former comparison remains to be clarified for this model; see Remark~\ref{remark:withoutSegal}. 
The latter requires additional compatibility with the cyclic structure beyond the internal-cohomology statement above.
The higher deloopings in Remark~\ref{remark:Higherdelooping} suggest that our trace may extend to a natural transformation of additive invariants. 
With the appropriate normalization, such an extension could identify it with the standard Dennis trace through the dg counterpart of Blumberg-Gepner-Tabuada’s universal characterization \cite{BGT}.
The same deloopings also provide a starting point for studying compatibility with products induced by suitable biexact dg functors. 
We leave these developments to future work.

The paper is organized as follows. 
Section~\ref{section:algebraicKfordgmodulesandcdgcomodules} introduces the Waldhausen categories associated with dg modules and cdg comodules and establishes the $K$-theory comparisons induced by derived Koszul equivalence. 
Section~\ref{sec:dgK-theory} develops the dg Waldhausen framework and the dg-enriched $S_\bullet$-construction. 
Section~\ref{sec:keyconstruction} introduces the thickened nerve, compares the resulting $K$-theory model with the usual one, and constructs the chain map $\tau$. 
Section~\ref{sec:explicitDennistrace} proves the delooping theorem, defines the Dennis trace map, and establishes its compatibility with derived Koszul equivalence. 
Appendices~\ref{appendix:chk-modelcategory} and~\ref{sec:Reedymodelcategories} provide the enriched model-category and Reedy arguments used in these constructions, while Appendix~\ref{appendix:cyclicnerve} treats the dg cyclic nerve and the homotopies needed for additivity.

\subsection*{Notations and assumptions}
Throughout the paper, $k$ denotes a field of characteristic zero. 
Dg categories are denoted by underlined letters. 
The dg nerve functor is denoted $\underline{N}_\bullet$.
For a dg category $\underline{\mathcal C}$, its underlying ordinary category is $\mathcal C=Z^0\underline{\mathcal C}$, with closed degree-zero morphisms; its homotopy category is denoted by $H^0(\underline{\mathcal C})$. In Section~4, the notation $\underline{\mathrm w\mathcal M}$ denotes the maximal Kan subcomplex of the dg nerve, and should be distinguished from a dg category. 
The corresponding convention applies to $\underline{\mathrm wS}_n\underline{\mathcal M}$.

For a dg module $M$ or a cdg comodule $N$, the underlying graded objects are denoted by $M^\#$ and $N^\#$. Let $\mathrm{Ch}(k)$ denote the category of unbounded cochain complexes. The homological grading of the Hochschild total complexes is the regrading fixed in Section~\ref{sec:tau}; all totalizations in Sections~\ref{sec:tau} and~\ref{sec:explicitDennistrace} use direct sums and the unnormalized complexes specified there. Finally, $\Delta$ denotes the category of nonempty finite ordinals $[n]=\{0,\ldots,n\}$, for $n\geq0$, and order-preserving maps.

\begin{ack}
The author acknowledges the assistance of OpenAI’s ChatGPT (GPT-5.0 to 6.0, and 6~Astra) in mathematical discussions and in developing and revising the manuscript. 
The initial idea of using the thick nerve originated with the author. 
The author takes full responsibility for the mathematical content, proofs, and references.
I thank Sichuan University and Chongqing University of Technology for their hospitality during my visits while this work was in progress. 
I am especially grateful to Professor Xiaojun Chen for his support, encouragement, helpful discussions, and valuable suggestions for improving the manuscript.
This project is supported by NSFC 12271377.
\end{ack}

\tableofcontents

\section{Algebraic $K$-theory for certain dg modules \\and cdg-comodules}\label{section:algebraicKfordgmodulesandcdgcomodules}

Our goal is to define $K$-theory for certain cofibrant-fibrant cdg-comodules of a counital conilpotent cdg-coalgebra.
This definition mimics \cite[Section~4.4]{BM2}, the $K$-theory for certain dg modules of a unital dg algebra.

\subsection{Preliminaries on cdg comodules}
We recall the definitions and results from \cite{Posit} needed below.
We use cohomological $\Bbb{Z}$-grading for the cdg and dg objects.
A \emph{cdg-coalgebra} is a triple $C = (C,d_C,h)$ with $C$ a graded coassociative counital coalgebra, $d_C$ a coderivation of degree $1$, and a homogeneous linear map $h:C\to k$ of degree $2$ satisfying $hd_C = 0$ and $d_C^2(c) = \sum h(c_{(1)})c_{(2)} - \sum c_{(1)}h(c_{(2)})$ writing $\Delta(c) = \sum c_{(1)}\otimes c_{(2)}$, where $h$ is called the curvature.
A cdg-coalgebra is \emph{conilpotent} if it is coaugmented $w:k\to C$ such that $d_Cw = 0$ and $hw=0$, and each $c$ in the kernel of the counit is killed by some iterate of the reduced coproduct $\overline{\Delta}(c) = \Delta(c) - w(1)\otimes c - c\otimes w(1).$
A (left) \emph{cdg-comodule} is a triple $M = (M,\lambda,d_M)$ with $M$ a graded left $C$-comodule, $d_M$ a homogeneous endomorphism of degree $1$, and $\lambda:M\to C\otimes M$ such that $\lambda d_M = (d_C\otimes \id_M + \id_C \otimes d_M)\lambda$ and $d_M^2(m) = \sum h(m_{(-1)})m_{(0)}$ with the coaction $m\mapsto m_{(-1)}\otimes m_{(0)}.$
If $h = 0$, then a cdg-coalgebra is a dg coalgebra and cdg-comodules are dg comodules.

Suppose throughout this paper that all coalgebras are counital conilpotent.
Let $\underline{\mathrm{CM}}_C$ be the dg category of left cdg-comodules over $C$ so that for any homogeneous $f\in \underline{\mathrm{CM}}_{C}(L,M)$, we have $df = d_M f - (-1)^{|f|}fd_L.$
The category $H^0\underline{\mathrm{CM}}_C$ has the same objects but has morphism spaces $H^0(\underline{\mathrm{CM}}_C(L,M))$ for any $L,M$.
Write $\mathrm{CM}_C = Z^0\underline{\mathrm{CM}}_C$ for the abelian category of left cdg comodules and closed degree-zero morphisms.
Next, consider the smallest full triangulated subcategory $\mathrm{Acycl}^{\mathrm{co}}\underline{\mathrm{CM}}_C$ of $H^0 \underline{\mathrm{CM}}_C$ such that it contains the totalizations $\mathrm{Tot}(L\to M\to N)$, which are cdg comodules, of all short exact sequences $0\to L\to M\to N\to 0$ in $\mathrm{CM}_C$, and is closed under arbitrary direct sums.
Then the \emph{coderived category} is defined as the following Verdier quotient
\[D^{\mathrm{co}}\underline{\mathrm{CM}}_C = H^0 \underline{\mathrm{CM}}_C/\mathrm{Acycl}^{\mathrm{co}}\underline{\mathrm{CM}}_C.\]

\begin{lemma}[\text{\rm{Positselski \cite[Theorem~4.4(c)]{Posit}}}]\label{lemma:Verdierquotientcategoryeq}
Let $\underline{\mathrm{CM}}_C^{\mathrm{inj}}$ denote the full dg subcategory consisting of those cdg comodules whose underlying graded comodule is injective.
The following composition of functors is an equivalence of triangulated categories
\[H^0\underline{\mathrm{CM}}_C^{\mathrm{inj}} \to H^0\underline{\mathrm{CM}}_C \to D^{\mathrm{co}}\underline{\mathrm{CM}}_C.\]
\end{lemma}

We need the following model category structure
\begin{proposition}[\text{\rm{Positselski \cite[Theorem~8.2]{Posit}}}]\label{prop:modelstructureforcdgcomodules}
There exists a model category structure on $\mathrm{CM}_C$:
\begin{enumerate}[label=\textup{(\alph*)},noitemsep,leftmargin=*,nosep]
    \item A morphism is a \emph{cofibration} if and only if it is injective. 
    \item A morphism is a \emph{fibration} if and only if it is surjective and the kernel is graded injective.
    \item A morphism is a \emph{weak equivalence} if and only if its cone is in the subcategory $\mathrm{Acycl}^{\mathrm{co}}\underline{\mathrm{CM}}_C.$
    \item Every object is cofibrant. An object is \emph{fibrant} if and only if it is graded injective.
\end{enumerate}
\end{proposition}

Let $(A,d_A)$ be a unital dg algebra and $\underline{\mathrm{M}}_A$ be the dg category of left dg-modules over $A$.
Put $\mathrm{M}_A \coloneqq Z^0\underline{\mathrm{M}}_A$.
This category admits the following projective model structure.
\begin{proposition}[\text{\rm{Positselski \cite[Theorem~8.1]{Posit}}}]\label{prop:projectivemodelstructureformodules}
There exists a model category structure on $\mathrm M_A$:
\begin{enumerate}[label=\textup{(\alph*)},noitemsep,leftmargin=*,nosep]
\item A morphism $f$ is a \emph{cofibration} if and only if it is injective and its cokernel $\mathrm{Coker}(f)$ is both graded projective and such that $\underline{\mathrm{M}}_A(\mathrm{Coker}(f),N)$ is acyclic for every acyclic dg $A$-module $N$;
\item A morphism is a \emph{fibration} if and only if it is surjective;
\item A morphism is a \emph{weak equivalence} if and only if it is a quasi-isomorphism;
\item Every object is fibrant. An object $L$ is \emph{cofibrant} if and only if it is both graded projective and such that $\underline{\mathrm{M}}_A(L,N)$ is acyclic for every acyclic dg $A$-module $N$.
\end{enumerate}
\end{proposition}

The Koszul functors reviewed below connect two categories $\underline{\mathrm M}_A$ and $\underline{\mathrm{CM}}_C$. 
Let $v:A\to k$ be a homogeneous $k$-linear retraction, that is, its composition with the unit map is the identity map of $k$.
The \emph{bar construction} $\mathrm{B}_v(A)= T^c(\mathrm{Ker}(v)[1])$ has a cdg-coalgebra structure: \begin{itemize}[noitemsep,leftmargin=*,nosep]
\item The degree-one coderivation $d_{\mathrm B}$ is given by
\[\textstyle d_{\mathrm{B}}([a_1|\cdots|a_n]) = \sum_{r=1}^{\min\{2,n\}}\sum_{i=0}^{n-r} (-1)^\epsilon[a_1|\cdots|p_r([a_{i+1}]| \cdots |[a_{i+r}])|\cdots |a_n],\]
where for the Koszul sign $(-1)^\epsilon$, $\epsilon = \sum_{j=1}^i (|a_j|-1)$ with the degree $|a_i|$ of $a_i$ in $A$ and
\[\begin{cases}p_1([a_1]) = -[d_{\mathrm{Ker}}(a_1)] & n = 1;\\
p_2([a_1| a_2]) = (-1)^{|a_1|-1}[m_{\mathrm{Ker}}(a_1\otimes a_2)] & n = 2;\\
p_n = 0& n\neq 1,2,\end{cases}\]
where $d_\mathrm{Ker}:\mathrm{Ker}(v)\to \mathrm{Ker}(v)$ and $m_\mathrm{Ker}:\mathrm{Ker}(v)^{\otimes 2}\to \mathrm{Ker}(v)$ are given by precomposing $A\to \mathrm{Ker}(v);\,\,a\mapsto a-v(a)$ with respectively the differential and the multiplication.
\item The curvature $h_{\mathrm{B}}: \mathrm{B}_v(A)\to k$ writes \[\begin{cases}
h_{\mathrm{B}}([a_1]) = -d_k(a_1) & n = 1;\\
h_{\mathrm{B}}([a_1|a_2]) = (-1)^{|a_1|-1}m_{k}(a_1\otimes a_2) & n = 2;\\
h_{\mathrm{B}}|_{\text{Ker}(v)[1]^{\otimes n}} = 0 & n \neq 1,2,
\end{cases}\]
where $d_k:\mathrm{Ker}(v)\to k$ and $m_k:\mathrm{Ker}(v)^{\otimes 2}\to k$ are given by precomposing $v$ with respectively the differential and the multiplication. Moreover, $h_B$ vanishes if and only if $v$ is a dg-algebra augmentation.
\end{itemize}
On the other hand, for a coaugmented cdg-coalgebra $C$ with $w$ being the coaugmentation, its \emph{cobar construction} $\mathrm{Cob}_w(C) = T(\mathrm{Coker}(w)[-1])$ has a dg-algebra structure. Indeed, one may dually define a degree-one derivation $d_\mathrm{Cob}$ which turns out to be a differential, and a curvature element which turns out to be zero (See \cite[p.\ 62]{Posit}).

\begin{lemma}\label{lemma:dg}
Let $(C,d_C,h)$ be a cdg-coalgebra. The triple $(C,\Delta,d_C)$ is a left cdg $C$-comodule iff $h_C = 0$.
For a dg-algebra $A$ and a homogeneous $k$-linear retraction $v$ of its unit, $v$ is a dg-algebra augmentation if and only if $h_{\mathrm{B}_v(A)} = 0$. 
\end{lemma}

For a coaugmented cdg-coalgebra $C$ and a dg-algebra $A$, a homogeneous linear map $\tau:C\to A$ of degree one satisfying $\tau^2+\tau\circ d + d\circ \tau - 1_Ah_C =0$ is called a \emph{twisting cochain}. 
Here the square $\tau^2$ is given by the convolution product on $\mathrm{Hom}_k(C,A)$ and $1_A$ denotes the unit of $A$.
A twisting cochain $\tau:C\to A$ such that $\tau\circ w = 0$ bijectively corresponds to a morphism of dg-algebras $f_\tau:\mathrm{Cob}_w(C)\to A$.
The twisting cochain $\tau$ satisfying $\tau\circ w = 0$ is called \emph{acyclic} if $f_\tau$ is a quasi-isomorphism.
For a left cdg-comodule $M$ over $C,$ the pair $A\otimes^\tau M = (A\otimes_k M,d)$ with $d(a\otimes x) = d_{A}a \otimes x + (-1)^{|a|}a\otimes d_Mx + (-1)^{|a|-1}\sum a\tau(x_{(-1)})\otimes x_{(0)}$ is a left dg-module over $A$.
For a left dg-module $N$ over $A$, the $3$-tuple $C\otimes^\tau N = (C\otimes_k N,\lambda,d)$ with the  coaction $\lambda(c\otimes y) = \sum c_{(1)}\otimes (c_{(2)}\otimes y)$ and $d(c\otimes y) = d_C(c) \otimes y + (-1)^{|c|}c\otimes d_N y + \sum (-1)^{|c_{(1)}|}c_{(1)}\otimes \tau(c_{(2)})y$ is a left cdg-comodule over $C$. 
\begin{proposition}[\text{\rm{Positselski \cite[Theorem~6.5, and Sections~6.2 and 8.4]{Posit}}}]\label{prop:derivedKoszulequivalence}
Let $C$ be a counital conilpotent cdg-coalgebra with the coaugmentation $w$ and $A$ a unital dg-algebra. 
Let $\tau$ be an acyclic twisting cochain with $\tau \circ w = 0$.
Then there are adjoint dg functors $A\otimes^\tau - \dashv_{\mathrm{dg}} C\otimes^{\tau} -$:
\[A\otimes^\tau - :\underline{\mathrm{CM}}_C\rightleftarrows \underline{\mathrm{M}}_A:C\otimes^{\tau} - .\]
On the underlying categories of closed degree-zero morphisms, they form a Quillen equivalence for the model structures of Propositions~\text{\rm \ref{prop:modelstructureforcdgcomodules}} and \text{\rm \ref{prop:projectivemodelstructureformodules}}.
Consequently, their total derived functors give mutually inverse equivalences of triangulated categories \[A\otimes_{\mathrm{L}}^\tau - :D^{\mathrm{co}}\mathrm{CM}_C\rightleftarrows D\mathrm{M}_A:C\otimes_\mathrm{R}^{\tau} - .\]
\end{proposition}
The following canonical twisting cochains satisfy the hypotheses of Proposition~\ref{prop:derivedKoszulequivalence}:
\begin{itemize}[noitemsep,leftmargin=*,nosep]
    \item the map $\tau$ being $(-1)$ times the composition $C\coloneqq\mathrm{B}_v(A)\to \mathrm{Ker}(v)[1]\to \mathrm{Ker}(v)\hookrightarrow A$;
    \item the map $\tau$ being the composition $C\to \mathrm{Coker}(w)\to \mathrm{Coker}(w)[-1]\to A\coloneqq\mathrm{Cob}_w(C)$,
\end{itemize}
where both middle arrows are the degree shift by $+1$;

\subsection{Definition of $K$-theory}

For a cdg-coalgebra $C$ and a left cdg $C$-comodule $M$, let $\mathrm{thick}_C(M)$ be the smallest full triangulated subcategory of $D^\mathrm{co}\underline{\mathrm{CM}}_C$ containing the image of $M$ and closed under retracts.
Let $\mathcal{M}_{C}(M)$ be the full subcategory of $\mathrm{CM}_C$ consisting of cofibrant-fibrant objects whose images in $D^{\mathrm{co}}\underline{\mathrm{CM}}_C$ belong to $\mathrm{thick}_{C}(M).$

\begin{lemma}\label{lemma:MCMisaWaldhausencategory}
For any left cdg $C$-comodule $M$, equip $\mathcal M_C(M)$ with the cofibrations and weak equivalences inherited from the model structure of $\mathrm{CM}_C$, and then $\mathcal M_{C}(M)$ is a Waldhausen category.\end{lemma}

\begin{proof}
The zero object is in $\mathcal M_C(M).$
Because of the model structure, it suffices to show: (1) closure under pushouts along cofibrations; and (2) the gluing axiom.

For (1), we show that the pushout $P\coloneqq Y\coprod_X Z$ is fibrant given cofibration $X\hookrightarrow Y$ and a map $X\to Z$ in $\mathcal M_C(M).$
Put $Q\coloneqq \mathrm{coker}(X\hookrightarrow Y).$
For the underlying graded comodules $X^\#$, $Y^\#$, and $Z^\#$ in the category of graded comodules, we have the exact sequence $X^\#\hookrightarrow Y^\# \twoheadrightarrow Q^\#$, and $X^\#$ is injective implies that this sequence is split, forgetting the differentials, and hence $Q^\#$ is injective. Similarly, the sequence $Z^\#\hookrightarrow P^\# \twoheadrightarrow Q^\#$ is split exact in the category of graded comodules, and hence $P$ is graded injective.
We also need to show that $P$ has image in $\mathrm{thick}_C(M).$
This follows from that $X\to Y\to Q\to X[1]$ and $Z\to P\to Q\to Z[1]$
are distinguished triangles in $D^\mathrm{co}\mathrm{CM}_C.$

For (2), with the notations above, suppose weak equivalences $X\to X'$, $Y\to Y'$, $Z\to Z'$, a cofibration $X'\to Y'$, a map $X'\to Z'$ and the corresponding commutativity.
The above maps induce the map between distinguished triangles $X\to Y\to Q\to X[1]$ and $X'\to Y'\to Q'\to X'[1]$, where $Q'= \mathrm{coker}(X'\to Y')$ and the map $Q\to Q'$ is a weak equivalence.
Similarly, for the induced map between $Z\to P\to Q\to Z[1]$ and $Z'\to P'\to Q'\to Z'[1],$ the map $P\to P'$ is a weak equivalence, as desired.\end{proof}

On the other hand, for a dg algebra $A$, and an $A$-module $N$ and the derived category $D\underline{\mathrm M}_A$, we may similarly define $\mathrm{thick}_{A}(N)$ and $\mathcal M_A(N),$ where the latter category consists of cofibrant objects and is defined using the projective model structure. The $K$-theory of $\mathrm{thick}_C(M)$ and $\mathrm{thick}_A(N)$ are defined by applying the $S_\bullet$-construction (See \cite[Section~2]{BM2} or \cite[Section II.9]{Weibel}) to the Waldhausen categories.

\begin{definition}[$K$-theory]
Let $C$ be a counital conilpotent cdg coalgebra over $k$ and $M$ a left cdg $C$-comodule.
The $K$-theory space $K(\mathcal M_C(M))$ of $\mathrm{thick}_C(M)$ is the pointed loop space of the $S_\bullet$-construction $\Omega |N_\bullet\mathrm{w}S_\bullet\mathcal{M}_C(M)|,$ where $\mathrm{w}S_\bullet\mathcal{M}_C(M)$ has the same objects but morphisms are weak equivalence only.
Let $A$ be a unital dg algebra over $k$ and $N$ a left dg $A$-module.
The $K$-theory space $K(\mathcal M_A(N))$ of $\mathrm{thick}_A(N)$ is the pointed loop space of the $S_\bullet$-construction $\Omega|N_\bullet\mathrm{w} S_\bullet\mathcal{M}_A(N)|.$
\end{definition}

\begin{theorem}[Waldhausen's approximation]\label{thm:approximation}
Let $F:\mathcal A\to \mathcal B$ be an exact functor between Waldhausen categories. 
Assume that weak equivalences in both categories satisfy the two-out-of-three property, and that $\mathcal A$ admits a cylinder functor satisfying the cylinder axiom.
Suppose moreover
\begin{enumerate}[label=\textup{(\alph*)},noitemsep,leftmargin=*,nosep]
\item A morphism $f$ in $\mathcal A$ is a weak equivalence if and only if $F(f)$ is a weak equivalence in $\mathcal B$;
\item For an object $X\in \mathcal A$ and every morphism $u:F(X)\to Y$ in $\mathcal B$, there exist a cofibration $i:X\hookrightarrow X'$ in $\mathcal A$ and a weak equivalence $v: F(X')\overset{\sim}{\to}Y$ in $\mathcal B$ such that $u = v\circ F(i).$
\end{enumerate}
Then $F$ induces a homotopy equivalence $K(\mathcal A)\overset{\simeq}{\to}K(\mathcal B).$
\end{theorem}

Following Blumberg-Mandell's idea in \cite[Remark~1.7]{BM1}, we may find a small category which may be treated as a subcategory of  $\mathcal{M}_C(M)$ but give the same $K$-theory, as below.

\begin{lemma}
Fix $C$ and $M$ as above. Let $\kappa$ be a cardinal $\geq \max\{\aleph_0,\mathrm{card}\,k,\mathrm{card}\,C,\mathrm{card}\,M\}$ and regard it as a set, where taking cardinal means that of the underlying set of the total graded vector space.
Let $\mathcal M_\kappa$ be the full subcategory of $\mathcal M_C(M)$ consisting of objects whose underlying sets are subsets of $\kappa$.
Equip this subcategory with the inherited cofibrations and weak equivalences.
Then $\mathcal M_\kappa$ is a small Waldhausen category, and the inclusion functor induces a weak homotopy equivalence $K(\mathcal M_\kappa) \overset{\simeq}{\to}K(\mathcal M_C(M))$.
\end{lemma}

\begin{proof}
The category $\mathcal M_\kappa$ is small: the possible cdg-comodule structures on subsets of $\kappa$ form a set, and its morphisms are functions between these underlying sets.
For a cofibration $A\hookrightarrow B$ and a morphism $A\to D$ in $\mathcal M_\kappa$, by Lemma~\ref{lemma:MCMisaWaldhausencategory} the pushout in the ambient category $P\in\mathcal M_C(M)$ writes $P = B\coprod_A D$. 
As $P$ is a quotient of $B\oplus D$, we have $\mathrm{card}\,P\leq \kappa.$
Choose a subset $S\subset \kappa$ and a bijection $e$ from $S$ to the underlying set $U(P)$ of $P$.
Transport the cdg-comodule structure and pushout map along $e$, obtaining an object $Q$ in $\mathcal M_\kappa$.
Hence $\mathcal M_\kappa $ inherits a Waldhausen category structure and the inclusion is exact.

It suffices to check the hypothesis of the Waldhausen approximation theorem.
The inclusion functor reflects weak equivalences, and the inherited weak equivalences in both categories satisfy the two-out-of-three property.
The mapping-cylinder construction for comodules, transported via the bijection of underlying sets contained in $\kappa$, gives a cylinder functor satisfying the cylinder axiom. 
It suffices to show the claim: for object $A$ in $\mathcal M_\kappa$ and a map $f:A\to X$ in $\mathcal M_{C}(M)$, there is a factorization $A\hookrightarrow A' \overset{\sim}{\to}X$ of $f$ as a composition of a weak equivalence in $\mathcal M_C(M)$ and a cofibration with $A'\in\mathcal M_\kappa.$

To construct $A'$, we first show that every object in the smallest full triangulated subcategory containing $M$ admits a cdg-comodule representative of cardinality at most $\kappa$.
Indeed, $M$ satisfies this bound. Shifts and finite direct sums preserve it. To check the closure under cones, we note that given $N\in \mathcal M_\kappa$, the fibrant replacement $J(N)$ whose underlying graded comodule writes $J(N)^\# = \bigoplus_{n\geq 0}\big(C^{\otimes_k (n+1)}\otimes_k N\big)[-n]$ constructed in \cite[Theorem~4.4]{Posit}, is fibrant and has cardinal at most $\kappa$, and the map $N\to J(N)$ has coacyclic cone.
Let $f: X\to Y$ be a morphism in the coderived category such that there are cdg-comodule representatives of cardinality at most $\kappa.$
Since the composition $X\to Y\to J(Y)$ is represented by a closed degree zero map $\tilde{f}:X\to J(Y)$, the mapping cone of $f$ is represented in the coderived category by that of $\tilde{f}.$
We have $\mathrm{card}\,\mathrm{cone}(\tilde{f}) = \mathrm{card}(J(Y)^\#\oplus X^\#[1])\leq \kappa,$ as desired.

We show: for any object $W$ in $\mathrm{thick}_C(M)$, there is a fibrant $C$-comodule $V$ satisfying $\mathrm{card}(V) \leq \kappa$ that is weak equivalent to $W$.
By the preceding claim, choose a cdg-comodule representative \(P\) of the object of which \(W\) is a retract, with \(\operatorname{card}(P)\leq\kappa\). 
Use the isomorphism $P\simeq J(P)$ in the coderived category to regard $W$ as a retract of $J(P).$
By Lemma~\ref{lemma:Verdierquotientcategoryeq}, the resulting idempotent is represented by a closed endomorphism \(e:J(P)\to J(P)\), with \(e^2\simeq e\).
Then $W$ is weak equivalent to the mapping telescope of $e$: $\textstyle\mathrm{tel}(e) = \mathrm{cone}\big(1-s_e:\bigoplus_{j\geq 0}J(P)\to \bigoplus_{j\geq 0}J(P) \big),$ where $s_e$ maps each summand to the next by $e$.
Put $V \coloneqq J(\mathrm{tel}(e))$.
As $V$ is fibrant, $\mathrm{card}(V)\leq \kappa,$ and represents $W$, the desired claim follows.

For $X$ above, choose \(Y\in\mathcal M_\kappa\) whose image in the coderived category is isomorphic to that of \(X\).
By Lemma~\ref{lemma:Verdierquotientcategoryeq}, we have weak equivalence $g:X\overset{\sim}{\rightarrow}Y$ and its homotopy inverse $h:Y\to X$. 
The two maps $hgf$ and $f$ from $A$ to $X$ are homotopic.
The homotopy gives the following map $H:\mathrm{cyl}(A)\to X$ such that for the two canonical map $i_0,i_1:A\to \mathrm{cyl}(A)$, we have $Hi_0 = hgf$ and $Hi_1 = f$.
Let $A'$ denote the following pushout.
\[\xymatrix@=1.5em{A\ar[r]^{i_0}\ar[d]_{gf} & \mathrm{cyl}(A)\ar[d]_{l}\ar[rdd]^H & \\
Y\ar[r]^{j}\ar[rrd]_h & A'\ar@{.>}[rd] & \\
&& X.}\]
The dotted arrow $q$ is the desired weak equivalence $A'\overset{\sim}{\to} X.$
Note that $j$ is a weak equivalence because it is the pushout of the trivial cofibration $i_0$.
As $h$ is also a weak equivalence, $q$ being a weak equivalence follows from the two-out-of-three property.
Because of $q$, the image of $A'$ in the derived category is in $\mathrm{thick}_C(M).$ 
The desired cofibration $A\hookrightarrow A'$ is the composition $l\circ i_1$ because $q\circ l \circ i_1 = H\circ i_1 = f.$
We check $\mathrm{card}\,A'^\#\leq \kappa$ and $A'$ is fibrant.
The cylinder object $\mathrm{cyl}(A)$ writes $A\oplus A\oplus A[1]$ and $A'$ writes $Y\oplus A\oplus A[1].$
Then $l\circ i_1$ maps $a$ in $A$ to $(0,a,0)$ in $Y\oplus A\oplus A[1]$. 
Transporting $A'$ and the factorization maps to an isomorphic copy with underlying set contained in $\kappa$ gives the desired factorization with $A'\in\mathcal M_\kappa.$ 
\end{proof}

Under the hypotheses of Proposition~\ref{prop:derivedKoszulequivalence}, we have the Koszul dual functors $A\otimes^\tau -$ and $C\otimes^\tau -$ induces a category equivalence $\mathrm{thick}_C(k)\simeq \mathrm{thick}_A(A).$ If moreover $C$ is dg and $A$ has a dg algebra augmentation $v:A\to k$ satisfying $v\circ \tau =0,$ then $\mathrm{thick}_{C}(C) \simeq \mathrm{thick}_{A}(k).$

Following Blumberg-Mandell \cite[Section~5]{BM2}, one may apply the derived Koszul equivalence and Waldhausen's approximation theorem to get equivalences of $K$-theories of algebras and coalgebras. 

\begin{proposition}\label{prop:KKoszulequivalence}
Under the hypotheses of Proposition~\ref{prop:derivedKoszulequivalence}, we have the following homotopy equivalence 
\[K(\mathcal M_C(k)) \simeq K(\mathcal M_A(A)),\]
and if $C$ is dg and $A$ has a dg algebra augmentation $v$ such that $v\circ \tau = 0$, then
\[K(\mathcal M_C(C)) \simeq K(\mathcal M_A(k)),\]
\end{proposition}

\begin{proof}[Sketch of proof]
We check the condition Theorem~\ref{thm:approximation}~(b). 
Put $F\coloneqq A\otimes^\tau - $ and $G \coloneqq C\otimes^\tau - $.
For any $u: F(X)\to Y$ of in $\mathrm{M}_A,$ the model structure gives factorization of its adjoint in $\mathrm{CM}_C$ into a cofibration and a trivial fibration $X\hookrightarrow X'\overset{p}{\underset{\sim}{\twoheadrightarrow}} G(Y).$
As $Y$ is fibrant and $G$ is right Quillen, $G(Y)$ and hence $X'$ are fibrant. 
Since the derived functor of $G$ identifies the corresponding thick subcategories and $p$ is a weak equivalence, $X'$ is in the required source Waldhausen category.
Combining the image of this composition and the counit gives the desired $F(X)\hookrightarrow F(X')\overset{\sim}{\twoheadrightarrow} FG(Y) \to Y$.
Indeed, the composition is $u$ by the adjointness.
For $q:F(X') \overset{\sim}{\twoheadrightarrow} FG(Y)\to Y$, its adjoint is the weak equivalence $p$. 
Since $X'$ is cofibrant, $Y$ is fibrant, and $F\dashv G$ is a Quillen equivalence, the composition $q$ is a weak equivalence.
\end{proof}

\section{$K$-theory of dg Waldhausen categories}\label{sec:dgK-theory}
Thereafter, let $\mathcal M$ denote the category $\mathcal M_A(N)$ or $\mathcal M_C(M)$ for a dg $A$-module $N$ and a cdg $C$-comodule $M$, with the Waldhausen structures and size conventions fixed in Section~\ref{section:algebraicKfordgmodulesandcdgcomodules}.
We note that $\mathcal M$ consists of cofibrant-fibrant objects.
Let $\underline{\mathcal M}$ be the full dg subcategory of $\underline{\mathrm M}_A$ or $\underline{\mathrm{CM}}_C$, respectively, on the objects of $\mathcal M$ so that $Z^0\underline{\mathcal M}=\mathcal M.$

Following the approach of \cite{CLMPZ}, we introduce dg Waldhausen categories and their enriched $S_\bullet$-construction. 
The compatibility between the dg enrichment and the Waldhausen structure in the definition below will be needed when the enriched $S_\bullet$-construction is applied.
\begin{definition}[\text{\rm{Following \cite[Definition~5.1]{CLMPZ}}}]\label{def:dgWald}
A \emph{dg Waldhausen category} consists of a dg category $\underline{\mathcal P}$ together with a Waldhausen structure on the associated category $\mathcal P \coloneqq Z^0\underline{\mathcal P}$ of closed degree-zero morphisms, satisfying the following conditions:
\begin{enumerate}[noitemsep,leftmargin=*,nosep]
\item In $\underline{\mathcal P}$, for each weak equivalence $u:X\to Y$ and any object $E$, the following two maps induced by pre- and post-composition are quasi-isomorphisms of chain complexes
\begin{equation}
u^*:\underline{\mathcal P}(Y,E) \to \underline{\mathcal P}(X,E)\text{ and }
u_*:\underline{\mathcal P}(E,X) \to \underline{\mathcal P}(E,Y);\label{formula:dgWaldweakequivalence}
\end{equation}
\item In $\underline{\mathcal P}$, for the pushout $W$ of $X\to Z$ along a cofibration $X\hookrightarrow Y$ and an object $E$, the induced squares 
\begin{align}\vcenter{\xymatrix@=1.5em{\underline{\mathcal P}(W,E) \ar[r] \ar[d] & \underline{\mathcal P}(Y,E)\ar[d]\\
\underline{\mathcal P}(Z,E)\ar[r] & \underline{\mathcal P}(X,E)}}\quad\vcenter{ \xymatrix@=1.5em{\underline{\mathcal P}(E,X) \ar[r] \ar[d] & \underline{\mathcal P}(E,Y)\ar[d]\\
\underline{\mathcal P}(E,Z)\ar[r] & \underline{\mathcal P}(E,W)}}\label{formula:dgWaldbicartesian}\end{align}
are homotopy bicartesian in $\mathrm{Ch}(k)$, the category of unbounded cochain complexes of $k$-modules.
\end{enumerate}
\end{definition}

Because $\mathcal P = Z^0\underline{\mathcal P}$, the zero object of $\mathcal P$ satisfies
$\underline{\mathcal P}(0,X) = 0 = \underline{\mathcal P}(X,0)\text{ for any }X.$
This is the condition \cite[Definition~5.1(1)]{CLMPZ}.

\begin{lemma}\label{lemma:waldcategories}
The dg category $\underline{\mathcal M}$ is a dg Waldhausen category. 
\end{lemma}

With the ambient model structures in Section~\ref{section:algebraicKfordgmodulesandcdgcomodules}, model-category homotopy between cofibrant-fibrant objects agrees with dg homotopy. 
Hence for a morphism $u$ in $\mathcal M$
\begin{equation*}
\begin{array}{ccl}
 \text{$u$ is a weak equivalence in $\mathrm w\mathcal M$} & \Leftrightarrow & \text{$u$ is invertible in $H^0\underline{\mathcal M}$}\\
 & \Leftrightarrow & \text{$u^*$ in (\ref{formula:dgWaldweakequivalence}) is a quasi-isomorphism for each $E$}\\
 & \Leftrightarrow & \text{$u_*$ in (\ref{formula:dgWaldweakequivalence}) is a quasi-isomorphism for each $E$}.
\end{array}
\end{equation*}
Indeed, either family of quasi-isomorphisms induces bijections on the corresponding \(H^0\)-Hom sets, so Yoneda implies that $u$ in $H^0\underline{\mathcal M}$ is invertible. The signs and the graded-splitting argument need no further changes.

\begin{proof}
We check (\ref{formula:dgWaldweakequivalence}).
Let $\underline{d}$ denote differentials of mapping complexes.  
As $\underline{\mathcal M}$ has objects being cofibrant-fibrant, a weak equivalence $u:X\to Y$ has a homotopy inverse $v$ and degree $(-1)$ endomorphisms $H_X$ and $H_Y$ so that 
$vu - \id_X = \underline{d}(H_X) = d_XH_{X} + H_Xd_X$ and $uv - \id_Y = d(H_Y).$
Then for any $f$ in $\underline{\mathcal M}(X,E)$ of homogeneous degree $p$, consider the degree $-1$ operator $K_X: f\mapsto (-1)^pfH_X$. 
The Leibniz rule writes
\[\underline{d}(K_Xf) = (-1)^p\underline{d}(f)H_X +  f \underline{d}(H_X).\]
Therefore, we check the desired equation
\begin{align*}(\underline{d}K_X + K_X\underline{d})(f) & = (-1)^p \underline{d}(f)H_X + f\underline{d}(H_X) - (-1)^p\underline{d}(f)H_X\\
& = f\underline{d}(H_X) = f(vu - \id_X) = fvu - f.
\end{align*}
Similarly, by putting $K_Y(f)= (-1)^{|f|}fH_Y,$ $L_X(g)=H_Xg$ on $\underline{\mathcal M}(E,X)$ and $L_Y(g) = H_Yg$ on $\underline{\mathcal M}(E,Y)$, we have $dK_Y + K_Yd = v^*u^*-\id,$ $dL_X + L_Xd = v_*u_* - \id$, and $dL_Y + L_Yd = u_*v_* - \id$. Hence $u^*$, $v^*$, $u_*$, and $v_*$ are chain homotopy inverses, and (\ref{formula:dgWaldweakequivalence}) follows.

Next, we check diagrams in (\ref{formula:dgWaldbicartesian}) are bicartesian. 
Set $i:X\hookrightarrow Y$ and $f:X\to Z$.
We have the following degreewise split exact sequence
\begin{align*}0\to X\xrightarrow{\makebox[1cm]{\scriptsize \text{$(i,-f$)}}} Y\oplus Z \to W\to 0.\end{align*}
For the degreewise splitting, if $\underline{\mathcal M}$ consists of dg modules, then the graded module map underlying $i:X\hookrightarrow Y$ is a split monomorphism.
If $\underline{\mathcal M}$ consists of cdg comodules, then as $X$ is fibrant, i.e., the underlying graded comodule $X^\#$ is injective, the graded comodule map underlying $i:X\hookrightarrow Y$ is a split monomorphism.
The degreewise split exactness of the above sequence implies the following two degreewise split exact sequences
\begin{align*} & 0\to \underline{\mathcal M}(W,E) \to \underline{\mathcal M}(Y,E)\oplus \underline{\mathcal M}(Z,E) \xrightarrow{\makebox[1cm]{\scriptsize \text{$i^*-f^*$}}} \underline{\mathcal M}(X,E) \to 0\text{ and}\\& 0\to \underline{\mathcal M}(E,X) \xrightarrow{\makebox[1cm]{\scriptsize \text{$(i_*,-f_*)$}}}\underline{\mathcal M}(E,Y)\oplus \underline{\mathcal M}(E,Z) \to \underline{\mathcal M}(E,W) \to 0\end{align*}
A short exact sequence of cochain complexes identifies its last term with the homotopy cofiber of its first map, and hence the squares in (\ref{formula:dgWaldbicartesian}) are homotopy bicartesian.
\end{proof}

Let $\mathrm{Ar}[n]$ be the category of pairs $(i, j)$ with $0 \leq i \leq j \leq n$ and morphisms are the obvious order-preserving maps.
For a Waldhausen category $\mathcal P$, recall that Waldhausen's $S_\bullet$-construction $S_\bullet \mathcal P$ consists for all $n\geq 0$ the categories $S_n\mathcal P$ whose objects are functors $\mathrm{Ar}[n]\to \mathcal P$ satisfying certain pushout conditions and morphisms are natural transformations.

Alternatively, put $\mathcal I_n=[n]\times[n]$. 
Extending $X$ by zero for $i>j$ identifies $S_n\mathcal P$ with the full dg subcategory of $\operatorname{Fun}(\mathcal I_n,\mathcal P)$ whose objects satisfy $X_{ij}=0$ for $i\geq j$ and the preceding cofibration and pushout conditions. The added structure maps are zero; 
compatibility with composition follows because any original structure map factoring through an added zero entry already factors through a diagonal zero entry. 
We use this description below, with the Reedy structure of Example~\ref{example:Sconstruction}, to introduce the dg $S_\bullet$-construction following \cite{CLMPZ}.

\begin{definition}\label{def:dgSconstruction}
Let $\underline{\mathcal P}$ be a dg Waldhausen category which is a full subcategory of a $\mathrm{Ch}(k)$-model category $\underline{\mathcal Q}$, and $\mathcal P \coloneqq Z^0\underline{\mathcal P}$, $\mathcal Q \coloneqq Z^0\underline{\mathcal Q}.$
Assume that the inclusion functor $\mathcal P\hookrightarrow\mathcal Q$ preserves the zero object, cofibrations, weak equivalences, and pushouts along cofibrations.
Use the size conventions fixed in Section~\ref{section:algebraicKfordgmodulesandcdgcomodules}.
Suppose that all objects in $S_{n}\mathcal{P}$ are both cofibrant and fibrant in the sense of the Reedy model structure on $\mathcal Q^{\mathcal I_n}$.
For each $n$, its $S_n$-\emph{construction} is the subcategory $\underline{S}_n\underline{\mathcal P}$ of $\underline{\mathcal Q}^{\mathcal I_n}$ whose objects are those of $S_{n} \mathcal{P}$, and whose mapping complex between any two objects $X$ and $Y$ is the strict end, where the equalizer is taken in $\mathrm{Ch}(k)$ and its two maps send $(a_r)_r$ to $(Y(\alpha)a_r)_\alpha$ and
$(a_sX(\alpha))_\alpha$, respectively
\[\underline{\mathcal Q}^{\mathcal I_n}(X,Y) = \mathrm{Eq}\left(\prod_{r\in \mathcal I_n} \underline{\mathcal Q}(X_r,Y_r) \rightrightarrows \prod_{(r\to s)\in \mathcal I_n} \underline{\mathcal Q}(X_r,Y_s)\right).\]
\end{definition}

By Theorem~\ref{thm:ChkReedy} and Lemma~\ref{lemma:cofib-fibreplacementcomputesderivedhom}, the Reedy cofibrancy-fibrancy requirement in the definition implies the strict end in $\underline{S}_{\bullet}\underline{\mathcal P}$ represents a \emph{derived} mapping complex; its zeroth cohomology gives the corresponding homotopy-category hom-set. 

For an order-preserving map $\theta:[m]\to[n]$, precomposition with
$\theta\times\theta$ defines a dg functor
\[\theta^*: \underline{S}_n\underline{\mathcal P}\longrightarrow\underline{S}_m\underline{\mathcal P},\quad(\theta^*X)_{ij}=X_{\theta(i),\theta(j)},\quad(\theta^*a)_{ij}=a_{\theta(i),\theta(j)}.\]
The $S$-conditions and the naturality equations for homogeneous families are preserved under reindexing. 
$\theta^*$ is a dg functor.
Indeed, for such a family $a$ in the mapping complex, and a composable family $b$, one has \[d(\theta^*a)=\theta^*(da),\quad\theta^*(ba)=(\theta^*b)(\theta^*a),\quad\theta^*(\id_X)=\id_{\theta^*X}.\]
For composable order-preserving maps $\eta:[l]\to[m]$ and $\theta:[m]\to[n]$, reindexing also gives $(\theta\circ\eta)^*=\eta^*\theta^*$ and $(\id_{[n]})^*=\id$.
Therefore, $\underline{S}_\bullet\underline{\mathcal P}$ is a simplicial dg category with faces and degeneracies respectively given by cofaces and codegeneracies in $\Delta$.

For an $S_n$-cofibration $f:X\to Y$ \cite[Section~2]{BM2}, the cofibrations $\textstyle X_{ij}\coprod_{X_{i,j-1}}Y_{i,j-1}\to Y_{ij}$ in $\mathcal Q$ for $i<j$ are the relative latching maps. 
Thus every $S_n$-cofibration is a Reedy cofibration in $\mathcal Q^{\mathcal I_n}$ as remaining relative latching maps are identities of zeros. 
Weak equivalences are objectwise. 
Pushouts along $S_n$-cofibrations are computed objectwise and, since the inclusion $\mathcal P\hookrightarrow\mathcal Q$ preserves pushouts along cofibrations, are preserved by the inclusion into $\mathcal Q^{\mathcal I_n}$.
\begin{lemma}
Resume the notations in the definition.
Assume $\mathcal Q$ is stable.
Then $\underline{S}_\bullet\underline{\mathcal P}$ is a simplicial object in dg Waldhausen categories.
\end{lemma}

\begin{proof}
By Appendix~\ref{sec:Reedymodelcategories}, $\underline{\mathcal Q}^{\mathcal I_n}$ is a stable $\mathrm{Ch}(k)$-model category.
For $E\in S_n\mathcal P,$ the following functors are right Quillen whose left adjoint, by the fibrancy and cofibrancy of $E$, are respectively $K\mapsto E^K$ (regarded as an object of $(\mathcal Q^{\mathcal I_n})^{op}$) and $K\mapsto K\otimes_k E$
\begin{equation*}
\begin{array}{c}
h_E:(\mathcal Q^{\mathcal I_n})^{op}\to \mathrm{Ch}(k);\,\,M\mapsto \underline{\mathcal Q}^{\mathcal I_n}(M,E),\quad h^E:(\mathcal Q^{\mathcal I_n})\to \mathrm{Ch}(k);\,\,M\mapsto \underline{\mathcal Q}^{\mathcal I_n}(E,M).
\end{array}
\end{equation*}
Brown's lemma shows that $h_E$ preserves weak equivalences between cofibrant objects of $\mathcal Q^{\mathcal I_n}$ (fibrant objects of $(\mathcal Q^{\mathcal I_n})^{\mathrm{op}}$), and $h^E$ preserves weak equivalences between fibrant objects of $\mathcal Q^{\mathcal I_n}$.
All objects in $S_n\mathcal P$ satisfy both conditions. Therefore the precomposition and postcomposition maps in (\ref{formula:dgWaldweakequivalence}) are quasi-isomorphisms.

For (\ref{formula:dgWaldbicartesian}), take a pushout square $Z\leftarrow X\hookrightarrow Y$ in $S_n\mathcal P$ with bottom-right vertex $W$. 
By the compatibility of the embedding functor, it is a pushout square in $\mathcal Q^{\mathcal I_n}$ along a Reedy cofibration between cofibrant diagrams. 
Hence it is a homotopy pushout and also a homotopy pullback by the stability of $\mathcal Q^{\mathcal I_n}$.
The right derived functor of $h_E$ sends this homotopy pushout to a homotopy pullback. 
The right derived functor of $h^E$ preserves the homotopy pullback in $\mathcal Q^{\mathcal I_n}$.
The ordinary mapping complexes already compute these derived values and hence both induced mapping-complex squares are homotopy pullbacks. 
Since $\mathrm{Ch}(k)$ is stable, both are homotopy bicartesian. The zero diagram is a dg zero object.
Finally, as $n$ varies, because the operators $\theta^*$ are dg and exact, this lemma follows. 
\end{proof}

\begin{remark}
In the notation of the above lemma, put $W=Y\amalg_X Z$, with $X\hookrightarrow Y$ a cofibration. 
The contravariant mapping-complex square is the strict pullback
\[\underline{\mathcal Q}^{\mathcal I_n}(W,E)\cong \underline{\mathcal Q}^{\mathcal I_n}(Y,E)\times_{\underline{\mathcal Q}^{\mathcal I_n}(X,E)}\underline{\mathcal Q}^{\mathcal I_n}(Z,E).\]
Since $E$ is Reedy fibrant, SM7 makes
$\underline{\mathcal Q}^{\mathcal I_n}(Y,E)\to\underline{\mathcal Q}^{\mathcal I_n}(X,E)$ a fibration. Thus this is also a homotopy
pullback in $\mathrm{Ch}(k)$.
\end{remark}

The following lemma gives examples

\begin{lemma}\label{lemma:cofibrancyfibrancyofSconstruction}
For the dg Waldhausen category $\underline{\mathcal M}$ with its ambient $\mathrm{Ch}(k)$-model category $\underline{\mathcal Q}\coloneqq \underline{\mathrm M}_A$ or $\underline{\mathcal Q}\coloneqq \underline{\mathrm{CM}}_C$, the objects in the ordinary $S_\bullet$-construction $S_\bullet\mathcal M$, extended by zero to $\mathcal I_n$, are both cofibrant and fibrant in the Reedy model category $\mathcal Q^{\mathcal I_n}$.
\end{lemma}

\begin{proof}
By Example~\ref{example:Sconstruction}, the desired Reedy cofibrancy follows from $X_{i,j-1}\to X_{ij}$ and $X_{i0} = 0.$
For each $i<j$, the $S_n$-construction give the short exact sequence $X_{i,i+1}\hookrightarrow X_{ij}\twoheadrightarrow X_{i+1,j}$ in $\mathcal Q$.
The map $X_{ij}\twoheadrightarrow X_{(i+1)j}$ is a fibration in the dg module case as it is surjective.
It is also a fibration in the cdg comodule case because the kernel $X_{i,i+1}$ is injective as a graded comodule. 
These are fibration criterion in Propositions~\ref{prop:modelstructureforcdgcomodules} and \ref{prop:projectivemodelstructureformodules}.
\end{proof} 

\begin{remark}
One may ask whether the Reedy-fibrancy hypothesis in Definition~\ref{def:dgSconstruction} can be weakened by using a dg model of homotopy-coherent $S_n$-diagrams.
Such a model would need compatible composition and simplicial operators,
together with a comparison recovering the Waldhausen $K$-theory space.
We do not pursue this extension here.
\end{remark}

\section{The key construction for the Dennis trace map}\label{sec:keyconstruction}
The dg version of McCarthy's $\tau$ map (\ref{formula:originalDennistrace}) is the key to giving a dg version of the Dennis trace map. 
The dg version of the cyclic nerve is well-developed. 
We come up with a dg version of the source of $\tau$ in Section~\ref{sec:thickenednerve} without changing the $K$-theory space and then dg-lize $\tau$ in Section~\ref{sec:tau}.
We show in Lemma~\ref{lemma:tauisachainmap} that our dg-lized $\tau$ is compatible with the differentials arising on both sides, and in Proposition~\ref{prop:taumonoidal} that $\tau$ is moreover monoidal.

\subsection{The thickened nerve and its geometric realization}\label{sec:thickenednerve}
We retain the notations in previous sections.
For the category $\underline{\mathcal M}$ in Section~\ref{sec:dgK-theory}, let $\underline{\mathrm w\mathcal M}$ be the maximal Kan complex inside $\underline{N}_\bullet(\underline{\mathcal M})$, the \emph{dg nerve} as an $\infty$-category.
For the simplicial degree $n$, the simplicial set $\underline{\mathrm w\mathcal M}$ has $n$-simplices to be maps $\sigma:\Delta^n \to \underline{N}_\bullet(\underline{\mathcal M})$ whose all edges are \emph{weak equivalences}.
Alternatively, the category $\underline{\mathrm w\mathcal M}$ consists of all objects of $\mathcal M$, $1$-simplices are weak equivalences, all homotopies between weak equivalences, and all related higher homotopies and compatibility data.

Next, let $\mathcal E^n$ denote the \emph{indiscrete groupoid} with objects $\{0,1,...,n\}$ and exactly one morphism for any pairing $(i,j)$ for $i,j=0,1,\ldots,n$ such that the map $i\to i$ is identity and $(i\to j\to k) = (i\to k)$ for each $i,j,k$.
Put $E^n = N_\bullet \mathcal E^n$, which is contractible.
As $n$ varies, $E^n$ form a cosimplicial simplicial set $E^\bullet:\Delta\to \mathbf{sSet}$ that is cocyclic, where $\mathbf{sSet}$ is the category of simplicial sets.
\begin{definition}\label{def:thicknerve}
    Let $\underline{\mathcal K}$ be a Kan simplicial set.
    The \emph{thickened nerve} is the simplicial set whose degree $n$ is $\Bbb{N}_n\underline{\mathcal K} = \mathbf{sSet}(E^n,\underline{\mathcal K})$ with $\theta^*F = F\circ E(\theta)$ for $\theta:[m]\to [n].$    
\end{definition}

We use the explicit description of the dg nerve in
\cite[Construction~1.3.1.6]{LurieHA}, with the cohomological grading and signs
of \cite[equation~(2.2.9)]{Faonte}.
For $\underline{\mathcal K} = \underline{\mathrm w\mathcal M},$ a thickened $n$-simplex $F:E^n\to \underline{\mathrm w\mathcal M}$ is described by objects $X_i \coloneqq F(i)$ for all $i$ and structure maps for all unrestricted integer $m\geq 1$
\[f_{i_0\cdots i_m}\in\underline{\mathcal M}(X_{i_0},X_{i_m})^{\,1-m}\text{ with }i_0,\ldots,i_m\in\{0,\ldots,n\}.\]
Here $f_{i_0\cdots i_m}$ denotes the highest component of the simplex $F(i_0,\ldots,i_m)$. 
The indices $i_0,\ldots,i_m$ need not be increasing or distinct. 
The maps $f_{ij}$ are closed degree-zero weak equivalences.
The degeneracy conditions are $f_{ii}=\mathrm{id}_{X_i}$ and $f_{i_0\cdots i_m}=0$ for $m\geq2$ whenever two adjacent indices coincide.
Nonadjacent repetitions need not give zero.
For $m\geq 2$, the higher components satisfy
    \begin{align}\textstyle df_{i_0\cdots i_m} = \sum_{j=1}^{m-1}(-1)^{j-1}f_{i_0\cdots \widehat{i_j}\cdots i_m}+\sum_{j=1}^{m-1}(-1)^{1+m(j-1)}f_{i_{j}\cdots i_m}\circ f_{i_{0}\cdots i_{j}},\label{formula:formulaofthickenednerve}\end{align}
where $d$ is the differential of the mapping complex $\underline{\mathcal M}(X_{i_0},X_{i_{m}}).$

\begin{example}
The coherence relation (\ref{formula:formulaofthickenednerve}) for $f_{i_0\cdots i_m}$ can be visualized as the boundary formula for an $(m-1)$-cube, whose faces correspond to the omission and composition terms in that relation. 
For $m=4,$ this gives the following cube for $f_{ijklm}$.
\[\begin{tikzpicture}[
    every node/.style={
        inner sep=2pt,
        font=\large
    },
    arrow/.style={
        -{Computer Modern Rightarrow[
            length=4.2pt,
            width=2.6pt
        ]},
        line width=0.4pt,
        shorten <=1pt,
        shorten >=1pt
    },
    crossing arrow/.style={
        preaction={
            draw=white,
            line width=2.5pt,
            -
        },
        arrow
    }
]

\node (A) at (0,2)
    {$f_{lm}f_{kl}f_{jk}f_{ij}$};

\node (B) at (4.2,2)
    {$f_{lm}f_{kl}f_{ik}$};

\node (C) at (2.1,1.2)
    {$f_{lm}f_{jl}f_{ij}$};

\node (D) at (6.3,1.2)
    {$f_{lm}f_{il}$};

\node (E) at (0,0.5)
    {$f_{km}f_{jk}f_{ij}$};

\node (F) at (4.2,0.5)
    {$f_{km}f_{ik}$};

\node (G) at (2.1,-.3)
    {$f_{jm}f_{ij}$};

\node (H) at (6.3,-.3)
    {$f_{im}$};

\node[anchor=east] at (-1,2.5)
    {\textcolor{red}{$f_{klm}f_{ijk}$}};

\node at (3.1,2.5)
    {\textcolor{blue}{$f_{lm}f_{ijkl}$}};

\node[anchor=east] at (-1,1)
    {\textcolor{green}{$f_{jklm}f_{ij}$}};

\node[anchor=west] at (7.3,1)
    {\textcolor{green}{$f_{iklm}$}};

\node at (3.1,-.8)
    {\textcolor{blue}{$f_{ijkm}$}};

\node[anchor=west] at (7.3,-.8)
    {\textcolor{red}{$f_{ijlm}$}};

\draw[arrow] (A) -- (B);
\draw[arrow] (E) -- (F);
\draw[arrow] (G) -- (H);

\draw[arrow] (A) -- (C);
\draw[arrow] (B) -- (D);
\draw[arrow] (E) -- (G);
\draw[arrow] (F) -- (H);

\draw[arrow] (A) -- (E);
\draw[arrow] (B) -- (F);

\draw[crossing arrow] (C) -- (G);
\draw[arrow] (D) -- (H);
\draw[crossing arrow] (C) -- (D);
\end{tikzpicture}\]
The vertices are composites of degree-zero maps, the edges represent degree-$(-1)$ homotopies, and the square faces carry the displayed degree-$(-2)$ maps. The red, green, and blue labels indicate the front/back, left/right, and top/bottom faces, respectively. Their signed sum is $df_{ijklm},$ with signs as in (\ref{formula:formulaofthickenednerve}).
\end{example}

\begin{lemma}\label{remark:cyclicsettrivialcase}
The thickened nerve $\Bbb{N}_\bullet \underline{\mathcal K}$ is a cyclic set.
Its faces and degeneracies are operators defined by precomposition: $d_jF\coloneqq F\circ E(\delta^j)$ and $s_jF \coloneqq F\circ E(\sigma^j)$ with the cofaces $\delta^j$ and codegeneracies $\sigma^j$ for $j = 0,\ldots,n$ in $\Delta$, 
With $\rho_n:[n]\to[n];\,\,\rho_n(j) = j+1 \mod n+1,$ this map induces a functor between the indiscrete groupoids and hence the map $E(\rho_n)$ on $E^n$, the cyclic operator is $t_nF\coloneqq F\circ E(\rho_n)$.
Besides the simplicial identities, the operators satisfy
\[t_n^{n+1} = \id,\quad d_jt_n = \begin{cases}t_{n-1}d_{j+1}& j<n;\\d_0\quad & j = n,
\end{cases}\quad  s_jt_n = \begin{cases}t_{n+1}s_{j+1} & j<n;\\
t^2_{n+1}s_0 & j=n.\end{cases}\]
\end{lemma}
The proof is given by checking the corresponding identities for $\delta^j$, $\sigma^j,$ and $\rho_n$, and then precomposing. 
The cyclic operator $t_2$ maps, for instance, $f_{12}$ to $f_{20}$ and $f_{120}$ to $f_{201}$.
The reduced simplicial $k$-module $\widetilde{k}[\Bbb{N}_\bullet\underline{\mathcal K}]$ has the simplicial boundary \begin{align}\textstyle
\partial F = \sum_{j=0}^n (-1)^j d_jF\text{ for degree }n\geq 1\text{ and }\partial = 0\text{ in degree $0$}.\label{formula:differentialforthickennerve}\end{align}

\begin{lemma}\label{lemma:homotopyequivalenceforthicknerve}
Let $\underline{\mathcal K}$ be a Kan complex.
Restriction along $\Delta^n\hookrightarrow E^n$ defines a natural trivial Kan fibration $r:\Bbb{N}_\bullet \underline{\mathcal K}\to \underline{\mathcal K}.$
The geometric realization $|\Bbb{N}_\bullet\underline{\mathcal K}|$ is homotopy equivalent to the geometric realization $|\underline{\mathcal K}|.$
\end{lemma}

\begin{proof}
The natural embedding $[n]\hookrightarrow \mathcal E^n$ induces a map $\Delta^n = N_\bullet [n] \to E^n.$
We show that the map $\Delta^\bullet \to E^\bullet$ is a trivial cofibration in the sense of Reedy model structure for cosimplicial simplicial sets.
As $\Delta^n$ and $E^n$ are both contractible, the embedding is a degreewise weak equivalence.
For the cofibration, we check the relative latching map $j_n:\Delta^n \coprod_{L^n \Delta^\bullet}L^nE^\bullet \to E^n$ is cofibration in the sense of Kan-Quillen, that is, in the present, injective as a map of sets.
The latching object $L^n \Delta^\bullet$ consists of all non-decreasing sequence $(i_0,\ldots,i_s)$ with $i_j\in\{0,...,n\}$ so that the associated map $[s]\to [n]$ is not surjective.
The latching object $L^n E^\bullet$ consists of all sequence $(i_0,\ldots,i_s)$ with $i_j\in\{0,...,n\}$ so that the associated map $[s]\to [n]$ is not surjective.
Hence the desired injectivity follows.

By \cite[Lemma~1.4]{Reedy}, every relative latching map $j_n$ is a trivial cofibration. 
A lifting problem against $\partial\Delta^n\hookrightarrow\Delta^n$ for the restriction map $\Bbb{N}_\bullet \underline{\mathcal K}\to \underline{\mathcal K}$ is precisely an extension problem $\Delta^n\amalg_{\partial\Delta^n}L^nE^\bullet \to \underline{\mathcal K}$ across $j_n$. 
This extension exists as $\underline{\mathcal K}$ is Kan. 
Hence $r$ is a trivial Kan fibration. 
Its geometric realization is a homotopy equivalence.
\end{proof}

Since $\mathcal M$ consists of cofibrant-fibrant objects, a morphism in $\mathcal M$ is a weak equivalence if and only if it admits a dg homotopy inverse (Lemma~\ref{lemma:waldcategories}). 
Thus the weak equivalences are precisely the edges of $\underline{\mathrm w\mathcal M}$, and the thickened nerve records them together with coherent inverse and higher homotopy data. 
Later, we apply $\Bbb{N}_\bullet$ degreewise to the maximal Kan complex $\underline{\mathrm wS}_\bullet\underline{\mathcal M}$ of the dg nerve in each $S_\bullet$-degree.

\subsubsection{A geometric realization using $|E^n|$}
Let $\underline{\mathcal K}$ be a Kan complex. It is natural to consider a geometric realization of the thickened nerve $\Bbb{N}_\bullet\underline{\mathcal K}$ where $|E^n|$ plays the role of $|\Delta^n|$. 
It turns out that this new geometric realization is homotopy equivalent to the ordinary geometric realization $|\Bbb{N}_\bullet\underline{\mathcal K}|$.
It allows an explicit description of the cyclic operator but will not be needed for the rest of this paper.

\begin{definition}
Let $\Bbb{N}_\bullet\underline{\mathcal K}$ be a thickened nerve.
Its $E$-\emph{geometric realization} writes 
    \[|\Bbb{N}_\bullet\underline{\mathcal K}|_E\coloneqq \left(\coprod_{n\geq 0}\Bbb{N}_n\underline{\mathcal K}\times |E^n|\right)/\sim\]
    where $|E^n|$ denotes the geometric realization of the simplicial set $E^n$ and where $\sim$ is generated by $(\theta^*\sigma,t) \sim (\sigma,|E(\theta)|(t))$ for every morphism $\theta:[m]\to[n]$ in $\Delta$, every $\sigma\in \Bbb{N}_n\underline{\mathcal K}$, and every $t\in|E^m|$. 
    Here $\theta^*: \Bbb{N}_n\underline{\mathcal K}\to\Bbb{N}_m\underline{\mathcal K}$ and $E(\theta)^m: E^m\to E^n$ are respectively the simplicial and cosimplicial structure maps.
\end{definition}

\begin{lemma}
Both realizations of $\Bbb{N}_\bullet \underline{\mathcal K}$ carry natural $S^1$-actions, and there is a natural $S^1$-equivariant map $R_{\underline{\mathcal K}}:|\Bbb{N}_\bullet |_E\rightarrow |\Bbb{N}_\bullet \underline{\mathcal K}|$
whose underlying map is a homotopy equivalence.
\end{lemma}

\begin{proof}
Note that $\Delta^n\hookrightarrow E^n$ induces, for every simplicial set $X$, a natural map $a_X:|X|\rightarrow |X|_E.$
This is a weak homotopy equivalence by induction on the skeleta of $X$: the attaching inclusions are $\partial\Delta^n\hookrightarrow\Delta^n$ on one side and $L^nE^{\bullet}\hookrightarrow E^n$ on the other; both are cofibrations, and the simplexwise comparison is a weak equivalence. 
Let $D_n(X)$ be the set of nondegenerate $n$-simplices of $X$ whose elements $\sigma$ satisfies $\sigma|_{\partial \Delta^n}:\partial \Delta^n\to \mathrm{sk}_{n-1}\,X.$
The gluing lemma for the following diagram and passage to the union of skeleta give the assertion
\[\xymatrix@=1.5em{
\textstyle \prod_{\sigma\in D_n(X)}|\Delta^n|\ar[d]_{\simeq} & \prod_{\sigma\in D_n(X)}|\partial \Delta^n|\ar@{_{(}->}[l]\ar[d]_{\simeq}\ar[r] & |\mathrm{sk}_{n-1}\,X|\ar[d]_{\simeq}\\
\prod_{\sigma\in D_n(X)}|E^n| & \prod_{\sigma\in D_n(X)}|L^n E^\bullet|\ar@{_{(}->}[l]\ar[r] & |\mathrm{sk}_{n-1}\,X|_E.
}\]
Both spaces are realizations of simplicial sets (take the defining quotient before realization), hence are CW complexes.
Thus $a_X$ is a homotopy equivalence.

Take $X=\Bbb{N}_\bullet \underline{\mathcal K}$. 
It remains to construct the equivariant comparison $R_{\underline{\mathcal K}}$ as a homotopy inverse of $a_{X}.$
For fixed $n$, each element $i$ in the simplicial set $E^n$ may be written as $i:[q]\to[n];\,\,j\mapsto i_j$ (need not be order-preserving) and the cyclic structure on $E^n$ in the simplex dimension $q$ explicitly writes $(i_0,\ldots,i_q)\mapsto(i_1,\ldots,i_q,i_0).$
Moreover, the canonical circle actions on $|E^n|$ descend through the quotient relations to an action on $|X|_E$.

Given $F:E^n\to \underline{\mathcal K}$, define a map of cyclic sets 
$\phi_F:E^n\rightarrow X,\,\,
(i_0,\ldots,i_q)\mapsto F\circ E(i).$
For $\theta:[m]\to[n]$ in $\Delta$, we also have $\phi_{\theta^*F}=\phi_F\circ E(\theta).$
Compatibility with faces and degeneracies follows from composition
of functions. 
Compatibility with cyclic rotation follows from $\phi_F(i\circ\rho_q)=F\circ E(i)\circ E(\rho_q)=t_q\phi_F(i).$
The realized maps $|\phi_F|$ therefore respect the relations defining
$|X|_E$ and give an $S^1$-equivariant map
\[
R_{\underline{\mathcal K}}:|X|_E\longrightarrow |X|,
\quad [F,z]\longmapsto |\phi_F|(z).
\]
More explicitly, if we write $z\in |E^n|$ as $[(a_0,\ldots,a_q),u]$ with
$u\in |\Delta^q|$, then $R_{\underline{\mathcal K}}([F,z])=[F\circ E(a),u].$
For the simplex $(0,\ldots,n)$, the function $a$ is the identity.
Consequently, 
$R_{\underline{\mathcal K}}\circ a_X=\id_{|X|}.$
Since $a_X$ is a homotopy equivalence, so is the underlying map of $R_{\underline{\mathcal K}}$.
\end{proof}

\subsubsection{The comparison}

By the following comparison proposition, one may use $\underline{\mathrm wS}_\bullet\underline{\mathcal M}$ to alternatively define the $K$-theory.
The rest of Section~\ref{sec:thickenednerve} is devoted to proving this comparison proposition by reducing it to Blumberg--Gepner--Tabuada's result \cite[Proposition~2.10]{BGT}, using the comparison between dg nerves and simplicial localizations.

\begin{proposition}[Comparison]\label{prop:comparison}
For the categories $\mathcal M$ and $\underline{\mathrm wS}_\bullet\underline{\mathcal M}$ as above, we have the following homotopy equivalence
\[|N_\bullet\mathrm{w}S_\bullet\mathcal M|\simeq |\Bbb{N}_\bullet\underline{\mathrm wS}_\bullet\underline{\mathcal M}|.\]
\end{proposition}

\noindent \text{\bf Step~0: Reduction to cofibrant-fibrant diagrams.}
By Lemma~\ref{lemma:homotopyequivalenceforthicknerve}, degreewise in the $S_\bullet$-direction and then by the diagonal lemma, we have $|\Bbb{N}_\bullet \underline{\mathrm wS}_\bullet\underline{\mathcal M} | \simeq |\underline{\mathrm wS}_\bullet\underline{\mathcal M}|.$
We need to compare the latter with $|N_\bullet\mathrm w S_\bullet\mathcal M|$.
For each $m$, the canonical map from the ordinary nerve to the dg nerve
restricts to
\[
 c_m:N_\bullet\mathrm w S_m\mathcal M\longrightarrow
 \bigl(\underline{N}_\bullet(\underline{S}_m\underline{\mathcal M})\bigr)^{\simeq}
 =\underline{\mathrm wS}_m\underline{\mathcal M}.
\]
Indeed, weak equivalences between the Reedy cofibrant-fibrant $S_m$-diagrams are dg homotopy equivalences. 
The maps $c_m$ commute with the $S_\bullet$-operators, since those operators are dg functors.
We will show that each $c_m$ is a weak homotopy equivalence.

Let $[\overline{m}]$ be the category $\{1\to \cdots\to m\}$ and $[\overline{0}] = \varnothing.$ 
Let $\mathcal Q$ denote $\mathrm M_A$ or $\mathrm CM_C$ and $\underline{\mathcal Q}$ denote the dg enrichment $\underline{\mathrm M}_A$ or $\underline{\mathrm{CM}}_C.$
Consider the auxiliary category \[\mathcal D_m\coloneqq \left\{X\in \mathcal Q^{[\overline{m}]}\,\left|\, \begin{aligned}&\text{the image of }X_i\text{ in the (co)derived}\\
& \text{category is in }\mathrm{thick}_C(M)\text{ or }\mathrm{thick}_A(N)\end{aligned}\right.\right\}.\]
This auxiliary category is motivated by \cite[Proposition~2.10]{BGT}.
We use throughout the size conventions of Section~\ref{section:algebraicKfordgmodulesandcdgcomodules}. 
Weak equivalences in $\mathcal D_m$ are objectwise. 
Let $\mathcal B$ denote the full subcategories of cofibrant-fibrant objects, in the ambient projective Reedy model structure on $\mathcal Q^{[\overline{m}]}$ with the Reedy structure $[\overline{m}]^+ = [\overline{m}]$ and the discrete subcategory $[\overline{m}]^-$. 
Let $\underline{\mathcal B}$ be the full dg subcategory of the dg diagram category on $\mathcal B$ with strict mapping complexes.

By the latching-object description, an object of $\mathcal B$ is a composition $0\hookrightarrow X_1\hookrightarrow\cdots\hookrightarrow X_m$ with $X_i\in\mathcal M.$
Restriction to $X_{01},\ldots,X_{0m}$ gives an equivalence
$S_m\mathcal M\simeq\mathcal B$, with inverse given by adjoining
$X_0=0$ and forming the quotients $X_j/X_i$.
These quotients belong to $\mathcal M$: their images in the derived or coderived category remain in the thick subcategory, and they are cofibrant-fibrant by Section~\ref{section:algebraicKfordgmodulesandcdgcomodules}.
Moreover, a homogeneous map between objects in $\mathcal B$ induces unique maps on the quotients and is compatible with differentials. 
Thus restriction identifies the mapping complexes and moreover gives a dg equivalence
$\underline{S}_m\underline{\mathcal M}\simeq\underline{\mathcal B}$.

We are to check the following weak homotopy equivalences for each $m$:
\begin{align}N_\bullet \mathrm{w}S_m\mathcal{M} \simeq (N_\bullet^{\mathrm{hc}}((L^H\mathcal B)^{\mathrm{fib}}))^\simeq\simeq \underline{\mathrm wS}_m\underline{\mathcal M},\label{formula:comparison}\end{align}
where $L^H\mathcal B = L^H(\mathcal B,\mathrm w\mathcal B)$ denotes the hammock simplicial localization; the upper index $\mathrm{fib}$ means taking the fibrant replacement in the Bergner model structure on simplicial categories; $N_\bullet^{\mathrm{hc}}$ denotes the homotopy coherent nerve; and $(-)^\simeq$ is the maximal Kan complex of the resulting $\infty$-category, whose simplices have edges being equivalences.

\noindent \text{\bf Step~1: The Blumberg-Gepner-Tabuada comparison.}
For the below equivalences, the left follows from \(S_m\mathcal M\simeq\mathcal B\). 
We now prove the right equivalence.
\begin{align}N_\bullet\mathrm{w}S_m\mathcal M\simeq N_\bullet \mathrm{w}\mathcal B \simeq (N_{\mathrm{hc}}((L^H\mathcal B)^{\mathrm{fib}}))^\simeq.\label{formula:comparison1}\end{align} 
For the right weak homotopy equivalence in (\ref{formula:comparison1}), we apply \cite[Proposition~2.10]{BGT} to $(\mathcal B,\mathrm w\mathcal B)$ with the size conventions of Section~\ref{section:algebraicKfordgmodulesandcdgcomodules}. 
Weak equivalences satisfy the two-out-of-three property by the ambient model structure. 
We claim $(\mathcal B,\mathrm w\mathcal B)$ admits a homotopy calculus of two-sided fractions by checking Dwyer-Kan's criterion \cite[Proposition~8.2]{DwyerKan}. 
The constructions below are functorial and verify the three conditions in loc.\ cit. 
\begin{itemize}[noitemsep,leftmargin=*,nosep]
    \item 
Let their $W_1$ and $W_2$ consist respectively of maps in $\mathcal B$ that are trivial cofibrations and trivial fibrations in $\mathcal Q^{[\overline{m}]}$, and choose functorial factorizations in the ambient diagram model category, using the constructions in the proofs of \cite[Theorems~8.1 and 8.2]{Posit} (also \cite[Theorem~5.1.3]{Hovey} for the diagram case). 
Let $Q$ and $R$ be the resulting functorial cofibrant and fibrant replacements.
    \item
Given $X'\xleftarrow{u}X\xrightarrow{f}Y$ in $\mathcal B$ with $u\in W_1$, let $P\coloneqq X'\amalg_XY$ be the ambient pushout and take a fibrant replacement $P\overset{\simeq}{\hookrightarrow} RP$ by a trivial cofibration. Since $Y\to P\to RP$ is a composite of trivial cofibrations and $Y$ is cofibrant, $RP$ is cofibrant-fibrant. 
Moreover, $RP$ is weakly equivalent to $Y$, so its entries remain in the prescribed thick subcategory. 
Thus $RP\in\mathcal B$, and the commutative square with bottom-right object $RP$ has its map $Y\to RP$ in $W_1$.
If $f$ is a weak equivalence, so is $X'\to RP$ by two-out-of-three.
    \item
Next, given dually  $X\xrightarrow{g}Y\xleftarrow{v}Y'$ with $v\in W_2$, put now $P\coloneqq X\times_Y Y'$ and take a cofibrant replacement $QP\overset{\simeq}{\twoheadrightarrow} P$ by a trivial fibration. 
Since $QP\to P\to X$ is a composite of trivial fibrations and $X$ is fibrant, $QP$ is cofibrant-fibrant.
As $QP$ is weakly equivalent to $X$, it belongs to $\mathcal B$.
The resulting square has $QP\to X$ in $W_2$. 
If $g$ is a weak equivalence, then $QP\to Y'$ is also one by two-out-of-three.
    \item 
Finally, factor a weak equivalence $X\to Y$ functorially as a cofibration $X\to Z$ followed by a trivial fibration $Z\to Y$.
The first map is also a weak equivalence by two-out-of-three.
Since $X$ is cofibrant and $Y$ is fibrant, $Z$ is cofibrant-fibrant and since moreover $Z$ is weakly equivalent to $X$, it belongs to $\mathcal B$. 
\end{itemize} 
Moreover, $\mathrm w\mathcal B$ is closed in $\mathcal B$ in the sense of Dwyer-Kan, that is, a morphism $f$ of $\mathcal B$ belongs to $\mathrm w\mathcal B$ if and only if its image in $\mathcal B[(\mathrm w\mathcal B)^{-1}]$ is an isomorphism. 
Indeed, the inclusion into the ambient diagram model category $\mathcal Q^{[\overline{m}]}$ induces a functor $\mathcal B[(\mathrm w\mathcal B)^{-1}] \longrightarrow \mathrm{Ho}(\mathcal Q^{[\overline{m}]}).$
If $f$ becomes invertible in the source, it becomes invertible in $\mathrm{Ho}(\mathcal Q^{[\overline{m}]})$, hence is an ambient weak equivalence by \cite[Theorem~1.2.10(iv)]{Hovey}. 
Thus $f\in \mathrm w\mathcal B$ and the converse holds by the definition of localization.
Consequently, \cite[Proposition~2.10]{BGT} shows that the canonical map in the right of (\ref{formula:comparison1}) is a weak homotopy equivalence, as desired.

\noindent \text{\bf Step~2: The dg nerve and hammock localization.}
Next, by \cite[Proposition~1.3.1.17]{LurieHA}, we have the equivalence of $\infty$-categories $\underline{N}_\bullet(\underline{S}_\bullet\underline{\mathcal M})\simeq  N_\bullet^\mathrm{hc}((\underline{S}_\bullet\underline{\mathcal M})_\Delta).$
For any dg category $\underline{\mathcal C}$, we have $\underline{\mathcal C}_\Delta(X,Y) = \Gamma(\underline{\mathcal C}(X,Y)^{\leq 0})$ with degree-zero term $Z^0\mathcal C(X,Y)$ and the functor $\Gamma$ given by the Dold-Kan correspondence; composition uses the Alexander-Whitney construction.
The outcome mapping spaces are Kan.
By the dg equivalence of $\infty$-categories $\underline{\mathcal B} \simeq \underline{S}_m\underline{\mathcal M}$, we have the equivalence  $N^\mathrm{hc}_\bullet((\underline{S}_m\underline{\mathcal M})_\Delta) \simeq N^\mathrm{hc}_\bullet (\mathcal B_\Delta).$
We are to check the following rightmost equivalence of $\infty$-categories \begin{align}\underline{N}_\bullet(\underline{S}_\bullet\underline{\mathcal M})\simeq  N_\bullet^\mathrm{hc}((\underline{S}_\bullet\underline{\mathcal M})_\Delta)\simeq N^\mathrm{hc}_\bullet (\mathcal B_\Delta) \simeq N_\bullet^\mathrm{hc}((L^H\mathcal B)^\mathrm{fib}).\label{formula:comparison2}\end{align}

We are to follow the argument in \cite[Section~1.4]{Hinich} to show the following three DK-equivalences of simplicial categories (that is a functor between simplicial categories which induces weak homotopy equivalences between all hom-simplicial sets and induces an essentially surjective functor between homotopy categories):
\begin{align}
L^H\mathcal B \overset{\simeq}{\longrightarrow} L^H(\mathcal B_\Delta)\overset{\simeq}{\longleftarrow}\mathcal B_\Delta\label{formula:hammockversionofthesimpliciallocalization}
\end{align}
Here: $L^H$ means the diagonal of the degreewise hammock localization; 
the left functor $L^HU$ is induced by the inclusion as constant simplices functor $U$; and the right functor is the localization functor $j$.

The left DK-equivalence in (\ref{formula:hammockversionofthesimpliciallocalization}): Let $\mathcal B_n$ be the ordinary category with the same object as in $\mathcal B$ and with $\mathcal B_n(X,Y) = \mathcal B_\Delta(X,Y)_n$ that is all morphism with simplicial degree $n$ (which is \emph{not} the dg degree). 
By \cite[1.4.4]{Hinich}, We have two adjoint functors $U_n:\mathcal B\rightleftarrows \mathcal{B}_n:\Lambda_n,$ $U_n\dashv \Lambda_n,$ defined as below: $U_n$ is identity on objects and sends an ordinary morphism to its constant $n$-simplex; 
The functor $\Lambda_n$ maps $Y$ to its cotensor $\Lambda_nY = Y^{C_\Delta},$ taken objectwise in the ambient diagram category (See Remark~\ref{remark:weakpathobj}), where $C_\Delta \coloneqq \mathrm{Norm}_\bullet(k[\Delta^n])$ denotes the normalized chains, regarded as a cochain complex by equipping $\mathrm{Norm}_q(k[\Delta^n])$ with degree $-q$. 
The cotensor $\Lambda_nY$ is still in $\mathcal B.$ Indeed, as $C_\Delta$ is bounded and finite free, we have $\Lambda_nY \cong  C_\Delta^\vee\,\otimes_k Y$, and complexes $C_\Delta,$ $C_\Delta^\vee$ are a cofibrant.
By the dg Reedy model structure, tensoring with $C_\Delta^\vee$ preserves cofibrant objects and cotensoring with $C_\Delta$ preserves fibrant objects. 
Moreover, the augmentation $C_\Delta\to k$ is a chain homotopy equivalence, so $\Lambda_nY$ is dg homotopy equivalent to $Y$. 

For any morphism $f:Y\to Z$ in $\mathcal B_n$, $\Lambda_n(f)$ is the unique morphism $\theta_{\Lambda_nY,Z}(f\circ \varepsilon_Y)$ such that $\varepsilon_Z\circ U_n(\Lambda_n(f)) = f\circ \varepsilon_Y$, 
Here, $\varepsilon_Y \coloneqq \theta_{\Lambda_nY,Y}^{-1}(\id_{\Lambda_nY})$ is the counit with $\theta_{X,Y}$ denoting the following explicit isomorphism 
\begin{align}\mathcal B_n(U_nX,Y) \cong \Gamma(\underline{\mathcal B}(X,Y)^{\leq 0})_n\cong \mathrm{Ch}(k)(N_\bullet(k[\Delta^n]),\underline{\mathcal B}(X,Y))\cong \mathcal B(X,\Lambda_nY).\label{formula:DKadjunction}\end{align}
The isomorphism is natural in $X$, and thus defines the right adjoint on morphisms by representability.
Let $v_i:\mathcal B_n \to \mathcal B$ be the functor that is the identity on the objects and taking the value on $i$-th vertex on the morphisms.
Let $\mathrm w\mathcal B_n$
be the subcategory with all objects and those morphisms $f$ for which
$v_i(f)$ is a weak equivalence for every $i$. Equivalently, it suffices
to check one vertex: the vertices of an $n$-simplex represent the same
class in $\pi_0\mathcal B_\Delta(X,Y)=H^0\underline{\mathcal B}(X,Y)$.
We are to apply \cite[Corollary~3.6]{DwyerKan} to show $L^HU_n:L^H\mathcal B\to L^H\mathcal B_n$ is a DK-equivalence.
For this we need to check the hypothesis:
\begin{itemize}[noitemsep,leftmargin=*,nosep]
    \item The unit $\eta:\id_{\mathcal B}\to \Lambda_nU_n$ on $X$ writes $\eta_X = a_n^*:X^k \to X^{C_\Delta} = \Lambda_nX$ which is induced by the augmentation $a_n:C_\Delta\to k.$
    By (\ref{formula:DKadjunction}), the counit is an $n$-simplex in $\mathcal B_n(U_n\Lambda_nY,Y)$ so that $v_i(\varepsilon_Y):\Lambda_nY\to Y$ is induced by the $i$-th vertex inclusion $s_i:k\to C_\Delta$.
    The contraction of $\Delta^n$ to its $i$-th vertex implies the maps of complexes: $a_n s_i = \id_k$ and $s_ia_n \simeq \id_{C_{\Delta}}$.
    They induce $v_i(\varepsilon_Y)\circ \eta_Y = \id_Y$ and $\eta_Y\circ v_i(\varepsilon_Y)\simeq \id_{\Lambda_nY}.$
    Hence $\eta_X$ and $\epsilon_Y$ are weak equivalences in respectively $\mathrm w\mathcal B$ and $\mathrm w\mathcal B_n.$
    \item The functor $U_n$ preserves weak equivalence by its definition.
    If $f\in \mathrm w\mathcal B_n$, apply $v_i$ to the equation defining $\Lambda_n$ and we have $v_i(\varepsilon_Z)\circ \Lambda_n(f) = v_i(f)\circ v_i(\varepsilon_Y).$
    We already knew that $v_i(\varepsilon_Z),$ $v_i(f)$, and $v_i(\varepsilon_Y)$ are all weak equivalences and hence desired homotopy equivalence follows.
\end{itemize}
Consequently, loc.\ cit. gives the desired DK-equivalence $L^HU_n.$
The functor $U_n$ commutes with all simplicial operators: indeed, as $U_n$ outputs a constant $n$-simplex on a morphism, for any order preserving $\theta:[n]\to [n']$ induced $\theta^*:\mathcal B_{n'}\to \mathcal B_n,$ we have $\theta^*\circ U_{n'} = U_n$ and hence $(L^H\theta^*)\circ (L^HU_{n'}) = L^H(\theta^*\circ U_{n'}) = L^HU_n$ as $L^H$ is functorial.
The categories $\mathrm w\mathcal B_n$ also assemble into a simplicial subcategory
$\mathrm w\mathcal B_\bullet\subset\mathcal B_\Delta$. For fixed objects $X,Y$, the functors $L^HU_n$ assemble to be a bisimplicial map 
\[A_{n,l} = L^H\mathcal B(X,Y)_l\quad\text{and}\quad B_{n,l}=L^H\mathcal B_n(X,Y)_l\]
such that $A_{n,l}$ is constant in the $n$-direction and $A_{n,\bullet}\to B_{n,\bullet}$ are weak homotopy equivalence for all $n$.
By the diagonal lemma, we have the weak homotopy equivalence \[L^H\mathcal B(X,Y) = \mathrm{diag}\,A_{\bullet,\bullet}\to \mathrm{diag}\,B_{\bullet,\bullet} = L^H(\mathcal B_\Delta,\mathrm w\mathcal B_\bullet)(X,Y).\]
Since \(L^HU\) is the identity on objects, it also induces an essentially surjective functor on homotopy categories.
All these imply the desired left DK-equivalence in (\ref{formula:hammockversionofthesimpliciallocalization}).

The rightmost DK-equivalence in (\ref{formula:hammockversionofthesimpliciallocalization}):
We retain the notations above.
Let $j:\mathcal B_\Delta \to L^H(\mathcal B_\Delta,\mathrm w\mathcal B_\bullet)$ be the canonical map.
The desired DK-equivalence follows from \cite[2.4(ii)]{BK}.
For its hypothesis, every morphism in $\mathrm w\mathcal B_\bullet$ represents
an invertible morphism in
$\pi_0\mathcal B_\Delta\cong H^0\underline{\mathcal B}$: weak equivalences between
these cofibrant-fibrant diagrams admit dg homotopy inverses.
Thus $\mathrm w\mathcal B_\bullet$ is neglectable in $\mathcal B_\Delta$ in the sense of loc.\ cit. and hence $j$ is a DK-equivalence.

\noindent \text{\bf Step 3: The canonical comparison and realization.}
Overall, for an integer $m\geq 0$, combining (\ref{formula:comparison1}) and the maximal Kan complexes of (\ref{formula:comparison2}), we have weak homotopy equivalence $N_\bullet \mathrm wS_m\mathcal M \simeq \underline{\mathrm w S}_m\underline{\mathcal M}.$
We need to relate this equivalence with the canonical map $c_m$.

Let $\ell:\mathcal B\to L^H\mathcal B$ be the localization map,
regarding $\mathcal B$ as a discrete simplicial category. Naturality
of localization gives the commutative square
\[\xymatrix@=1.5em{\mathcal B\ar[r]^U\ar[d]_\ell & \mathcal B_\Delta\ar[d]^j\\L^H\mathcal B \ar[r]^{L^HU\hspace{.7cm}}& L^H(\mathcal B_\Delta,\mathrm w\mathcal B_\bullet)}\]
The left comparison is (\ref{formula:comparison1}), the canonical map $N_\bullet\mathrm w\mathcal B\longrightarrow(N_\bullet^\mathrm{hc}((L^H \mathcal B)^{\mathrm{fib}}))^\simeq$ that is a weak homotopy equivalence. 
Since $U$ and $\ell$ send weak equivalences to equivalences, the square above induces, after functorial fibrant replacement at the localized corners, a commutative square of simplicial sets with upper-left corner $N_\bullet\mathrm w\mathcal B$ and the other three corners the corresponding maximal Kan complexes.
Its bottom and right maps are weak homotopy equivalences by the two DK-equivalences in (\ref{formula:hammockversionofthesimpliciallocalization}). 
Hence so is the map induced by $U$,
\[N_\bullet\mathrm w\mathcal B\longrightarrow      (N^\mathrm{hc}(\mathcal B_\Delta))^\simeq,\]
where $\mathcal B_\Delta$ already has
Kan mapping spaces.

Restriction $\underline{S}_m\underline{\mathcal M}\simeq \underline{\mathcal B}$ and the natural comparison of \cite[Proposition~1.3.1.17]{LurieHA} are compatible with the ordinary nerve maps. Hence the preceding conclusion says precisely that the canonical map
\[c_m:N_\bullet\mathrm w S_m\mathcal M\longrightarrow (\underline{N}_\bullet(\underline{S}_m\underline{\mathcal M}))^\simeq =\underline{\mathrm w S}_m\underline{\mathcal M}\]
is a weak homotopy equivalence. Since the maps $c_m$ commute with the $S_\bullet$-operators, the diagonal lemma gives the desired homotopy equivalence
\[|N_\bullet\mathrm w S_\bullet\mathcal M|\simeq |\underline{\mathrm w S}_\bullet\underline{\mathcal M}|\simeq |\Bbb{N}_\bullet\underline{\mathrm w S}_\bullet\underline{\mathcal M}|,\]
where the last comparison is the application of Lemma~\ref{lemma:homotopyequivalenceforthicknerve} at the beginning of the proof. 
This completes the proof.

\subsection{The key map $\tau$}\label{sec:tau}
Suppose $\underline{\mathcal P}$ to be a small $k$-linear dg category with a strict zero object and finite biproduct.
Its cyclic nerve is $\underline{\mathrm{CN}}_\bullet(\underline{\mathcal P})$ is the following 
\begin{align}
\underline{\mathrm{CN}}_n(\underline{\mathcal P}) = \textstyle\bigoplus_{\substack{X_0,\ldots,X_n\text{ are}\\\text{objects in $\underline{\mathcal P}$}}}\underline{\mathcal P}(X_n,X_0)\otimes_k \underline{\mathcal P}(X_{n-1},X_n) \otimes_k\cdots \otimes_{k}\underline{\mathcal P}(X_{0},X_{1}),\label{formula:Hochschildcomplex}\end{align}
Together with the internal differential and the Hochschild differential, this gives a double complex.
We refer to Appendix~\ref{appendix:cyclicnerve} for the detail of a modified but equivalent simplicial structure that we will use in this section:
If $d_j^\mathrm H$ is the ordinary $j$-th face, then $j$-th modified face $\bar{d}_{j}^{\mathrm{H}}$ at the degree $n$ equals $d_{n-j}^{\mathrm{H}}$ for each $j=0,\ldots,n$, and for the ordinary differential $b_n$, the differential $\bar b_n \coloneqq \sum_{j=0}^{n}(-1)^j\bar d_j^{\mathrm{H}} = (-1)^{n}b_n$ with $\bar b_0 = b_0 = 0.$

The Dennis trace map will be obtained by applying suitable functors to the following $\tau$ map.
With the dg Waldhausen category $\underline{\mathcal M}$ above, we refine McCarthy's map $\tau$ in \cite[Section~4.5]{McCarthy} at the chain level:
\begin{align}\tau_{\underline{\mathcal M}}:\Bbb{C}_\bullet(\underline{\mathrm w\mathcal M}) \to \underline{\mathrm C}_\bullet(\underline{\mathcal M}),\label{formula:tau}\end{align}
Here $\Bbb{C}_\bullet(\underline{\mathrm w\mathcal M})$ denotes the \emph{unnormalized} reduced simplcial chain complex $\widetilde k[\Bbb{N}_\bullet\underline{\mathrm w\mathcal M}]$ with differential $\partial$ in (\ref{formula:differentialforthickennerve}) and the target $\underline{\mathrm C}_\bullet(\underline{\mathcal M})$ denotes the direct sum totalization of the double complex defined below.

With the internal differential $d_{\mathrm{int}}$ as in (\ref{formula:internaldifferential}),  using the homological grading for two degrees, put:
\[\underline{B}(\underline{\mathcal M})_{p,q}\coloneqq \underline{\mathrm{CN}}_p(\underline{\mathcal M})^{-q}\text{ for any $p\geq 0$ and $q\in \Bbb{Z}$ so that }\begin{cases}d_{\mathrm{int}}:\underline{B}(\underline{\mathcal M})_{p,q}\to \underline{B}(\underline{\mathcal M})_{p,q-1};\\
\bar b:\underline{B}(\underline{\mathcal M})_{p,q}\to \underline{B}(\underline{\mathcal M})_{p-1,q}.\end{cases}\]
Here, the upper index $-q$ denotes the sum of degrees of $p+1$ mapping complexes factors.
By Lemma~\ref{lemma:dgcyclicnerve}, $d_\mathrm{int}$ commutes with $\bar b.$
The \emph{direct sum} total complex therefore has
\[\textstyle \underline{\mathrm C}_n(\underline{\mathcal M})\coloneqq (\mathrm{Tot}\,\underline{B}(\underline{\mathcal M}))_n = \bigoplus_{p+q=n}\underline{B}(\underline{\mathcal M})_{p,q},\quad D|_{\underline{B}(\underline{\mathcal M})_{p,q}}=d_\mathrm{int} + (-1)^{q}\bar b.\]
The total differential that lower the total degree by $1$.
This homological convention is the regrading of the cohomological total complex in Appendix~\ref{appendix:cyclicnerve}, where elements in $\underline{B}(\underline{\mathcal M})_{p,q}$ has cohomological total degree $-(p+q).$
Consider $\Phi = (-1)^{p(p+1)/2}\cdot \id_{\underline{B}(\underline{\mathcal M})_{p,q}}$ on $\underline{B}(\underline{\mathcal M})_{p,q}$ for each $p,q$.
For the ordinary total differential $D_b|_{B_{p,q}} = d_\mathrm{int}+(-1)^qb$, it defines a chain isomorphism $(\mathrm{Tot}\,\underline{B}(\underline{\mathcal M}),D_b)\to (\mathrm{Tot}\,\underline{B}(\underline{\mathcal M}),D)$.

Fix an integer $n$ and a $n$-simplex $F:E^n\to \underline{\mathrm w\mathcal M}$ as a simplicial map.
Write $X_j=F(j)$.
In general, for each partition $0\leq \pi(0)<\pi(1)<\ldots<\pi(s)\leq n$ of $0,\ldots,n$, the data of $F$ give tensors of the following form. Write
\[\begin{aligned}a_j&=f_{\pi(j)\cdots \pi(j+1)}&&\in\underline{\mathcal M}(X_{\pi(j)},X_{\pi(j+1)})^{-(\pi(j+1)-\pi(j)-1)}\text{ for }0\leq j<s,\\
a_s&=f_{\pi(s)\cdots n 0\cdots \pi(0)}&&\in\underline{\mathcal M}(X_{\pi(s)},X_{\pi(0)})^{-(\pi(0)-\pi(s)+n)},\end{aligned}\]
where the upper index $-(\pi(j+1)-\pi(j)-1)$ denotes the degree of the mapping complex.
Then the corresponding summand in $\mathrm{Tot}\,\underline{B}(\underline{\mathcal M})$ is 
\begin{equation}
R_{s+1}^{n}(\pi)\,a_s\otimes a_{s-1}\otimes\cdots\otimes a_0,
\label{formula:tauexplicit}
\end{equation}
where $R_{s+1}^{n}(\pi)\in\{\pm1\}$. 
The term in (\ref{formula:tauexplicit}) belongs to the direct summand indexed by $(X_{\pi(0)},\ldots,X_{\pi(s)}).$
Consequently
\[\textstyle T_r(F)=
\sum_{\substack{\pi\subseteq[n]\\|\pi|=n-r+1}}R_{n-r+1}^{n}(\pi)\cdot a_{n-r}\otimes\cdots\otimes a_0\text{\quad for }0\leq r\leq n.\]
Thus $\tau_{\underline{\mathcal M}}(F)$ is supported only on summands indexed by nonempty subsequences of $(X_0,\ldots,X_n)$, including the full tuple.
The signs and the precise summation defining $\tau$ will be specified in the proof of Lemma~\ref{lemma:tauisachainmap}. 
For $s=0$, only the closing factor $a_0$ occurs. In general, (\ref{formula:tauexplicit}) has $s+1$ factors and total internal degree $s-n$.
It therefore lies in $B_{s,n-s}$. According to this double grading, we write
\[\textstyle \tau_{\mathcal M}(F)=\sum_{r=0}^{n}T_r(F),\quad
T_r(F)\in \underline{B}(\underline{\mathcal M})_{n-r,r}.\]
In particular, the term in (\ref{formula:tauexplicit}) contributes to $T_r(F)$ with $r=n-s$.

\begin{example}
Suppose $s=n$ so that $\pi(j) = j$ for all $j$.
Set $R_{n+1}^{n}(\pi) = 1.$
The corresponding internal degree-zero tensor is
\[f_{n0}\otimes f_{(n-1)n} \otimes \cdots \otimes f_{01}\in \underline{\mathcal M}(X_{n},X_0)^0 \otimes \underline{\mathcal M}(X_{n-1},X_n)^0 \otimes \cdots \otimes \underline{\mathcal M}(X_0,X_1)^0.\]
Suppose $A$ (resp. $C$) is an ordinary algebra (coalgebra). 
If $\underline{\mathrm w\mathcal M}$ is the category consists of finitely generated projective modules (injective comodules) and isomorphisms between (co)modules, then maps $f_{i(i+1)}$ for $i=0,\ldots,n-1$ determine $F$. 
The other maps $f_{ij}$ are given by their compositions and inverses, and all higher homotopies are trivial. 
As $f_{n0} = (f_{(n-1)n}\cdots f_{01})^{-1}$ and $T_r=0$ for $r>0$, now  the map $\tau$ recovers that of McCarthy \cite[Section~4.5]{McCarthy}. 
\end{example}

\begin{lemma}\label{lemma:tauisachainmap}
Resume the above notations.
Then there are signs $R_{s+1}^n(\pi)$ for any $s\coloneqq n - r$ and any $(s+1)$-partition $\pi$ of $0,\ldots,n$ so that the following equations 
\begin{align}d_{\mathrm{int}}T_{r+1}(F) + (-1)^{r}\bar bT_r(F) = T_r(\partial F)\text{\quad for $0\leq r\leq n-1,$}\label{formula:chainmapE_r}\tag{$\text{\rm{Eq}}_{n,r}$}\end{align} and hence $D\tau = \tau \partial.$
\end{lemma}

This lemma says that if we treat $T_0$ as the initial map, then $T_0$ need not be a chain map from $\Bbb{C}_\bullet(\underline{\mathrm w\mathcal M})$ to the Hochschild chain complex $\underline{\mathrm C}_\bullet(\underline{\mathcal M})$.
The obstruction, namely $T_0(\partial F) - \bar bT_{0}(F)$, is a boundary with respect to the internal differential $d_\mathrm{int}.$
This results in the correction term $T_1(F).$
But $T_0+T_{1}$ need not yet be a chain map.
In general, $T_{r+1}$ corrects the remaining defect of $T_0+\cdots+T_r$ in the next bidegree. For an $n$-simplex, the process ends with $T_n(F)$.

\begin{proof}
We prove the desired equation by double induction.
Fix $F:E^n\to\underline{\mathrm w\mathcal M}$ to be an $n$-simplex with ordered vertices $X_{0},\ldots,X_{n}$, where $0,\ldots,n$ denote its source-vertex labels. 
Every tensor below is regarded as an element of the corresponding direct summand of $\underline{\mathrm CN}_\bullet(\underline{\mathcal M})$, with the inclusion understood following (\ref{formula:Hochschildcomplex}). 
The formulas are defined on individual simplices and extended $k$-linearly. 
Faces and degeneracies delete or repeat the corresponding labels.

\noindent \text{\bf Step~1: Initial step for induction.} 
Set $R_{n+1}^n(\pi) = 1$ for the full cut set $\pi = (0,\ldots,n)$ so that 
\[T_0(F) = f_{n0}\otimes f_{n-1,n}\otimes \cdots \otimes f_{01}.\]
For $n=0$, this means $T_0(F) = \id_{X_0}$.
Since the edges are closed of degree zero, $d_{\mathrm{int}}T_0(F) = 0$.
Assume now that $n\geq 1$. 
For $1\leq j\leq n-1$, we have
\begin{align}\bar d_j^\mathrm{H}T_0(F) = f_{n0}\otimes \cdots \otimes (f_{j,j+1}f_{j-1,j})\otimes \cdots \otimes f_{01}.\label{formula:facecomputationofHH}\end{align}
whereas the sign $R_n^{n-1}$ on the full cut set in dimension $n-1$ gives 
\begin{align}
T_0(d_jF) = f_{n0}\otimes \cdots \otimes f_{j-1,j+1}\otimes \cdots \otimes f_{01}.\label{formula:facecomputationofthickenednerve}\end{align}
Let $h_j(F)$ be obtained from (\ref{formula:facecomputationofHH}) by replacing the displayed composite by $f_{j-1,j,j+1}.$
At the two ends, put 
\begin{equation*}
\begin{array}{rcl}
h_0(F) & = & f_{n,0,1}\otimes f_{n-1,n}\otimes \cdots \otimes f_{1,2},\\
h_n(F) & = & f_{n-1,n,0}\otimes f_{n-2,n-1} \otimes \cdots \otimes f_{01}. 
\end{array}
\end{equation*}
When $n=1$, empty strings of tensor factors are omitted.
The formula (\ref{formula:formulaofthickenednerve}) for the internal differential, with $df_{ijk} = f_{ik} - f_{jk}f_{ij}$, gives
\[d_{\mathrm{int}}h_j(F) = T_0(d_jF) - \bar d_{j}^{\mathrm H}T_0(F)\text{\quad for }0\leq j\leq n.\]
Hence if we put $T_1(F) = \sum_{j=0}^n(-1)^j h_j(F),$ then ($\mathrm{Eq}_{n,0}$) holds.
This means: we set $R_n^n(\{0,\ldots,\widehat{j},\ldots,n\}) = (-1)^j.$

\noindent \text{\bf Step~2: Preliminary computations for induction.} 
We then come to the induction step.
The outer induction is on $n$; for fixed $n$, the inner induction is on $r$, i.e., downward on the number of cuts.
Suppose the coefficients $R^{m}_{r}(\pi)$ have been constructed for $m<n$, and that $T_0(F),\ldots, T_r(F)$ have been constructed in dimension $n.$
With $s = n-r$, the signs $R_s^{n-1}(\pi)$ in $T_r(\partial F)$ and $R_{s+1}^{n}(\pi)$ are therefore already determined for any possible $\pi$.
We need to check that we can choose $R_{s}^n(\pi)$ to define $T_{r+1}(F)$ for all possible $\pi$.
By an $s$-partition we mean a choice of $\pi= \{\pi(0),\ldots,\pi(s-1)\}$ with $0\leq \pi(0)< \cdots <\pi(s-1)\leq n.$
For the cyclic block, we will use the following conventions throughout 
\[i_{u+n+1} = i_u,\quad \pi(s) = \pi(0)+n+1,\quad a_j = f_{\pi(j)\cdots \pi(j+1)}\text{ for }0\leq j<s.\]
Then the tensor given by $\pi$ is $a_{s-1}\otimes \cdots\otimes a_0$.
If $l$ avoids $\pi$, i.e., $l\notin \pi$ for all $j$, put 
\[\underline{l}=\begin{cases}l+n+1 & l<\pi(0);\\
l & l>\pi(0),\end{cases}\text{\quad and $h$ to be the index such that }\pi(h)<\underline{l}<\pi(h+1).\]
This $h$ is unique.
We would use $k$ for index and similarly put $\underline{k}$ and $j$ when $k$ avoids $\pi$.

For a cut $\pi$ avoiding $l$, let $\overline{\pi}$ denote its reindexing after deleting $l$.
Equivalently, for $\delta^l:[n-1]\hookrightarrow [n]$ skipping $l$, $\delta^l\overline{\pi} = \pi.$
The right hand side of (\ref{formula:chainmapE_r}) is
\begin{equation}
\begin{array}{c}
\displaystyle T_r(\partial F) = \sum_{l=0}^n \sum_{\substack{\text{$\pi:$ $s$-partition of $[n]$}\\
l\notin \pi}}(-1)^l R_s^{n-1}(\overline{\pi}) \cdot a_{s-1}\otimes \cdots \otimes f_{\pi(h)\cdots \widehat{\underline{l}}\cdots \pi(h+1)}\otimes \cdots \otimes a_0.
\end{array}\label{formula:rightsideofcompatibility}
\end{equation}
When $h = s-1$, this displayed factor replaces $a_h$ and is the first tensor factor and in particular includes the cyclic tensor factor.

For the left side of (\ref{formula:chainmapE_r}), we write the corresponding internal differential explicitly following (\ref{formula:formulaofthickenednerve})
\begin{equation}
\begin{array}{rcl}
\displaystyle df_{\pi(j)\cdots\pi(j+1)} & = & \displaystyle \sum_{k=\pi(j)+1}^{\pi(j+1)-1}(-1)^{k-\pi(j)-1}f_{\pi(j)\cdots \hat k \cdots \pi(j+1)}\\
&& \displaystyle + \sum_{k=\pi(j)+1}^{\pi(j+1)-1}(-1)^{1+(\pi(j+1)-\pi(j))(k-\pi(j)-1)}f_{k\ldots \pi(j+1)}f_{\pi(j)\cdots k}.
\end{array}\nonumber
\end{equation}
The indices in the cyclic factor are understood using the above periodic convention.
The tensor-product differential contributes the additional sign
\[\textstyle (-1)^{\sum_{u=j+1}^{s-1}|a_u|}\text{ with }\sum_{u=j+1}^{s-1}|a_u| = s-2-j+\pi(j+1)-\pi(0)-n.\]
In particular, for $j = s-1$, this sum is zero, as expected.
Then the left side of (\ref{formula:chainmapE_r}) is the sum of the following three expressions:
\begin{align}
&\displaystyle\sum_{\text{$\pi$: $s$-partition of $[n]$}}\sum_{l\notin \pi}R_s^n(\pi)(-1)^{\pi(0)-\pi(h+1)-\pi(h)+\underline{l}+n-s+h+1}\nonumber\\
&\hspace{5cm}\cdot a_{s-1}\otimes \cdots \otimes f_{\pi(h)\cdots \widehat{\underline{l}}\cdots \pi(h+1)}\otimes \cdots\otimes a_0,\label{formula:deletion}\\
&\displaystyle\sum_{\text{$\pi$: $s$-partition of $[n]$}}\sum_{k\notin\pi}R_s^n(\pi)(-1)^{\pi(0)-\pi(j+1)+(\pi(j+1)-\pi(j))(\underline{k}-\pi(j)-1)+n-s+j+1}\nonumber\\
&\hspace{5cm}\cdot a_{s-1}\otimes \cdots \otimes (f_{\underline{k}\cdots \pi(j+1)}f_{\pi(j)\cdots \underline{k}})\otimes \cdots \otimes a_0\text{ and}\label{formula:splitting}\\
&\displaystyle(-1)^r\bar b\bigg(\sum_{\text{$\pi$: $(s+1)$-partition of $[n]$}}R_{s+1}^n(\pi)a_s\otimes\cdots\otimes a_0\bigg).\label{formula:Hochschildterms}
\end{align}
By the formula for $\bar b$, the Hochschild boundary in (\ref{formula:Hochschildterms}) is
\begin{equation*}
\begin{array}{cl}
&\displaystyle(-1)^r \sum_{\text{$\pi$: $(s+1)$-partition of $[n]$}}R_{s+1}^n(\pi)\bigg((-1)^{(r+1)(1+\pi(0)-\pi(1))}(a_0a_s)\otimes a_{s-1}\otimes \cdots \otimes a_{1}\\
&\displaystyle\hspace{4.7cm}+\sum_{j=0}^{s-1}(-1)^{j+1}a_s\otimes\cdots\otimes(a_{j+1}a_j)\otimes\cdots\otimes a_0\bigg).
\end{array}
\end{equation*}
The sign $(-1)^{(r+1)(1+\pi(0)-\pi(1))}$ in the expression comes from moving $a_0$ past all the other factors.
Precisely, since $|a_0| = 1+\pi(0)-\pi(1)$ and $\sum_{j=0}^{s}|a_j| = s-n = -r$, we have
\[\textstyle|a_0|(\sum_{j=1}^{s}|a_j|) = |a_0|(-r-|a_0|) \equiv (r+1)|a_0|\mod 2\] and hence the sign as above.

It remains to choose $R_{s}^n(\pi)$ so that the splitting terms (\ref{formula:splitting}) cancel the Hochschild term (\ref{formula:Hochschildterms}), and then to show that the deletion terms (\ref{formula:deletion}) agree with (\ref{formula:rightsideofcompatibility}).
The signs $R_{s+1}^n$ and $R_{s}^{n-1}$ needed for these comparisons are already known by the two induction hypotheses.

\noindent \text{\bf Step~3: Construction of the sign $H_n(\pi,k)$.} 
For an $s$-partition $\pi$ and $k\notin\pi$, define
\[R_s^n(\pi):=H_n(\pi,k)R_{s+1}^n(\pi\cup\{k\}),\]
where the sign $H_n(\pi,k)$ is given by
\begin{equation}
H_n(\pi,k)=\begin{cases}(-1)^{\pi(0)-\pi(j+1)+1+(\pi(j+1)-\pi(j))(k-\pi(j)-1)},&k>\pi(0),\\
(-1)^{1+(\pi(0)-k)(\pi(s-1)-k+s+1)},&k<\pi(0).\end{cases}\tag{H}\label{formula:constructionofsignH}\end{equation}
In the first case, $\pi(j)<k<\pi(j+1)$, including the possibility $j=s-1$. 
Since $H_n(\pi,k)^2=1$, the defining relation is equivalently
\[R_{s+1}^n(\pi\cup\{k\})=H_n(\pi,k)R_s^n(\pi).\]

For $k>\pi(0)$, comparison of (\ref{formula:splitting}) with the corresponding noncyclic face in (\ref{formula:Hochschildterms}) gives the first sign in (\ref{formula:constructionofsignH}). For $k<\pi(0)$, put
\[u=\pi(0)-k,\quad v=n+1+k-\pi(s-1),\quad m=u+v,\quad r=n-s.\]
The splitting coefficient is $R_s^n(\pi)(-1)^{1+m(v-1)}$, and the corresponding cyclic Hochschild coefficient is $R_{s+1}^n(\pi\cup\{k\})(-1)^{r+(r+1)(1-u)}$.
Their cancellation requires
\[\begin{aligned}H_n(\pi,k)&=(-1)^{1+m(v-1)+(r+1)u}\\
&=(-1)^{1+u(v+r)}=(-1)^{1+(\pi(0)-k)(\pi(s-1)-k+s+1)},\end{aligned}\]
which is the second sign in (\ref{formula:constructionofsignH}). Thus this definition makes (\ref{formula:splitting}) cancel with (\ref{formula:Hochschildterms}), provided it is independent of the choice of $k$.

We now show that $R_s^n(\pi)$ depends only on the $s$-partition $\pi$. 
It suffices to check, for distinct $k,l\notin\pi$, that
\begin{equation}
H_n(\pi,k)H_n(\pi\cup\{k\},l)=H_n(\pi,l)H_n(\pi\cup\{l\},k).\tag{Indep}\label{formula:independentchoice}\end{equation}
Indeed, any two orders of inserting the missing cuts differ by adjacent interchanges. 
Starting with $R_{n+1}^n=1$, (\ref{formula:independentchoice}) therefore makes the recursive construction well defined. 
For a cut set missing only $j$, (\ref{formula:constructionofsignH}) gives $R_n^n=(-1)^j$, in agreement with the base case.

Suppose first that $\pi(0)<k<l$. If $k$ and $l$ lie in different blocks, each insertion leaves the other sign unchanged. 
If $\pi(j)<k<l<\pi(j+1)$, put
\[m=\pi(j+1)-\pi(j),\quad u=k-\pi(j)-1,\quad v=l-\pi(j)-1.\]
The left-hand side of (\ref{formula:independentchoice}) is $(-1)^{mu+(m-u-1)(v-u-1)}$, and the right-hand side is
\[(-1)^{\pi(j+1)+l+mv+(v+1)u}=(-1)^{m+v+1+mv+(v+1)u}.\]
Using $u^2\equiv u\pmod2$, these exponents agree. 
This also covers two positions in the part of the closing block after $\pi(s-1)$.

It remains to check (\ref{formula:independentchoice}) when an insertion changes the smallest cut. 
Write $A=\pi(0)$ and $B=\pi(s-1)$, and assume $k<l$.
If $k<A<l$, put $u=A-k$. Directly from (\ref{formula:constructionofsignH}), both ratios
\[\frac{H_n(\pi\cup\{k\},l)}{H_n(\pi,l)}\quad\text{and}\quad\frac{H_n(\pi\cup\{l\},k)}{H_n(\pi,k)}\]
are $(-1)^u$ when $l<B$, and $(-1)^{u(l-B+1)}$ when $l>B$. 
Hence (\ref{formula:independentchoice}) holds in either case.
If $k<l<A$, put \[u=A-k,\quad v=A-l,\quad w=B-A+s+1.\]
The two sides of (\ref{formula:independentchoice}) have respective exponents
\begin{equation*}\begin{array}{cl}&[1+u(w+u)]+[1+u(u-v)],\\
&[1+v(w+v)]+[1+(u-v)(w+u+1)].\end{array}\end{equation*}
Both are congruent to $u(w+v)$ modulo $2$. This proves (\ref{formula:independentchoice}) in all cases.

\noindent \text{\bf Step~4: Compatibility with deletion part (\ref{formula:deletion}).} 
\noindent\emph{}
It suffices to show that, with $R_s^n(\pi)$ defined above, the deletion part (\ref{formula:deletion}) agrees with the right-hand side (\ref{formula:rightsideofcompatibility}). Comparing their coefficients, the required relation is
\begin{equation}\begin{array}{rcl}
R_s^{n-1}(\overline{\pi})&=&V_s(\pi,l)R_s^n(\pi),\\
V_s(\pi,l)&:=&(-1)^{\pi(0)-\pi(h+1)-\pi(h)+\underline l-l+n-s+h+1}.
\end{array}\tag{V}\label{formula:constructionofsignV}\end{equation}
Here $l\notin\pi$, and $h$ is determined by $\pi(h)<\underline l<\pi(h+1)$. 
In particular, when $l<\pi(0)$, this gives
\[V_s(\pi,l)=(-1)^{n-\pi(s-1)}.\]

For $s=n$, the only cut set avoiding $l$ is $[n]\setminus\{l\}$. Then $R_n^n=(-1)^l$, $R_n^{n-1}=1$, and $V_n=(-1)^l$, so (\ref{formula:constructionofsignV}) holds. For $s<n$, let $\pi'$ be an $s$-partition avoiding $l$, choose $k\notin\pi'$ with $k\ne l$, and set $\pi=\pi'\cup\{k\}$. 
Let $\overline{\pi}'$ and $\overline{\pi}$ be their reindexings after deleting $l$, and put
\[
\overline{k}=(\delta^l)^{-1}(k)=\begin{cases}k,&k<l,\\k-1,&k>l.\end{cases}
\]
Consider the following diagram:
\usetikzlibrary {arrows.meta}
\begin{equation}
\begin{minipage}{.5\textwidth}
\begin{tikzpicture}

\node (A) at (0,1.2) {$R_s^{n-1}(\bar{\pi}')$};
\node (B) at (5,1.2) {$R_{s+1}^{n-1}(\bar{\pi})$};
\node (C) at (0,0) {$R_s^n(\pi')$};
\node (D) at (5,0) {$R_{s+1}^n(\pi)$};

\draw[arrows = {-Stealth[inset=0pt, length=5pt, angle'=30]}] (A) -- node[above] {\footnotesize $H_{n-1}(\overline{\pi}',\overline{k})$} (B);
\draw[arrows = {-Stealth[inset=0pt, length=5pt, angle'=30]}] (C) -- node[above] {\footnotesize $H_n(\pi',k)$} (D);

\draw[dotted, arrows = {-Stealth[inset=0pt, length=5pt, angle'=30]}] (C) -- node[left] {\footnotesize $V_s(\pi',l)$} (A);
\draw[arrows = {-Stealth[inset=0pt, length=5pt, angle'=30]}] (D) -- node[right] {\footnotesize $V_{s+1}(\pi,l)$} (B);
\end{tikzpicture}\end{minipage}\label{formula:commutativediagraminduction}\end{equation}
The arrows mean multiplication by $\pm1$. 
All solid arrows are known by the two induction hypotheses and the definition of $R_s^n$. 
Thus the desired dotted arrow follows once we check that the square commutes, that is,
\[\frac{H_{n-1}(\overline{\pi}',\overline{k})}{H_n(\pi',k)}=\frac{V_s(\pi',l)}{V_{s+1}(\pi,l)}.\]

Suppose first that $\pi'(j)<l<k<\pi'(j+1)$ and $k>\pi'(0)$. Put
\[A=\pi'(0),\quad a=\pi'(j),\quad b=\pi'(j+1),\quad m=b-a,\quad u=k-a-1.\]
Here $\overline{k}=k-1$, and the right endpoint of the block becomes $b-1$. 
From (\ref{formula:constructionofsignH}),
\[\frac{H_{n-1}(\overline{\pi}',\overline{k})}{H_n(\pi',k)}=(-1)^{[A-(b-1)+1+(m-1)(u-1)]-[A-b+1+mu]}=(-1)^{m+u}=(-1)^{b-k-1}.\]
From (\ref{formula:constructionofsignV}), similarly,
\[\frac{V_s(\pi',l)}{V_{s+1}(\pi,l)}=(-1)^{[A-b-a+n-s+j+1]-[A-k-a+n-(s+1)+j+1]}=(-1)^{b-k-1}.\]
If $a<k<l<b$, then $\overline{k}=k$, and the same calculation gives
\begin{equation*}
\begin{array}{rcl}
\displaystyle\frac{H_{n-1}(\overline{\pi}',\overline{k})}{H_n(\pi',k)}&=&(-1)^{[A-(b-1)+1+(m-1)u]-[A-b+1+mu]}=(-1)^{k-a},\\
\displaystyle\frac{V_s(\pi',l)}{V_{s+1}(\pi,l)}&=&(-1)^{[A-b-a+n-s+j+1]-[A-b-k+n-(s+1)+(j+1)+1]}
=(-1)^{k-a}.
\end{array}    
\end{equation*}

We also spell out the case $k<l<\pi'(0)$, where the smallest cut changes. Write $A=\pi'(0)$, $B=\pi'(s-1)$, and put
\[
u=A-k,\quad w=B-k+s+1.
\]
Then $\overline{k}=k$, $\overline{\pi}'(0)=A-1$, and $\overline{\pi}'(s-1)=B-1$. The second formula in (\ref{formula:constructionofsignH}) gives
\[
\frac{H_{n-1}(\overline{\pi}',\overline{k})}{H_n(\pi',k)}
=(-1)^{[1+(u-1)(w-1)]-[1+uw]}
=(-1)^{A+B+s}.
\]
Since $l<A$ before insertion, but $k<l<A$ afterwards, (\ref{formula:constructionofsignV}) gives
\[\frac{V_s(\pi',l)}{V_{s+1}(\pi,l)}=(-1)^{(n-B)-(n-s-A)}=(-1)^{A+B+s}.\]

For completeness, all cases are recorded on the next page. 
Each entry is the common exponent, modulo $2$, of the two ratios. 
They follow by substituting the reindexed positions in (\ref{formula:constructionofsignH}) and (\ref{formula:constructionofsignV}). 
The cases in which $l<\pi'(0)$ use the lifted position in (\ref{formula:constructionofsignV}).

\noindent
Write $A=\pi'(0)$ and $B=\pi'(s-1)$. In the first six rows, $k>A$ and $a=\pi'(j)<k<b=\pi'(j+1)$; if $j=s-1$, then $b=A+n+1$. Throughout, $k,l\notin\pi'$ and $k\ne l$.
\begin{center}
\renewcommand{\arraystretch}{1.12}
\begin{tabular}{@{}ll@{}}
\toprule
Position of the inserted and deleted vertices & Common exponent\\
\midrule
$a<l<k<b$ & $b-k-1$\\
$a<k<l<b$ & $k-a$\\
$A<l<a$ & $1$\\
$b<l\leq n$ & $0$\\
$l<A$, $j<s-1$ & $0$\\
$l<A$, $j=s-1$ & $k-B$\\
\midrule
$l<k<A$ & $0$\\
$k<l<A$ & $A+B+s$\\
$k<A<l<B$ & $A-k$\\
$k<A$ and $B<l$ & $0$\\
\bottomrule
\end{tabular}
\end{center}
Thus the square (\ref{formula:commutativediagraminduction}) commutes in every case. The downward induction on $s$ proves (\ref{formula:constructionofsignV}), so the deletion part (\ref{formula:deletion}) agrees with (\ref{formula:rightsideofcompatibility}). Together with the cancellation of (\ref{formula:splitting}) and (\ref{formula:Hochschildterms}), this proves $\mathrm{Eq}_{n,r}$ and completes the double induction.

Finally, $T_r(F)$ has bidegree $(n-r,r)$, and the total differential on that bidegree is $D=d_{\mathrm{int}}+(-1)^r\bar b$. 
Since $d_{\mathrm{int}} T_0(F)=0$, $\bar bT_n(F)=0$ in cyclic degree zero, and $T_n(\partial F)=0$ by degree, summing $\mathrm{Eq}_{n,r}$ gives the following finite sum
\[D\tau(F)=\tau(\partial F),\quad\tau(F)=\sum_{r=0}^{n}T_r(F).\]
The constant simplex at the strict zero object maps to zero, so the construction descends to $\Bbb{C}_\bullet(\underline{\mathrm w\mathcal M})$.
\end{proof}

For the double complex $\underline{B}(\mathcal M),$ by filtering this complex with degree $p,$ the first page of the associated spectral sequence \[E_{p,q}^1 = H_\mathrm{int}^{-q}(\underline{B}(\underline{\mathcal M})_{p,\bullet}) \coloneqq H^{-q}(\underline{\mathrm{CN}}_p(\underline{\mathcal M}))\text{ with $d^1 = (-1)^{q}\overline{b}$.}\]
Here $H_\mathrm{int}$ means internal cohomology, taken separately in each cyclic degree.

Besides the faces $\bar d^\mathrm H_j$, use the degeneracies $\bar s^\mathrm H_j=s_{n-j}^{\mathrm H}$ in degree $n$ and the forward cyclic operator $t_n^\mathrm H$ of Appendix~\ref{appendix:cyclicnerve}. 
By Lemma~\ref{lemma:dgcyclicnerve}, all these operators are dg maps and satisfy the simplicial and cyclic identities. They therefore induce a cyclic $k$-module structure on $H_\mathrm{int}^0(\underline{B}(\underline{\mathcal M}))$.

As the image under $T_0$ is closed with respect to the internal differential, taking its internal cohomology class defines \[\overline{T}_0:\Bbb{C}_\bullet(\underline{\mathrm w\mathcal M}) \to H_{\mathrm{int}}^0(\underline{B}(\underline{\mathcal M}));\quad F\mapsto [T_0(F)].\]
Equip the source with the filtration by simplicial degree. 
Since $T_r(F)$ has cyclic degree $n-r$, the map $\tau$ preserves these filtrations, and its filtration-preserving component is $T_0$.
Hence $\bar T_0$ is the map on the row $q=0$ of the first page. The source $\Bbb{C}_\bullet(\underline{\mathrm w\mathcal M})$ is then regarded as the underlying cyclic $k$-module $\widetilde k[\Bbb N_\bullet\underline{\mathrm w\mathcal M}]$, with the cyclic structure in Lemma~\ref{remark:cyclicsettrivialcase}.
\begin{corollary}\label{Cor:cyclicmap}
    The induced map 
    $\overline{T}_0$ is strict morphism of cyclic $k$-modules.
\end{corollary}

\begin{proof}
For the faces, by Step~1 in the proof of Lemma~\ref{lemma:tauisachainmap}, we have for $F$ being an $n$-simplex for $n\geq 1$
\[T_0(d_jF)-\bar d_j^\mathrm HT_0(F) = d_\mathrm{int}h_j(F).\]
This also includes the two end faces.
A degeneracy $s_j$ on the thick nerve repeats a vertex and the newly added edge is the identity.
Hence $T_0(s_j F)=\bar s_j^\mathrm H T_0(F).$
Finally, the forward cyclic operator on the thick nerve rotates the vertex labels.
For $n\geq 1,$ with $t_n^\mathrm H$ in Appendix~\ref{appendix:cyclicnerve}, we have
\[T_0(t_nF) = f_{01}\otimes f_{n0}\otimes \cdots \otimes f_{12} = t_n^\mathrm HT_0(F).\]
The Koszul sign here is trivial because all factors have internal degree zero. 
In cyclic degree zero, both cyclic operators are the identity.
The map in the lemma is well-defined on the reduced linearization because the constant simplex at the strict zero object maps to zero, by $\id_0 = 0$. 
\end{proof}

We finish this section by showing that the map $\tau$ is moreover monoidal.
Choose the biproducts suitably so that $X\mapsto X\oplus0$ and $X\mapsto 0\oplus X$ are the identity functor.

\begin{lemma}
Suppose $\underline{\mathcal P}$ to be a small $k$-linear dg category with a strict zero object and finite biproduct.
We have the following diagram for the additivity that is commutative up to homotopy
\[\xymatrix@=1.5em{
\underline{\mathrm C}_\bullet(\underline{\mathcal P}\times \underline{\mathcal P})\ar[r]^{p\hspace{.3cm}}_{\simeq\hspace{.3cm}}\ar[d]_{\underline{\mathrm C}_{\bullet}(\oplus)} & \underline{\mathrm C}_\bullet(\underline{\mathcal P})\oplus \underline{\mathrm C}_\bullet(\underline{\mathcal P})\ar[ld]^{+}\\
\underline{\mathrm C}_\bullet(\underline{\mathcal P}),}\]
where $+$ means the sum of two chains inside the cyclic nerve as an abelian group, and $\underline{\mathrm C}_{\bullet}(\oplus)$ is given by the map $\underline{\mathcal P}\times \underline{\mathcal P}((X_1,X_2),(Y_1,Y_2))\to \underline{\mathcal P}(X_1\oplus X_2,Y_1\oplus Y_2).$
\end{lemma}

\begin{proof}
By Proposition~\ref{prop:additivity}, the map $p$ has a homotopy inverse $\iota$ with $p\iota=\id$, where $\iota$ is induced by functors $X\mapsto (X,0)$ and $X \mapsto (0,X)$ from $\underline{\mathcal P}$ to $\underline{\mathcal P}\times \underline{\mathcal P}$. 
Let $K$ be the total-complex homotopy constructed there so that $DK+KD=\iota p-\id.$
By the definition of these maps and that $X\mapsto X\oplus 0$ and $X\mapsto 0\oplus X$ are the identity,
$\underline{\mathrm C}_\bullet (\oplus)\iota=+$. Therefore, with $L=-\underline{\mathrm C}_\bullet (\oplus)K$, we have
$DL+LD=\underline{\mathrm C}_\bullet (\oplus)-+\,p.$
This is the required homotopy. In particular, the comparison is a chain homotopy, rather than merely a quasi-isomorphism of complexes.
\end{proof}

\begin{proposition}\label{prop:taumonoidal}
The following diagram commutes
\[\xymatrix@=1.5em{
\Bbb{C}_\bullet(\underline{\mathrm w\mathcal M}\times \underline{\mathrm w\mathcal M})\ar[rr]^{\tau_{\underline{\mathcal M}\times \underline{\mathcal M}}}\ar[d]_{\Bbb{C}_\bullet(\underline{\mathrm w}\oplus)} && 
\underline{\mathrm C}_\bullet(\underline{\mathcal M}\times \underline{\mathcal M})\ar[r]^{p\hspace{.3cm}}_{\simeq\hspace{.3cm}}\ar[d]_{\underline{\mathrm C}_{\bullet}(\oplus)} & \underline{\mathrm C}_\bullet(\underline{\mathcal M})\oplus \underline{\mathrm C}_\bullet(\underline{\mathcal M})\ar[ld]^{+}\\
\widetilde{k}[\Bbb{N}_\bullet\underline{\mathrm w\mathcal M}]\ar[rr]^{\tau_{\underline{\mathcal M}}} && \underline{\mathrm C}_\bullet(\underline{\mathcal M}).}\]
In particular, the map $\tau$ is monoidal up to homotopy.
\end{proposition}

\begin{proof}
By the above lemma, it suffices to check the commutativity of the left square. 
For $\Bbb{C}_\bullet(\underline{\mathrm w}\oplus)$, if $F,G:E^n \to \underline{\mathrm w\mathcal M}$ are two $n$-simplices, then we have $(F\oplus G)(i_0,\ldots, i_q) = f_{i_{0}\cdots i_{q}}\oplus g_{i_0\cdots i_q}.$
Hence for fixed $s$ and $\pi$ in (\ref{formula:tauexplicit}),  $\tau_{\underline{\mathcal M}}\circ \Bbb C_\bullet(\underline{\mathrm w}\oplus)(F,G)$ gives the summand $R_{s+1}^n(\pi )((a_s\oplus b_s)\otimes \cdots \otimes (a_0\oplus b_0)).$
On the other hand, we have $\tau_{\underline{\mathcal M}\times \underline{\mathcal M}}(F,G)$ gives the summand $R_{s+1}^n((a_s, b_s)\otimes \cdots (a_0, b_0))$ and the image under $\underline{\mathrm C}_\bullet(\oplus)$ is the same summand.
\end{proof}

\begin{remark}\label{remark:withoutSegal}
This enables us, after realization, to extend the map induced by $\tau$ through Segal's group-completion construction for direct sum. 
Here one uses the naturality of $\tau$ on the finite direct-sum diagrams, together with the finite-product homotopy equivalences of Appendix~\ref{appendix:cyclicnerve}. 
We will not pursue this route further, although McCarthy's comparison diagram tells a beautiful story \cite[Section 4.5]{McCarthy}. 
Our reason is that we have not established a comparison between Segal's geometric realization and $K$-theory space. 
Such a comparison would be an interesting question for further study.
\end{remark}

\section{The explicit Dennis trace map}\label{sec:explicitDennistrace}
We are ready to define the dg Dennis trace map using $\tau$.
We hope to follow McCarthy's proof of his Dennis trace map being identical to the ordinary one and give a refined dg version of his commutative diagram (See Theorem~\ref{thm:Dennistracediagram}).
As we mentioned in Remark~\ref{remark:withoutSegal}, we will not include Segal's geometric realization.
For this we need the first delooping, that is Lemma~\ref{lemma:delooping} for the first delooping. 
The explicit map $\tau$ also allows us to extend the derived Koszul equivalence to a commutative diagram in Theorem~\ref{thm:compatibility}: the Dennis trace is compatible with the derived Koszul equivalence.

\subsection{The delooping}\label{sec:delooping}
Throughout this section, all total complexes use direct sums. 
Retain the notation $(\underline{\mathrm C}_\bullet(\underline{\mathcal M}),D)$ of Section~\ref{sec:tau} for the \emph{unnormalized} total complex of the dg cyclic nerve. 
Put
\[\textstyle \underline{\mathrm C}_S(\underline{\mathcal M}):=\mathrm{Tot}\,\underline{\mathrm C}_\bullet(S_\bullet\underline{\mathcal M}),\quad\underline{\mathrm C}_S(\underline{\mathcal M})_n=\bigoplus_{p+q=n}\underline{\mathrm C}_q(\underline{S}_p\underline{\mathcal M}),\]
where total complex is in the $\underline{S}_\bullet$-direction and $\underline{\mathrm C}_\bullet$-direction. The differential is
\[\textstyle D_S=D+(-1)^q b^S\quad\text{on }\underline{\mathrm C}_q(S_p\underline{\mathcal M}),\quad b^S=\sum_{j=0}^{p}(-1)^j\underline{\mathrm C}_\bullet(d_j^S),\] where $d_j^S$ denotes the faces in the $\underline{S}_\bullet$-direction (See Definition~\ref{def:dgSconstruction}) and $b^S=0$ when $p=0$. 

\begin{lemma}\label{lemma:delooping}
For the dg categories $\underline{\mathcal M}$ considered above, there is a quasi-isomorphism
\[\Sigma \underline{\mathrm C}_\bullet(\underline{\mathcal M}) \simeq \underline{\mathrm C}_S(\underline{\mathcal M}),\]
where $\Sigma \coloneqq [-1]$ is the homological suspension that is the degree shifting functor so that $(\Sigma M)_n = M_{n-1}$ and $d_{\Sigma M} = -d_M$ for any chain complex $(M,d).$
\end{lemma}

The key tool is derived Morita invariance of Hochschild homology; see \cite[Sections~3.8 and~5.3]{Keller}. 
These Hochschild homology result works over a general commutative base ring, but requires that the relevant dg categories are $k$-flat which is automatic over a field.

\begin{proof}
Recall the additive bar construction $B_p(M) = M^{\oplus p}$ and $B_0(M)$ for a complex $M$ which is a simplicial chain complex such that for an order-preserving $\theta:[m]\to[n]$ we have 
\begin{align}B(\theta):M^{\oplus n}\to M^{\oplus m};\,\, 
(x_1,\ldots,x_n) \mapsto \bigg(\sum_{\theta(0)< i\leq \theta(1)}x_i,\ldots, 
\sum_{\theta(m-1)<i\leq \theta(m)}x_i\bigg)\label{formula:barsimplicials}\end{align}
and empty sums are zero.
The degenerate subcomplex of $B_\bullet(\underline{\mathcal M})$ is contractible, and the normalization satisfying $\mathrm{Norm}_p^B B(M) = M$ if $p=1$ and $=0$ otherwise. This contraction gives a natural chain homotopy equivalence $j_M:\Sigma M\overset{\simeq}{\to}\mathrm{Tot}\,B_\bullet(M)$ that maps $x\in (\Sigma M)_{q+1} = M_q$ to $(-1)^{q}x$, where the total differential on $B_p(M)_q$ is $d_M+(-1)^qb_B$ with the alternating sum $b_B$ of bar faces. 

Let $e_i$ be the dg functor $e_i^n:\underline{\mathcal M}\to \underline{S}_{n}\underline{\mathcal M};\,\, X\mapsto (0\to \cdots \to 0\to X\overset{\id}{\to}\cdots \overset{\id}{\to}X)$ so that
\[(e_i^nX)_{ab} = \begin{cases}X& a<i\leq b;\\0& \text{otherwise}.\end{cases}\]
This functor on a morphism $f$ equals $f$ on each nonzero component and zero elsewhere.
Put
\[\alpha_n:\textstyle\bigoplus_{i=1}^n \underline{\mathrm C}_\bullet(\underline{\mathcal M})\to \underline{\mathrm C}(\underline{S}_n \underline{\mathcal M});\quad (x_1,\ldots,x_n)\mapsto \sum_i \underline{\mathrm C}(e_i^n)(x_i)\text{ and $\alpha_0 = 0$}.\]
We are to show that for 
\[\alpha_\bullet:B_\bullet(\underline{\mathrm C}_\bullet(\underline{\mathcal M}))\to \underline{\mathrm C}_\bullet(\underline{S}_\bullet\underline{\mathcal M})\] 
\begin{enumerate}[label=\textup{(\arabic*)},noitemsep,leftmargin=*,nosep]
\item \label{claim:simplicial} is a map of simplicial chain complexes for each $n$; 
\item \label{claim:quasiisomorphism} is a levelwise quasi-isomorphism.
\end{enumerate}
These claims implies that $\mathrm{Tot}(\alpha_\bullet)$ is a quasi-isomorphism, and the desired quasi-isomorphism is the composition $\mathrm{Tot}(\alpha_\bullet)\circ j:\Sigma\underline{\mathrm C}_\bullet(\underline{\mathcal M})\to \mathrm{Tot}\,B_\bullet(\underline{\mathrm C}_\bullet(\underline{\mathcal M}))\to\underline{\mathrm C}_S(\underline{\mathcal M}).$
Indeed, although the mapping complex is unbounded, because the total complex uses the direct sum, each individual element only have finite nonzero component. 
Using this, one can show, by an inductively on the simplicial degree that $\mathrm{cone}(\mathrm{Tot}(\alpha_\bullet))\cong \mathrm{Tot}\,\mathrm{cone}(\alpha_\bullet)$ is acyclic.

For the Claim~\ref{claim:simplicial}, note that all faces and degeneracies for $\underline{S}_\bullet$-construction could be written as 
\begin{align}(\theta^* Y)_{ab} = Y_{\theta(a),\theta(b)}\text{ for some order preserving }\theta:[m]\to [n],\label{formula:Sconstructionfacesanddegeneracies}\end{align}
where $Y_{\bullet,\bullet}$ is an object in the $\underline{S}_\bullet$-construction.
Compute \[(\theta^* e_i^n X)_{ab} = \begin{cases}X & \theta(a)<i\leq \theta(b);\\ 0 & \text{otherwise},\end{cases}.\]
If $\theta(0)<i\leq \theta(m)$, there is a unique $k \in\{1,\ldots,m\}$ such that $\theta(k-1)<i\leq \theta(k),$ and hence $\theta(a)<i\leq \theta(b) \Leftrightarrow a<k\leq b.$
Hence 
\[\theta^* e_i^n = \begin{cases}e_k^m & \theta(k-1)<i\leq \theta(k)\text{ given }1\leq k\leq m;\\
0 & i\leq \theta(0)\text{ or }i>\theta(m).\end{cases}\]
This results in
\begin{align*}
\textstyle\underline{\mathrm C}(\theta^*)\alpha_n(x_1,\ldots,x_n) & \textstyle = \sum_{i=1}^n \underline{\mathrm C}(\theta^* e_i^n)(x_i) = \sum_{k=1}^m\sum_{\theta(k-1)<i\leq \theta(k)}\underline{\mathrm C}(e_k^m)(x_i)\\
& = \textstyle\sum_{k=1}^m \underline{\mathrm C}(e_k^m)\big(\sum_{\theta(k-1)<i\leq \theta(k)}x_i\big)\\
& = \alpha_m B(\theta)(x_1,\ldots,x_n)\text{ (using (\ref{formula:barsimplicials}))}.
\end{align*}
As for the internal $\underline{\mathrm C}_\bullet$-differential, each $\alpha_n$ is a chain map as $e_i^n$ are dg functors. \ref{claim:simplicial} follows.

We show the Claim \ref{claim:quasiisomorphism}. 
The $n=0$ is immediate as $\underline{\mathrm C}_\bullet(S_0\underline{\mathcal M}) = 0.$
Fix $n\geq 1$ and consider the full dg subcategory $\underline{\mathcal U}_n$ of $\underline{S}_n\underline{\mathcal M}$ whose objects are $e_i^n(X)$ for all $1\leq i\leq n$ and $X\in\underline{\mathcal M}.$
For any $m\geq 0$ and $Y_r \coloneqq e_{i_r}^n(X_r)$, we have 
\[\textstyle\underline{\mathrm{CN}}_m(\underline{\mathcal U}_n) = \bigoplus_{Y_0,\ldots,Y_m}\underline{\mathcal U}_n (Y_m,Y_0)\otimes_k \cdots \otimes_k \underline{\mathcal U}_n (Y_0,Y_1),\]
with the usual factor $\underline{\mathcal U}(Y_0,Y_0)$.
As
\[\underline{S}_n\underline{\mathcal M}(e_i^nX,e_j^nX') = \begin{cases}0 & i<j;\\
\underline{\mathcal M}(X,X') & i\geq j,\end{cases}\]
the component in $\underline{\mathrm{CN}}(\underline{\mathcal U}_n)$ is nonzero if $i_m\geq i_0\geq i_1 \geq \cdots \geq i_m$, which means $i_m = i_0 = \ldots = i_{m-1}.$
After summing over $i = 0,\ldots,n$, one obtains the following isomorphism of cyclic objects in complexes of $k$-modules  $\bigoplus_{i=1}^n \underline{\mathrm{CN}}_\bullet(\underline{\mathcal M}) \cong \underline{\mathrm{CN}}_\bullet(\underline{\mathcal U}_n).$
Hence after taking the unnormalized total complexes, we have an isomorphism
\[\beta_n:\textstyle\bigoplus_{i=1}^n\underline{\mathrm C}_\bullet(\underline{\mathcal M}) \cong \underline{\mathrm C}_\bullet(\underline{\mathcal U}_n).\]
For the embedding functor $J_n:\underline{\mathcal U}_n \to \underline{S}_n\underline{\mathcal M}$, we have $\alpha = \underline{\mathrm C}_\bullet(J_n)\circ \beta_n$.

We show that $J_n$ is a Morita equivalence.
For this we use Keller \cite[Lemma~3.10]{Keller} (also the notations there), which claims the desired Morita equivalence.
The functor $J_n$ defines a $\underline{\mathcal U}_n$-$\underline{S}_n\underline{\mathcal M}$-bimodule $M(Z,U) \coloneqq  \underline{S}_n\underline{\mathcal M}(-,J_n(-)).$ 
Put $h_{Y} = \underline{S}_n\underline{\mathcal M}(-,Y)$
The conditions of this lemma are checked by
\begin{itemize}[noitemsep,leftmargin=*,nosep]
\item For the first condition, as $M(-,U) = h_{J_nU}$ is representable, it is perfect. 
\item As the functor $J_n$ is fully faithful, the second condition holds: \[\underline{\mathcal U}_n(X,X') = \underline{S}_n\underline{\mathcal M}(J_n(X),J_n(X')) \cong \underline{\mathrm{Mod}}(\underline{S}_n\underline{\mathcal M})(h_{J_n(X)},h_{J_n(X')}),\]
where $\underline{\mathrm{Mod}}(\underline{S}_n\underline{\mathcal M})$ denotes the dg $\underline{S}_n\underline{\mathcal M}$-module category.
\item 
Since representable modules compactly generate $D(\underline{S}_n\underline{\mathcal M})$, it suffices to check that $h_Y$ for each $Y$ in $\underline{S}_n\underline{\mathcal M}$ belongs to the thick subcategory generated by the modules $h_{J_nU}$.
\end{itemize}
For each $Y = (0 = Y_0\hookrightarrow Y_1\hookrightarrow \cdots \hookrightarrow Y_n),$ consider the filtration
\[F^rY = (0  = Y_0\hookrightarrow Y_1 \hookrightarrow \cdots \hookrightarrow Y_r \overset{\id}{\hookrightarrow} Y_r \overset{\id}{\hookrightarrow} \cdots \overset{\id}{\hookrightarrow}Y_r).\]
On the full $S_n$-diagram, we have $(F^rY)_{ab} = Y_{\min\{a,r\},\min\{b,r\}}$ and hence all entries still in $\underline{\mathcal M}$.
This gives cofibration sequence 
\[0\to F^{r-1}Y \to F^{r}Y \to e_r^n(Q_r) \to 0.\]
As $Y_{r-1}$, $Y_{r}$, and $Q_r$ are cofibrant-fibrant, they are graded projective when $\underline{\mathcal M}=\underline{\mathrm M}_A$ and is graded injective when $\underline{\mathcal M}=\underline{\mathrm{CM}}_C$).
Hence the exact sequence $Y_{r-1}\hookrightarrow Y_r\twoheadrightarrow Q_r$ is split with a graded $A$-linear (or $C$-colinear) section $s_r:Q_r\to Y_r$ and hence a section (need not be a closed morphism) $\widetilde s_r:e_r^n (Q_r) \to F^rY$ by using $s_r$ at the first row (i.e., the $(a=0)$-row) with $b\geq r$ and zero when $a<r.$
If $a<r\leq b,$ then $(e_r^n(Q_r))_{ab} \to (F^rY)_{ab}$ is $Q_r \overset{s_r}{\to} Y_r\to Y_r/Y_a.$
These maps form a natural transformation of the underlying graded (co)module $S_n$-diagrams.

Applying the mapping complexes, for any $Z\in \underline{S}_\bullet \underline{\mathcal M}$, we have the following short exact sequence of mapping complexes
\[0\to \underline{S}_\bullet\underline{\mathcal M}(Z,F^{r-1}Y) \to \underline{S}_\bullet\underline{\mathcal M}(Z,F^{r}Y) \to \underline{S}_\bullet\underline{\mathcal M}(Z,e_{r}^nQ_r)\to 0.\]
Postcomposition with $\widetilde s_r$ splits the last map in every graded degree.
As this sequence is functorial in $Z$, it gives a short exact sequence of right $\underline{S}_\bullet\underline{\mathcal M}$-modules and hence the following distinguished triangle in $D(\underline{S}_\bullet\underline{\mathcal M})$
\[h_{F^{r-1}Y}\to h_{F^{r}Y} \to h_{e_{r}^nQ_r}\to h_{F^{r-1}Y}[1].\]
Inductively, $h_{F^rY}$ is in the thick subcategory generated by $h_{e_{s}^n(Q_s)}$ for $1\leq s \leq r$ and hence $h_Y = h_{F^nY}$ is in the thick subcategory generated by $h_{J_nU}$.
This checks the third condition.
Hence $J_n$ is a derived Morita equivalence.

Finally, as $\underline{S}_n\underline{\mathcal M}$ is $k$-flatness, so does its full dg subcategory $\underline{\mathcal U}_n$, the cyclic-nerve total complexes used here compute their derived Hochschild complexes. 
Keller's Morita-invariance theorem \cite[Section~5.3, Theorem~5.2(a)]{Keller} therefore implies that $\underline{\mathrm C}_\bullet(J_n):\underline{\mathrm C}_\bullet(\underline{U}_n)\longrightarrow\underline{\mathrm C}_\bullet(\underline{S}_n\underline{\mathcal M})$ is a quasi-isomorphism. 
Since $\alpha_n=\underline{\mathrm C}_\bullet(J_n)\beta_n$ and $\beta_n$ is an isomorphism, $\alpha_n$ is a quasi-isomorphism. 
This proves claim (2) and completes the proof of Lemma~\ref{lemma:delooping}.
\end{proof}

\begin{remark}[Higher delooping]\label{remark:Higherdelooping}
Suppose that $\underline{\mathcal P}$ is a small $k$-linear dg Waldhausen category in the setting of Section~3, with a strict zero object and finite biproducts, so that the iterated $\underline{S}_\bullet$-constructions are defined. 
Following McCarthy's additivity and path-construction argument \cite[Sections~3.4 to 3.6, especially Corollary~3.6.3]{McCarthy}, we obtain a quasi-isomorphism
\[\Sigma\underline{\mathrm C}_S(\underline{\mathcal P})\simeq\underline{\mathrm C}_{SS}(\underline{\mathcal P})\text{ with }\underline{\mathrm C}_{SS}(\underline{\mathcal P}):=\mathrm{Tot}\,\underline{\mathrm C}_\bullet(\underline{S}_\bullet \underline{S}_\bullet\underline{\mathcal P}).\]
Here the totalization uses the two independent $\underline{S}_\bullet$-directions and the preceding sign convention. 
The functors (including $\theta^*$ induced by order preserving $\theta:[m]\to [n],$ and the projection functor on the extension category) and simplicial homotopies in this argument are dg, so they commute with the internal differential and pass to the direct sum total complexes. In particular, this applies to the dg categories $\underline{\mathcal M}$ considered above.
\end{remark}

\subsection{The explicit Dennis trace map}\label{sec:explicitDennistraceconstruction}
The goal is to define the Dennis trace map.
For a homological chain complex $E$, we have topological spaces $\mathscr H(E) \coloneqq |\Gamma(E^{\geq 0})| \simeq \Omega^\infty H(E)$, where $H(E)$ is the corresponding Eilenberg-Mac Lane spectrum, $\Gamma$ is the Dold-Kan functor from nonnegatively graded chain complexes to simplicial $k$-modules, and the upper index $(-)^{\geq 0}$ denotes the truncation with degree zero term $\mathrm{ker}(d:E_0\to E_{-1})$ and positive degree terms being those of $E.$

Retain the notations and the categories $\underline{\mathcal M}$, $\underline{\mathrm w\mathcal M}$, and $\underline{\mathrm wS}_\bullet\underline{\mathcal M}$ in Sections~\ref{sec:keyconstruction} and \ref{sec:delooping}.
Put $\mathcal{H}(\underline{\mathcal M})\coloneqq \mathscr H(\underline{\mathrm C}_\bullet(\underline{\mathcal M}))$ and $\mathcal H_S(\underline{\mathcal M})\coloneqq \mathscr H(\underline{\mathrm C}_S(\underline{\mathcal M})).$
Consider the following maps, and their chain complex version $\eta_{\underline{\mathcal M}}$ and $\beta_{\underline{\mathcal M}}$
\[\eta_{\underline{\mathcal M}}^\mathrm t:|\Bbb N_\bullet \underline{\mathrm w\mathcal M}| \to \Omega |\Bbb{N}_\bullet (\underline{\mathrm wS}_\bullet\underline{\mathcal M})|\text{ and }\beta_{\underline{\mathcal M}}^\mathrm t:\mathcal H(\underline{\mathcal M}) \to \Omega\mathcal H_{S}(\underline{\mathcal M}).\]
The bisimplicial set $\Bbb{N}_\bullet (\underline{\mathrm wS}_\bullet\underline{\mathcal M})$ has two simplicial directions, and both the thick nerve direction and the $\underline{S}_\bullet$-direction have the ordinary realization.
The zero object supplies the
basepoint.
They all induced by \[j_1:\underline{\mathcal M} \to \underline{S}_1 \underline{\mathcal M};\quad X\mapsto (0\hookrightarrow X).\]
This induces the chain map \[\Sigma\eta_{\underline{\mathcal M}}:\Sigma\Bbb{C}_\bullet(\underline{\mathrm w\mathcal M})\to \Bbb{C}_S(\underline{\mathrm wS}_\bullet\underline{\mathcal M});\quad \Sigma(x)\mapsto (-1)^q\Bbb{C}_\bullet(j_1)(x)\text{ for $x\in \Bbb{C}_q(\underline{\mathrm w\mathcal M})$}.\]
This is a chain map: their $\underline{S}_\bullet$-direction differentials vanish because $\Bbb{C}_\bullet(\underline{S}_0\underline{\mathcal M}) = 0$, and the factor $(-1)^q$ is for the suspension differential $d_\Sigma=-d$.
We have $\Bbb{N}_\bullet(j_1)$ on $\Bbb{N}_\bullet \underline{\mathrm w\mathcal M}$ and $\eta_{\underline{\mathcal M}}^\mathrm t$ is adjoint to
\[\Sigma |\Bbb{N}_\bullet\underline{\mathrm w \mathcal M}|\to |\Bbb{N}_\bullet (\underline{\mathrm wS}_\bullet\underline{\mathcal M})|;\quad [x,t]\mapsto [\Bbb{N}_\bullet(j_1)(x),t].\]
This is well defined on the reduced suspension: in the target, both endpoints of each path parametrized by $|\Delta^1|$ lie in $|N_\bullet(\underline{wS}_0\underline{\mathcal M})|=*,$ and the interval over the zero-object basepoint maps constantly to the basepoint, since the corresponding edge in the $\underline{S}_\bullet$-direction is degenerate.

The functor $j_1$ induces the quasi-isomorphism in Lemma~\ref{lemma:delooping}: \begin{align}\Sigma\beta_{\underline{\mathcal M}}:\Sigma\underline{\mathrm C}_\bullet(\underline{\mathcal M})\overset{\simeq}{\to}\underline{\mathrm C}_S(\underline{\mathcal M});\quad \Sigma y\mapsto(-1)^q\underline{\mathrm C}_\bullet(j_1)(y)\text{ for $y\in\underline{\mathrm C}_q(\underline{\mathcal M})$}\label{formula:betaM}\end{align} 
with image in $\underline{S}_\bullet$-degree one. 
Applying $\mathscr H$ and the natural suspension-loop identification gives $\beta_{\underline{\mathcal M}}^\mathrm t$, which is a weak equivalence.
Write $\tau_{\underline{\mathcal M}}^{\mathrm t}:|\Bbb{N}_\bullet\underline{\mathrm w\mathcal M}|\to\mathcal H(\underline{\mathcal M})$ and
$\tau_{\underline{S}_\bullet\underline{\mathcal M}}^{\mathrm t}:|\Bbb{N}_\bullet\underline{\mathrm wS}_\bullet\underline{\mathcal M}|\to\mathcal H_S(\underline{\mathcal M})$ for the space maps induced by the chain maps below. The maps defined above give the following diagram:
\begin{align}\vcenter{\xymatrix@=1.5em{
|\Bbb{N}_\bullet \underline{\mathrm w\mathcal M}| \ar[r]^{\eta_{\underline{\mathcal M}}\hspace{.5cm}}\ar[d]_{\tau^\mathrm t_{\underline{\mathcal M}}} & \Omega|\Bbb{N}_\bullet \underline{\mathrm wS}_\bullet\underline{\mathcal M} |\ar[d]^{\Omega\tau^\mathrm t_{\underline{S}_\bullet\underline{\mathcal M}}}\\
\mathcal H(\underline{\mathcal M})\ar[r]^{\simeq\hspace{.5cm}}_{\beta_{\underline{\mathcal M}}\hspace{.5cm}}& \Omega\mathcal H_{S}(\underline{\mathcal M}),
}}\label{formula:Dennistracediagram}\end{align}
where the bottom map is a homotopy equivalence by Lemma~\ref{lemma:delooping}.

\begin{theorem}\label{thm:Dennistracediagram}
Let $\underline{\mathcal M}$ denote the dg enhancement of the Waldhausen category $\mathcal M_A(N)$ or $\mathcal M_C(M)$. 
Let $\underline{\mathrm w\mathcal M}$ and $\underline{\mathrm wS}_\bullet\underline{\mathcal M}$ be the maximal Kan complexes of respectively the dg nerves $\underline{N}_\bullet\mathcal M$ and $\underline{N}_\bullet\underline{S}_\bullet\underline{\mathcal M}.$
Then the above diagram \text{\rm (\ref{formula:Dennistracediagram})} is commutative up to homotopy.
It is natural for exact $k$-linear dg functors preserving the chosen zero objects.
\end{theorem}

\begin{proof}
The $\tau$ map in (\ref{formula:tau}) gives, for every $n\geq 0$, a map  
\[\tau_n\coloneqq \tau_{\underline{S}_n\underline{\mathcal M}}:\Bbb{C}_\bullet(\underline{\mathrm wS}_n\underline{\mathcal M})\to \underline{\mathrm C}(\underline{S}_n\underline{\mathcal M})\]
By Lemma~\ref{lemma:tauisachainmap}, they are chain maps.
In particular, the explicit formula (\ref{formula:tauexplicit}) applies to the case $n=1$.
To apply $\underline{S}_\bullet$-construction to the maps $\tau_n$, we need to check the compatibility between $\tau$ and the simplicial structure on the $\underline{S}_\bullet$-direction.

An order-preserving map $\theta:[m]\to [n]$ induces a dg functor $\theta^*:\underline{S}_n\underline{\mathcal M} \to \underline{S}_m\underline{\mathcal M}$. 
We claim
\begin{lemma}\label{lemma:thicknerveofdgfunctor}
With the notations above, we have
\[\tau_{m}\circ \Bbb{C}_\bullet(\theta^*) = \underline{\mathrm C}_\bullet(\theta^*)\circ \tau_{n}.\]
\end{lemma}
\begin{proof}[Proof of the lemma]
With the notation of Section~4.1, let $G:\underline{\mathcal C}\to\underline{\mathcal C}'$ be a $k$-linear dg functor preserving the chosen zero objects and weak equivalences, e.g., $G=\theta^*$. 
For a $p$-simplex $F$ of the thick nerve, the induced simplex $G_*F$ has
\[(G_*F)(i) = G(F(i)) = G(X_i)\text{ and }(G_*F)(i_0,\ldots,i_r) = G(F(i_0,\ldots,i_r)) = G(f_{i_0\cdots i_r}).\]
On the other hand, for the cyclic nerve, we have 
\[\underline{\mathrm C}_\bullet(G)(x_q\otimes \cdots \otimes x_0) = G(x_q)\otimes\cdots\otimes G(x_0).\]
Hence, for the element in (\ref{formula:tauexplicit}), its image under $\underline{\mathrm C}_\bullet(G)$ writes $R_{s+1}^{q}(\pi)\cdot G(a_s)\otimes\cdots\otimes G(a_{0}).$
We have $\underline{\mathrm C}_\bullet(G)\tau_{\underline{\mathcal C}}(F) = \tau_{\underline{\mathcal C'}} (G_*F)$ for each simplex and the reduced $k$-linearization gives $\underline{\mathrm C}(G)\tau = \tau (\Bbb{C}_\bullet(G)).$
Explicitly, for $G = \theta^*,$ each $f_{i_0\ldots i_r}$ that is a morphism of $\underline{S}_n$-diagrams has two indices on the $\underline{S}_\bullet$-direction; $\theta^*$ reindexes these two and leaves the nerve indices unchanged.
The desired equation follows.
\end{proof}

Put $\tau_{\underline{S}_\bullet\underline{\mathcal M}} \coloneqq \mathrm{Tot}\,(\tau_{\bullet}).$
Two total complexes of the source and the target use the same sign convention: on a component of $\Bbb{C}_\bullet$- or $\underline{\mathrm C}_\bullet$-chain degree $q$, the total differential is the chain differential plus $(-1)^q$ times the alternating sum of the $\underline{S}_\bullet$-faces.

We show that the following diagram, which is the chain complex diagram in the category of chain complexes and of (\ref{formula:Dennistracediagram}), is commutative 
\begin{equation}\vcenter{\xymatrix@=1.5em{\Sigma \Bbb{C}_\bullet(\underline{\mathrm w\mathcal M})\ar[r]^{\Sigma\eta_{\underline{\mathcal M}}}\ar[d]_{\Sigma\tau_{\underline{\mathcal M}}} & \Bbb{C}_S(\underline{\mathrm wS}_\bullet\underline{\mathcal M})\ar[d]^{\tau_{\underline{S}_\bullet\underline{\mathcal M}}}\\
\Sigma\underline{\mathrm C}(\underline{\mathcal M})\ar[r]^{\Sigma\beta_{\underline{\mathcal M}}}&
\underline{\mathrm C}_S(\underline{\mathcal M}),}}\label{formula:chainDennis}\end{equation}
where $\Bbb{C}_S(\underline{\mathrm wS}_\bullet\underline{\mathcal M})\coloneqq \mathrm{Tot}\,\Bbb{C}_\bullet(\underline{\mathrm wS}_\bullet\underline{\mathcal M}).$
With the sign in $\Sigma\eta_{\underline{\mathcal M}}$ and $\Sigma\beta_{\underline{\mathcal M}}$, the horizontal maps $\eta_{\underline{\mathcal M}}$ and $\beta_{\underline{\mathcal M}}$ are chain maps because the $\underline{S}_\bullet$-direction differential vanishes. 
For any homogeneous element $x$ in $\Bbb{C}_q(\underline{\mathrm w\mathcal M})$, we have the desired
\[\tau_{\underline{S}_\bullet\underline{\mathcal M}}\circ \Sigma\eta_{\underline{\mathcal M}}(\Sigma x) = (-1)^q\tau_{1}(\Bbb{C}_\bullet(j_1)(x)) = (-1)^q\underline{\mathrm C}_\bullet(j_1)(\tau_{\underline{\mathcal M}}(x)) = (\Sigma\beta_{\underline{\mathcal M}})(\Sigma \tau_{\underline{\mathcal M}}(x)).\]

Finally, we pass to spaces using reduced chains and the Dold-Kan correspondence, with truncation after totalization. 
More explicitly, for a based CW space $W$ there is a natural identification
\[
[W,\mathscr H(E)]_*
\cong D(\mathrm{Ch}(k))\bigl(\widetilde C_*(W;k),E\bigr).
\]
The natural normalized-unnormalized and Eilenberg-Zilber comparisons give
\[\widetilde C_*(|\Bbb{N}_\bullet\underline{\mathrm w\mathcal M}|;k)\simeq \Bbb{C}_\bullet(\underline{\mathrm w\mathcal M})\text{ and }\widetilde C_*(|\Bbb{N}_\bullet(\underline{\mathrm wS}_\bullet\underline{\mathcal M})|;k)\simeq \Bbb{C}_S(\underline{\mathrm wS}_\bullet\underline{\mathcal M}).\]
Under these comparisons, $\tau_{\underline{\mathcal M}}$ and $\tau_{\underline{S}_\bullet\underline{\mathcal M}}$ define
$\tau_{\underline{\mathcal M}}^\mathrm t$ and $\tau_{S_\bullet\underline{\mathcal M}}^\mathrm t$, respectively, and the signed $\underline{S}_\bullet$-degree-one map $\Sigma\eta_{\underline{\mathcal M}}$ represents the delooping of $\eta_{\underline{\mathcal M}}^\mathrm t$. 
Its sign is the Koszul sign for placing the interval coordinate after a degree-$q$ chain, as in our totalization convention. 
Applying the suspension-loop adjunction to (\ref{formula:chainDennis}) therefore gives diagram (\ref{formula:Dennistracediagram}), commuting up to based homotopy. 
All these constructions are natural in the stated dg functors.
\end{proof}

\begin{definition}\label{def:dgDennistracemap}
For the dg Waldhausen category $\underline{\mathcal M}$ and the maximal Kan complex $\underline{\mathrm wS}_\bullet\underline{\mathcal M}$, the following map is called the \emph{Dennis trace map}
\[\xymatrix{\Omega\tau_{\underline{S}_\bullet\underline{\mathcal M}}^\mathrm t:\Omega|\Bbb{N}_\bullet \underline{\mathrm wS}_\bullet\underline{\mathcal M}| \ar[rr]&& 
\Omega\mathcal H_{S}(\underline{\mathcal M})}.\] 
\end{definition}

Taking homology, using Proposition~\ref{prop:comparison} and Lemma~\ref{lemma:delooping}, the right vertical map results in $K_\bullet(\underline{\mathcal M}) \to \mathrm{HH}_\bullet(\underline{\mathcal M}),$ denoted $\mathrm{tr}_{\underline{\mathcal M}}.$

\subsection{Derived Koszul equivalence and Dennis trace map}
We finish this paper with the compatibility between the derived Koszul equivalence and the Dennis trace map in the previous subsection. Retain the notations of Sections~\ref{sec:tau} to \ref{sec:explicitDennistraceconstruction}, including the unnormalized direct sum total complexes and the $k$-flatness hypothesis used in Lemma~\ref{lemma:delooping}.

\begin{theorem}\label{thm:compatibility}
With same assumptions in Proposition~\text{\rm \ref{prop:derivedKoszulequivalence}}, put 
\begin{equation}(\underline{\mathcal P},\underline{\mathcal Q}) = \begin{cases}(\underline{\mathcal M}_C(k),\underline{\mathcal M}_A(A)) & \text{without further condition};\\
(\underline{\mathcal M}_C(C),\underline{\mathcal M}_A(k)) & \left(\begin{split}&\small \text{if $C$ is dg and $A$ has a dg augmentation}\\[-1ex] 
&\small \text{$v:A\to k$ satisfying $v\circ \tau = 0$}\end{split}\right).
\end{cases}\nonumber\end{equation}
The twisting-tensor functors restrict to $L\coloneqq A\otimes^\tau - :\underline{\mathcal P}\rightleftarrows \underline{\mathcal Q}: C\otimes^\tau - \eqqcolon R.$ 
The functor $L$ induces the following diagram, which is commutative up to based homotopy and whose horizontal maps are weak equivalences:
\[\xymatrix@=1.5em{\Omega|\Bbb{N}_\bullet \underline{\mathrm wS}_\bullet\underline{\mathcal P}| \ar[d]_{\Omega\tau_{\underline{S}_\bullet\underline{\mathcal P}}^\mathrm t}\ar[rrr]_{\simeq}^{\Omega|\Bbb{N}_\bullet \underline{\mathrm wS}_\bullet(L)|} &&& \Omega|\Bbb{N}_\bullet \underline{\mathrm wS}_\bullet\underline{\mathcal Q}|\ar[d]^{\Omega\tau_{\underline{S}_\bullet\underline{\mathcal Q}}^\mathrm t}\\
\Omega\mathcal H_{S}(\underline{\mathcal P})\ar[rrr]_{\simeq}^{\Omega\mathcal H_{S}(L)} &&& \Omega\mathcal H_{S}(\underline{\mathcal Q}),}\]
where the vertical maps are the Dennis trace maps of Definition~\text{\rm{\ref{def:dgDennistracemap}}}.
\end{theorem}

\begin{proof}
We retain the notation in the proof of Theorem~\ref{thm:Dennistracediagram}.
We show the strict commutativity of the following diagram.
\begin{align}\vcenter{\xymatrix@=1.5em{\Bbb{C}_\bullet(\underline{\mathrm wS}_n\underline{\mathcal P})\ar[rr]^{\Bbb{C}_\bullet(\underline{S}_n L)}\ar[d]_{\tau_{\underline{S}_n}}&&\Bbb{C}_\bullet(\underline{\mathrm wS}_n\underline{\mathcal Q})\ar[d]^{\tau_{\underline{S}_n}}\\
\underline{\mathrm C}_\bullet(\underline{S}_n\underline{\mathcal P})\ar[rr]^{\underline{\mathrm C}_\bullet(\underline{S}_nL)} && \underline{\mathrm C}_\bullet(\underline{S}_n\underline{\mathcal Q}).}}\label{formula:KoszulDennischaincomplexescommutativity}\end{align}
Note that $L$ is a dg functor.
As $\underline{S}_nL$ is defined componentwise, the functor $\underline{S}_nL$ is also dg.
For any $p$-simplex $F$ in the thick nerve $\Bbb{N}_\bullet(\underline{S}_n\underline{\mathcal P})$, it gives higher homotopy $F(i_0,\ldots,i_r) = f_{i_0\ldots i_r}\in \underline{\mathrm wS}_n\underline{\mathcal P}(X_{i_0},X_{i_r})^{1-r}$ for $r\geq 1$. 
$(\underline{S}_nL)_*$ on the thick nerve gives the induced simplex $(\underline{S}_nL)_\bullet F$ that gives higher homotopy $\underline{S}_nL(f_{i_0\ldots i_r}).$ 
On the $(a,b)$-entry of the $\underline{S}_n$-diagram, this says that for $u\in A$, $x\in (X_{i_0})_{ab}$
\[((\underline{S}_nL)_*(F)(i_0,\ldots,i_r))_{ab}(u\otimes x) = (-1)^{(1-r)|u|}u\otimes ((f_{i_0\cdots i_r})_{ab}(x)).\]
Then, similar to the proof of Lemma~\ref{lemma:thicknerveofdgfunctor}, we have $\underline{\mathrm C}(\underline{S}_n L)\tau_{\underline{S}_n\underline{\mathcal P}} (F) = \tau_{\underline{S}_n\underline{\mathcal Q}}((\underline{S}_nL)_*(F)).$
After the reduced $k$-linearization, we have the desired commutativity in (\ref{formula:KoszulDennischaincomplexescommutativity}). 
For any $\theta:[m]\to [n],$ we have for any indices $(a,b)$ in the $\underline{S}_\bullet$-direction $(\theta^*\circ \underline{S}_nL (X))_{ab} = L(X_{\theta(a),\theta(b)}) = (\underline{S}_mL \circ \theta^*(X))_{ab}.$
The same identity holds on all homogeneous morphisms. 
Hence the dg functors equation $\theta^*\circ \underline{S}_nL = \underline{S}_mL\circ \theta^*$ and the commutative diagrams (\ref{formula:KoszulDennischaincomplexescommutativity}) for all $n$ assemble into one.

Taking direct sum totalization in the $\underline{S}_\bullet$- and $\underline{\mathrm C}_\bullet$-directions gives a strictly commutative square whose bottom map is $\underline{\mathrm C}_S(L):=\mathrm{Tot}\,\underline{\mathrm C}_\bullet(S_\bullet L):\underline{\mathrm C}_S(\underline{\mathcal P})\to\underline{\mathrm C}_S(\underline{\mathcal Q}).$
Both total complexes use the differential $d+(-1)^q b^S$ on $\underline{\mathrm C}_\bullet$-degree $q$. Finally, as in the proof of Theorem~\ref{thm:Dennistracediagram}, applying the natural Eilenberg-Zilber and reduced-chain comparisons for the bisimplicial sources, followed by the truncation, Dold-Kan, realization, and looping, gives the desired square of spaces, commutative up to based homotopy.

The upper equivalence is the one induced by $L$ in Proposition~\ref{prop:KKoszulequivalence}, translated in the thick-nerve model by Proposition~\ref{prop:comparison}.
It remains to justify the lower one.
We first use Proposition~\ref{prop:derivedKoszulequivalence} to show that the restricted dg functor $L:\underline{\mathcal P}\to\underline{\mathcal Q}$ is a quasi-equivalence (which means it induces quasi-isomorphisms on mapping complexes and being essentially surjective after passing to the homotopy category, cf. DK-equivalence in the proof of Proposition~\ref{prop:comparison}). 
Write $\eta:\id\to RL$ and $\varepsilon:LR\to\id$ for the unit and counit of the dg adjunction. 
Since all comodules are cofibrant and all modules are fibrant, the Quillen equivalence implies that
\[\eta_Y:Y\overset{\sim}{\to}RLY\text{ and }\varepsilon_V:LRV\overset{\sim}{\to}V\text{ for any $Y\in\underline{\mathcal P}$ and $V\in\underline{\mathcal Q}$}\]
are weak equivalences between cofibrant-fibrant objects in the ambient model structures and hence are homotopy equivalences. 
For $X,Y\in\underline{\mathcal P}$, using the dg adjunction, we have
\[\underline{\mathcal P}(X,Y)\xrightarrow{\eta_Y\circ-}
\underline{\mathcal P}(X,RLY)\cong\underline{\mathcal Q}(LX,LY).
\]
It is a quasi-isomorphism because postcomposition with a homotopy equivalence is a chain homotopy equivalence. 
The weak equivalence by the counit $\varepsilon_V$ shows that $H^0(L)$ is essentially surjective. 
This proves that $L$ is a quasi-equivalence.

In particular, $L$ is a derived Morita equivalence. With the $k$-flatness by $k$ being a field, Keller's Morita-invariance theorem \cite[Sections~4.6, 5.3, and Theorem~5.2(a)]{Keller} therefore implies that
\[
\underline{\mathrm C}_\bullet(L):\underline{\mathrm C}_\bullet(\underline{\mathcal P})\longrightarrow\underline{\mathrm C}_\bullet(\underline{\mathcal Q})
\]
is a quasi-isomorphism. The modified cyclic convention of Section~4.2 does not affect this conclusion, since its comparison with the ordinary convention is natural in dg functors.

The map $\Sigma\beta_{\underline{\mathcal M}}$ from Lemma~\ref{lemma:delooping} (See also (\ref{formula:betaM})) is natural in $L$, so the following square commutes strictly:
\[
\xymatrix@=1.5em{
\Sigma\underline{\mathrm C}_\bullet(\underline{\mathcal P})
\ar[rr]^{\Sigma\underline{\mathrm C}_\bullet(L)}_\simeq 
\ar[d]_{\Sigma\beta_{\underline{\mathcal P}}}^\simeq
&&\Sigma\underline{\mathrm C}_\bullet(\underline{\mathcal Q})
\ar[d]^{\Sigma\beta_{\underline{\mathcal Q}}}_\simeq\\
\underline{\mathrm C}_S(\underline{\mathcal P})\ar[rr]^{\underline{\mathrm C}_S(L)}&&\underline{\mathrm C}_S(\underline{\mathcal Q}).
}\] 
The bottom $\underline{\mathrm C}_S(L)$ is a quasi-isomorphism because of the other three maps being so and the two-out-of-three property. Consequently, $\Omega\mathcal H_S(L)$ is a weak equivalence.
\end{proof}

Taking homotopy groups of the diagram in the theorem now gives the stated group-level square in (\ref{formula:KKoszulhomotopygroup}). 
The naturality of $\Sigma\beta_{\underline{\mathcal M}}$ identifies its lower horizontal map with $\mathrm{HH}_n(L)$.

\appendix

\section{$\mathrm{Ch}(k)$-model category}\label{appendix:chk-modelcategory}
With the fixed ground field $k$, let $\mathrm{Ch}(k)$ denote the category of unbounded cochain complexes. 
It is a closed symmetric monoidal category with the tensor product $K\otimes_k L$ and the unit $k[0]$. 
It has the projective model structure with weak equivalences being quasi-isomorphisms, fibrations being degreewise surjective, and cofibrations being maps having left lifting property with respect to trivial fibrations.
Let $\underline{\mathcal Q}$ be a $\mathrm{Ch}(k)$-enriched category, i.e., we have $\underline{\mathcal{Q}}(X,Y) \in \mathrm{Ch}(k)$ together with the following composition and unit:
\[\underline{\mathcal Q}(Y,Z) \otimes_k \underline{\mathcal Q}(X,Y)\to \underline{\mathcal Q}(X,Z)\text{ and }k\to \underline{\mathcal Q}(X,X)\text{ for any objects }X,Y,Z.\]
Put $\mathcal Q \coloneqq Z^0\underline{\mathcal Q}.$ 
The dg category $\underline{\mathcal Q}$ is \emph{tensored and cotensored} if there are the following bifunctors
\begin{align*}& -\otimes- :\mathrm{Ch}(k) \times \mathcal Q\to \mathcal Q;\quad (K,X)\mapsto K\otimes X\text{ and}\\
&(-)^{(-)}:\mathrm{Ch}(k)^{op}\times \mathcal Q \to \mathcal Q;\quad (K,X)\mapsto X^K
\end{align*}
such that the enriched functoriality of these two actions, the associativity pentagon, left- and right-unit coherence holds, the cotensor coherence, and the following isomorphisms -- the tensor-cotensor adjunction hold
\begin{align}\mathcal Q(K\otimes X,Y)\cong \mathrm{Ch}(k)(K,\underline{\mathcal Q}(X,Y))\cong \mathcal Q(X,Y^K).\label{formula:tensor-cotensor}\end{align}
This tensor-cotensor adjunction has the following enriched version using the internal hom of $\mathrm{Ch}(k)$ ($\underline{\mathrm{Ch}}(k)$ denotes the internal hom), and these two version are equivalent:
\[\underline{\mathcal Q}(K\otimes X,Y)\cong \underline{\mathrm{Ch}}(k)(K,\underline{\mathcal Q}(X,Y))\cong \underline{\mathcal Q}(X,Y^K).\]

\begin{definition}\label{def:enrichedmodelcatgegory}
The $\mathrm{Ch}(k)$-enriched, tensored, and cotensored category $\underline{\mathcal Q}$ is a \emph{$\mathrm{Ch}(k)$-model category} if the following conditions holds:
\begin{itemize}[noitemsep,leftmargin=*,nosep]
    \item The underlying ordinary category $\mathcal Q$ is a closed model category;
    \item (SM7: compatibility of the model structure and the dg structure) For each cofibration $i:M\to N$ and fibration $p:X\to Y$ in $\mathcal Q,$ the following map is a fibration in $\mathrm{Ch}(k)$
    \begin{align}\Phi:\underline{\mathcal Q}(N,X)\to \underline{\mathcal Q}(M,X)\times_{\underline{\mathcal Q}(M,Y)}\underline{\mathcal Q}(N,Y)\label{formula:SM7compofmodelanddg}\end{align}
    and is a trivial fibration if either $i$ or $p$ is trivial, where the pullback is taken in $\operatorname{Ch}(k)$, and $\Phi(h)=(h\circ i,\;p\circ h).$
\end{itemize}
\end{definition}

\begin{lemma}[\text{\cite[Proposition~1.36]{Barwick}}]\label{lemma:cofib-fibreplacementcomputesderivedhom}
Let $\underline{\mathcal Q}$ be a $\mathrm{Ch}(k)$-model category and $Q,R:\underline{\mathcal Q} \to \underline{\mathcal Q}$ respectively cofibrant and fibrant replacement functors.
Then in $\mathrm{Ho}\,\mathrm{Ch}(k)$ we have
\[\mathbf{R}\underline{\mathcal Q}(X,Y)\cong \underline{\mathcal Q}(QX,RY).\]
\end{lemma}

\subsection{Example: the category $\underline{\mathrm M}_A$ of dg $A$-modules}\label{example:dgAmoduleChkmodelcategory}
In this subsection, we check that the dg category $\underline{\mathrm{M}}_A$ of dg $A$-modules is a $\mathrm{Ch}(k)$-model category.
The tensor and cotensor respectively write $K\otimes M = K \otimes_k M$ and $M^K = \underline{\mathrm{Ch}}(k)(K,M).$
It is straightforward to check $\underline{\mathrm{M}}_A$ is $\mathrm{Ch}(k)$-enriched tensored cotensored category.

By Proposition~\ref{prop:projectivemodelstructureformodules}, it remains to check the compatibility condition: the map $\Phi$ in the definition is a fibration (degreewise surjective) and trivial fibration (quasi-isomorphic and degreewise surjective) given $i$ or $p$ being moreover trivial. 
Suppose any degree $r$ element $(f,g)$ in the codomain of $\Phi$, where $f:M\to X$ and $g:N\to Y$ are homogeneous maps satisfying $pf=gi.$ 
As $i$ is a cofibration, $N/M$ is cofibrant, hence graded-projective, and the sequence splits as graded $A^\#$-modules.
Choose the retraction $r:N\to M$ and section \(s:P\to N\) from one graded direct-sum decomposition so that $ri=\mathrm{id}_M,$ $\pi s=\mathrm{id}_P,$ $ir+s\pi=\mathrm{id}_N.$
As $N/M$ is graded projective and $X\to Y$ is a fibration, there is a degree $r$ map $g'$ lifts $g s$ such that $pg'=gs$.
Let $h$ be the degree $r$ map $fr + g'\pi:N\to X.$
This is the preimage of $(f,g)$.
Indeed, one may straightforwardly check $hi = f$ and $ph = g,$ as desired.

We claim $\ker \Phi \cong \underline{\mathrm{M}}_A(N/M,\ker p).$
For a degree $|f|$ homogeneous map $f$ in the domain of $\Phi$, it is in $\ker \Phi$ means $fi = 0 = pf.$
The equation $fi=0$ means $f$ factors as $f' \circ \pi: N\to N/M\to X.$
As $p f' \pi = 0$ and $\pi$ is surjective, we have $pf' = 0$ and hence $f'$ uniquely factors as $\iota\tilde f: N/M\to \ker  p \hookrightarrow X$.
This defines $\ker \Phi\to \underline{\mathrm{M}}_A(N/M,\ker p);\,\,f\mapsto \tilde f$ and the inverse map is given by $g\mapsto \iota g\pi$. 
The desired isomorphism follows.

We are to use the fact that $\underline{\mathrm M}_A(P,T)$ is acyclic if $P$ is cofibrant and $T$ is acyclic \cite[Section~3.2]{Keller}.
If $p$ is trivial, then $\ker p$ is acyclic. 
As $N/M$ is cofibrant, $\underline{\mathrm{M}}_A(N/M,\ker p)$ is acyclic, as desired.
For the case $i$ being trivial, we are to find a degree $-1$ map $H$ so that $\delta H + H\delta = \id_{\underline{\mathrm{M}}_A(N/M,\ker p)}$ with the differential of the mapping complex $\delta(f)= d_{\ker p}f - (-1)^{|f|} fd_{N/M}$.
The condition implies the acyclicity of the complex $\underline{\mathrm{M}}_A(N/M,N/M)$ and hence there exists a degree $-1$ map $s\in \underline{\mathrm{M}}_A(N/M,N/M)$ such that $d_{N/M}s + sd_{N/M} = \id_{N/M}$.
It is routine to check the following $H$ is the desired one
\[H:\underline{\mathrm{M}}_A(N/M,\ker p) \to \underline{\mathrm{M}}_A(N/M,\ker p);\quad f\mapsto (-1)^{|f|}fs.\]

\subsection{Example: the category $\underline{\mathrm{CM}}_C$ of cdg $C$-comodules}\label{example:dgCcomoduleChkmodelcategory}
We check that the dg category $\underline{\mathrm{CM}}_C$ of cdg $C$-comodules is a $\mathrm{Ch}(k)$-model category.
To define the cotensor, for objects $K$ in $\mathrm{Ch}(k)$ and $M$ in $\underline{\mathrm{CM}}_C$, consider $\underline{\mathrm{Ch}}(k)^\#(K,M)$ with $\underline{\mathrm{Ch}}(k)^\#(K,M)^n = \prod_p k\text{-Vct}(K^p,M^{p+n})$, where $k\text{-Vct}$ denotes the category of $k$-vector spaces. 
There is a canonical injection of vector spaces
\begin{align}j:C\otimes \underline{\mathrm{Ch}}(k)^\#(K,M) \to \underline{\mathrm{Ch}}(k)^\#(K,C\otimes_k M);\quad j(c\otimes f)(l) = c\otimes f(l).\label{formula:Rationalsubspace}\end{align}
Put $\mathrm{Rat}_k(K,M) \coloneqq \{f\in \underline{\mathrm{Ch}}(k)^\#(K,M) \mid \rho_Mf\in \mathrm{Im}\,j\}$, the subspace of rational homogeneous morphisms, where $\rho_M$ denotes the comodule map of $M$. 
Here $f$ satisfies $\rho_Mf\in \mathrm{Im}\,j$ if there are homogeneous $c_i\in C$ and homogeneous $f_i:K\to M$ such that $\rho_M(f(l)) = \sum_i c_i\otimes f_i(l)$ for any $l\in K$ is a finite sum.

\begin{lemma}\label{lemma:comoduleChk-modelstructure}
The category $\underline{\mathrm{CM}}_C$ is a tensored and cotensored category with the tensor and cotensor operations respectively being $K\otimes_k M$ and $M^K\coloneqq   \mathrm{Rat}_{k}(K,M)$ with differentials 
\begin{align*}d_{\otimes}(l\otimes m) & = d_K l\otimes m + (-1)^{|l|}l\otimes d_M m\text{ and }\\
d_{\mathrm{R}}(f) & = d_Mf - (-1)^{|f|}fd_K.\end{align*}
\end{lemma}

\begin{proof}
For $K\otimes_k M$, the comodule structure writes \[\rho_\otimes(l\otimes m) = \sum(-1)^{|l||m_{(-1)}|}m_{(-1)}\otimes l\otimes m_{(0)}.\]
The graded symmetry commutes with the tensor differentials, so $\rho_\otimes$ is compatible with $d_\otimes$. 
By the linearity of the curvature $h$, we compute $d_{\otimes }^2(l\otimes m) = l\otimes d_M^2(m) = (h\otimes\id_{K\otimes M})\rho_\otimes(l\otimes m).$ Here only the sign $1$ survives the application of $h$ because $h(c)$ is nonzero only if $|c| = -2$. 
The proof of the isomorphisms corresponding to the tensor $K\otimes_k M$, including: the left side of (\ref{formula:tensor-cotensor}), the associativity pentagon, the left- and right-unit coherence, are similar to those for the module over a dg algebra, and hence omitted.

For $\mathrm{Rat}_k(K,M),$ it has the coaction $\rho_\mathrm{R}$ such that $\rho_\mathrm{R}(f) = \sum_i c_i \otimes f_i$ so that $j(\rho_\mathrm{R}(f))(l) = \rho_M(f(l)).$ 
We first show the map $j$ in (\ref{formula:Rationalsubspace}) is injective.
Indeed, for an element $\sum_i c_i\otimes f_i$ in the source with $c_i$ being $k$-linearly independent, if its image under $j$ vanishes, then applying linear functionals dual to $c_i$ shows all $f_i$ are zero. 
This injective implies $j^{-1}(\rho_Mf)$ and hence $\rho_\mathrm{R}(f)$ are independent of finite sum expressions.

We show that $f_i$ above is in $\mathrm{Rat}_k(K,M).$
To compute $\rho_M\circ f_j$, we use the following equation for any $l\in K$ given by the coassociativity $(\Delta\otimes \id_M)\rho_M = (\id_C \otimes \rho_M) \rho_M$ 
\[\textstyle \sum_{i=1}^r \Delta(c_i)\otimes f_i(l) = \sum_{i=1}^r c_i\otimes \rho_M(f_i(l)).\]
Choose the homogeneous coefficients $c_i$ to be $k$-linearly independent, and choose homogeneous linear functionals $\lambda_j:C\to k$ satisfying $\lambda_j(c_i)=\delta_{ji}.$
Put $e_{ji} \coloneqq (\lambda_j\otimes \id_C)\Delta(c_i).$
Apply $\lambda_j \otimes \id_C \otimes \id_M$ to the above equation and we have
\[\textstyle \sum_{i=1}^r e_{ji} \otimes f_i(l) = \rho_M\circ f_j(l),\]
and hence $f_j\in \mathrm{Rat}_k(K,M).$
Coassociativity and counitality of $\rho_\mathrm{R}$ follow by evaluation from those of $\rho_M$, using analogous injectivity of $j$ with $C$ replaced with $C\otimes C$ and with $k.$

For the operator $Df = d_M f - (-1)^{|f|}fd_K$, we need to show the stability, i.e., $Df$ is still in $\mathrm{Rat}_k(K,M)$, the compatibility with the coaction, and the curvature identity.
A straightforward computation shows \[\rho_M(Df) = \sum_i d_C(c_i)\otimes f_i + \sum_i (-1)^{|c_i|}c_i\otimes Df_i.\]
The finite expression shows the stability, which enables us to set $d_\mathrm{R} = D|_{\mathrm{Rat}_k(K,M)}.$
The coaction $\rho_\mathrm{R}$ is compatible with $d_\mathrm{R}$.
Moreover, we have $(d_\mathrm{R}^2f)(l) = d_M^2 f(l) - fd_K^2 = (h\otimes \id_M)\rho_M(f(l)) = \sum_i h(c_i)f_i(l),$ and hence $d_\mathrm{R}^2f = (h\otimes \id_\mathrm{R})\rho_\mathrm{R}(f)$ as desired.

We then show the right side isomorphism in  (\ref{formula:tensor-cotensor}).
If this isomorphism holds, the cotensor coherence follows straightforwardly.
Consider the map 
\[\Theta:\mathrm{Ch}(k)(K,\underline{\mathrm{CM}}_C(M,N))\to \mathrm{CM}_C(M,N^K);\quad \Theta(F)(m)(l) = (-1)^{|m||l|}F(l)(m)\]
We first need to show that $\Theta(F)(m)$ is rational.
Given $\rho_M(m) = \sum m_{(-1)}\otimes m_{(0)}$ (the sum is finite), as $F(l)$ is $C$-colinear, we have
\begin{align*}
\rho_N(\Theta(F))(m)(l) & = (-1)^{|m||l|}\rho_N(F(l)(m))\\
& = \sum (-1)^{|m||l|+|l||m_{(-1)}|}m_{(-1)}\otimes F(l)(m_{(0)})\\
& = \sum (-1)^{|m||l|+|l||m_{(-1)}|+|l||m_{(0)}|}m_{(-1)}\otimes \Theta(F)(m_{(0)})(l)\\
& = \sum m_{(-1)}\otimes \Theta(F)(m_{(0)})(l).
\end{align*}
As the sum above is finite, $\Theta(F)(m)$ is rational.
Because $\Theta(F)(m_{(0)})$ is also rational, using the definition of the rational coaction gives $\rho_{N^K}(\Theta(F)(m)) = \sum m_{(-1)}\otimes \Theta(F)(m_{(0)}),$ which means $\Theta(F)$ is a graded comodule map.
One may straightforwardly show the following map gives an inverse
\[\Psi:\mathrm{CM}_C(M,N^K) \to \mathrm{Ch}(k)(K,\underline{\mathrm{CM}}_{C}(M,N));\quad \Psi(G)(l)(m) = (-1)^{|m||l|}G(m)(l).\]
Indeed, by the colinearity identity for $G$:
\[\rho_N(\Psi(G)(l)(m)) = \sum (-1)^{|l||m_{(-1)}|}m_{(-1)}\otimes \Psi(G)(l)(m_{(0)}),\]
implying $\Psi(G)(l)\in \underline{\mathrm{CM}}_C(M,N)^{|l|}$.

Finally, we check the enriched functoriality.
For $u:K\to K'$ and $v:M\to M'$, we need the existence of the following morphism of complexes 
\begin{equation*}
\begin{array}{c}
\mathrm{Ch}(k)(K,K')\otimes \underline{\mathrm{CM}}_{C}(M,M')\to \underline{\mathrm{CM}}_{C}(K\otimes_k M,K'\otimes_k M'),\\
\mathrm{Ch}(k)(K',K)\otimes\underline{\mathrm{CM}}_C(M,M')\to \underline{\mathrm{CM}}_C(M^K,(M')^{K'}).
\end{array}    
\end{equation*}
that preserves enriched identities and composition.
A similar morphism for the cotensor should also hold.
We claim the following map is the desired one:
\[(u,v)\mapsto u\otimes v;\,\,\text{with }u\otimes v(l\otimes m) = (-1)^{(|v||l|)}u(l)\otimes v(m).\]
This map is $C$-colinear because both $\rho_{K'\otimes M'}(u\otimes v)(l\otimes m) = (\id_C\otimes (u\otimes v))\rho_{K\otimes M}(l\otimes m)$ equal 
\[\textstyle\sum (-1)^{|v||l| + (|l|+|u|+|v|)|m_{(-1)}|}m_{(-1)}\otimes u(l)\otimes m_{(0)}.\]
The compatibility with the differential, the preservation of identity, and compositions are straightforward.
As for the cotensor, we replace $K\to K'$ with $K'\to K$, still denoted $u$.
Given $v$ and $u$, the following map is the desired one:
\[\text{for }f\in M^K,\,\,f\mapsto R(u,v)f\coloneqq (-1)^{|u|(|v|+|f|)}vfu.\]
The $C$-colinearity follows from that of $v$.
We need to show the image is still rational.
If $\rho_{M^K}(f) = \sum_{i=1}^{r}c_i\otimes f_i$, by the $C$-colinearity, we have \begin{align*}
\textstyle \rho_{M'}\circ R(u,v)(f)(l) & \textstyle = \sum_{i=1}^r (-1)^{|u|(|v|+|f|)+|v||c_i|}c_i\otimes vf_iu(l)\\
& \textstyle  = \sum_{i=1}^{r}(-1)^{(|v|+|u|)|c_i|}c_i\otimes R(u,v)(f_i)(l)\end{align*} and hence the image is rational.
By this equation, the desired graded-colinearity identity also follows.
The compatibility with the differential, the preservation of identity and compositions are routine to check.
\end{proof}

By Proposition~\ref{prop:modelstructureforcdgcomodules}, it suffices to check the compatibility condition: the map $\Phi$ is degreewise surjective, or is degreewise surjective and a weak equivalence if $i$ or $p$ is moreover trivial.
The map $i$ is a degreewise injection.
The map $p$ gives an exact sequence $\ker p\overset{\iota}{\hookrightarrow} X\overset{p}{\twoheadrightarrow} Y$ with $\ker p$ being graded injective.
This sequence of graded comodules is split and let $s: Y\to X$ be the section of $p$.
For $(f,g)$ in the codomain of $\Phi$ such that $pf = gi$, the map $f-sgi$ has image in $\ker p$.
As $\ker p$ is graded injective and $M\to N$ is a cofibration, this map extend to $f':N\to \ker p$.
One may check that the map $\iota f'+sg$ is the preimage of $(f,g)$.

As in Example~\ref{example:dgAmoduleChkmodelcategory}, one may show the isomorphism $\ker \Phi \cong \underline{\mathrm{CM}}_C(N/M,\ker p)$ of complexes.
If $i$ is a trivial cofibration, then its pushout $0\to N/M$ along $M \to 0$ is a trivial cofibration, so $N/M$ is coacyclic. 
If $p$ is a trivial fibration, then its pullback $\ker p\to 0$ along $0\to Y$ is a trivial fibration, so $\ker p$ is coacyclic.
Suppose first that $i$ is trivial. 
Since $N/M$ is coacyclic and $\ker p$ is graded-injective, \cite[Theorem~4.4(a)]{Posit}implies that $\underline{\mathrm{CM}}_C(N/M,\ker p)$ is acyclic.
Suppose instead that $p$ is trivial. Applying the same theorem to $\ker p$ shows that $\underline{\mathrm{CM}}_C(\ker p,\ker p))$ is acyclic. Hence there is a degree-$(-1)$ graded comodule map $s:\ker p\to \ker p$ satisfying
$d_{\ker p}s+sd_{\ker p}=\mathrm{id}_{\ker p}.$
The operator $H(f)=sf$ on $\underline{\mathrm{CM}}_C(N/M,\ker p))$ satisfies
\[(\partial H+H\partial)(f)=(d_{\ker p}s+sd_{\ker p})f=f.\]
Thus $\ker\Phi$ is acyclic whenever either $i$ or $p$ is trivial. 
Since $\Phi$ is degreewise surjective, it is then a quasi-isomorphism. This verifies SM7.

\begin{remark}
The transposition formulas above extend to homogeneous maps of arbitrary degree and give natural isomorphisms of complexes
\[\underline{\mathrm{CM}}_C(K\otimes M,N)\cong\underline{\operatorname{Ch}}(k)\bigl(K,\underline{\mathrm{CM}}_C(M,N)\bigr)\cong\underline{\mathrm{CM}}_C(M,N^K).\]
Indeed, direct substitution into the mapping-complex differentials gives
$\partial\Theta(F)=\Theta(\partial F).$
Taking $Z^0$ recovers the ordinary adjunctions in (\ref{formula:tensor-cotensor}). The cotensor associativity and unit isomorphisms are induced by the corresponding tensor isomorphisms; their coherence follows by uniqueness of adjoints.
\end{remark}

\begin{remark}[The weak path functors]\label{remark:weakpathobj}
We will also need the categories $\underline{\mathrm M}_A$ and $\underline{\mathrm CM}_C$ to admit weak path functors, that is the following isomorphisms holds for any simplicial set $K$
\[\mathrm{M}_A(X,Y^K) \cong \mathbf{sSet}(K,(\mathrm{M}_A)_\Delta (X,Y))\text{ and }\mathrm{CM}_C(X,Y^K) \cong \mathbf{sSet}(K,(\mathrm{CM}_C)_\Delta (X,Y)),\]
where the lower index $\Delta$ means taking corresponding simplicial category.
By Lemma~\ref{lemma:comoduleChk-modelstructure}, for $C_K\coloneqq \mathrm{Norm}_\bullet(k[K])$ as a cochain complex and $Y^K\coloneqq \mathrm{Rat}_k(C_K,Y),$ we have 
\begin{align*}
\mathrm{CM}_C(X,Y^K) & \cong \mathrm{Ch}(k)(C_K,\underline{\mathrm{CM}}_C(X,Y))\\
& \cong \mathbf{sSet}(K,\Gamma(\underline{\mathrm{CM}}_C(X,Y))^{\leq 0}),
\end{align*}
where $\Gamma$ is the functor given by the Dold-Kan correspondence and the second isomorphism follows from this.
Here for a complex $C$ with cohomology grading, the upper index $\leq 0$ means taking truncation so that $(C^{\leq 0})^{0} = \ker(d:C^0\to C^1).$
The proof for the left one is similar and omitted.
\end{remark}

\section{Reedy model categories}\label{sec:Reedymodelcategories}
We recall some facts about Reedy categories and corresponding defined model structures.
\begin{definition}\label{def:Reedymodelstructure}
A \emph{Reedy category} is a small category $\mathcal{R}$ with a degree function $\mathrm{deg}:\mathrm{obj}(\mathcal R)\to \Bbb{N}$ and wide subcategories $\mathcal R^+$ and $\mathcal R^-$ such that 
\begin{itemize}[noitemsep,leftmargin=*,nosep]
\item every nonidentity morphism in $\mathcal R^+$ strictly raises degree, every nonidentity morphism in $\mathcal R^-$ strictly lowers degree; 
\item each map of $\mathcal R$ factors uniquely as $\alpha = \alpha^+\circ \alpha^-$ with $\alpha^+\in \mathcal R^+$ and $\alpha^-\in \mathcal R^-$.
\end{itemize}
Here by ``wide subcategories'', we mean the subcategories contain all objects and all identity morphisms.
\end{definition} 

\begin{proposition}
Given a model category $\mathcal Q$, the category $\mathcal Q^{\mathcal R} = \mathrm{Fun}(\mathcal R,\mathcal Q)$ has a model structure, named the \emph{Reedy model structure} described as below.
Let $f:X\to Y$ be a natural transformation as a morphism in $\mathcal Q^{\mathcal R}$.
\begin{itemize}[noitemsep,leftmargin=*,nosep]
    \item $f$ is a weak equivalence if $f_r:X_r\to Y_r$ is a weak equivalence for each object $r$ of $\mathcal R$;
    \item $f$ is a cofibration if $X_r\coprod_{L_rX}L_rY \to Y_r$ is a cofibration for all $r$. Here, $L_rX$ is the \emph{latching object} of $X$ at $r$ defined by the following colimit and $L_rY$ is defined similarly
    \[\textstyle L_rX = \mathrm{colim}_{\substack{s\to r\in \mathcal R^+\\s\to r\neq \id}}X_s.\]
    If there is no nonidentity morphism $s\to r$ in $\mathcal R^+$, then $L_rX$ is the initial object of $\mathcal Q$.
    \item $f$ is a fibration if $X_r\to Y_r\times_{M_rY}M_rX$ is a fibration for all $r$. Here, $M_rX$ is the \emph{matching object} of $X$ at $r$ defined by the following limit and $M_rY$ is defined similarly
    \[\textstyle M_rX = \lim_{\substack{r\to s\in \mathcal R^-\\r\to s\neq \id}}X_s.\]
    If there is no nonidentity morphism $r\to s$ in $\mathcal R^-$, then $M_rX$ is the terminal object of $\mathcal Q$.
\end{itemize}
\end{proposition}

\begin{theorem}[\text{Angeltveit \cite[Theorem~4.7]{Angel}}]\label{thm:ChkReedy}
For a $\mathrm{Ch}(k)$-Reedy category $\underline{\mathcal A}$ and a $\mathrm{Ch}(k)$-model category $\underline{\mathcal Q}$, the enriched functor category $\underline{\operatorname{Fun}}_{\mathrm{dg}}(\underline{\mathcal A},\underline{\mathcal Q})$ equipped with its Reedy model structure, is a $\mathrm{Ch}(k)$-model category. \end{theorem}

\noindent We apply the theorem to the dg linearization $k[\mathcal R]$, which has the same objects as $\mathcal R$, and
\[\textstyle k[\mathcal R](r,s) = \bigoplus_{\alpha\in \mathcal R(r,s)}k [\alpha]\cong \bigoplus_{\alpha\in \mathcal R(r,s)}k [0].\]
For this, we need to check the conditions in \cite[Definition~4.1]{Angel}. 
Indeed, the unique factorization composition in Definition~\ref{def:Reedymodelstructure} now is the isomorphism $k[0]\otimes_k k[0]\cong k[0].$
Moreover, the morphism complexes $k[\mathcal R]^+(r,s)$ and $k[\mathcal R^-](r,s)$ are direct sums of copies of $k[0]$ and hence cofibrant.
Let $\underline{\mathcal Q}^{\mathcal R}$ denote the obtained enriched functor category.

Each functor in $\mathrm{Fun}(\mathcal R,\mathcal Q)$ linearly extends to a dg functor $k[\mathcal{R}]\to\underline{\mathcal Q}$.
Ordinary transformation correspond to closed degree-$0$ dg transformation.
The mapping complex $\underline{\mathcal Q}^{\mathcal R}(X,Y)$ is the strict end whose degree-$n$ component consisting of $(f_r)\in \prod_{r\in \mathcal R}\underline{\mathcal Q}(X_r,Y_r)^n$ giving commutative diagrams.
For the underlying ordinary $\mathcal Q$ of $\underline{\mathcal Q}$, the enriched latching and matching objects in \cite[Definition~4.2]{Angel} agree with the ordinary ones. 
Hence the existed Reedy model structure on $\mathcal Q^{\mathcal R}$ is compatible with the dg enrichement.
For $K\in\mathrm{Ch}(k)$, tensors and cotensors are defined objectwise with $(K\otimes_k X)_r = K\otimes_k X_r$ and $(X^K)_r = X_r^K.$

\begin{remark}\label{remark:stable}
For $\underline{\mathcal Q}$ as in the theorem, the underlying category $\mathcal Q^\mathcal R$ is pointed.
Consider \[\Sigma:\underline{\mathcal Q}^{\underline{\mathcal R}}\rightleftarrows \underline{\mathcal Q}^{\underline{\mathcal R}}: \Omega\text{ with }X\mapsto \Sigma X = X\otimes k[1]\text{ and }X\mapsto \Omega X = X^{k[1]}.\]
By $\underline{\mathcal Q}^{\mathcal R}(\Sigma X,Y)\cong \underline{\mathcal Q}^{\mathcal R}(X,\Omega Y),$ the degree-$0$ cycles give the ordinary adjunction $\mathcal Q^\mathcal R\rightleftarrows \mathcal Q^\mathcal R.$
Since $k[1]$ is cofibrant, this is a Quillen adjunction. 
Moreover, $\Omega\cong k[-1]\otimes-$ is also left Quillen because $k[-1]$ is cofibrant. 
Thus the mutually inverse shifts $\Sigma$ and $\Omega$ induce mutually inverse derived functors, and the adjunction is a Quillen equivalence.
One may identify the derived functors of $\Sigma$ and $\Omega$ with suspension and loops.
Indeed, for $I = k[0]$, the cofibration $I\hookrightarrow E = \mathrm{cone}(\id_I),$ and a Reedy cofibrant diagram $X\in\mathcal Q^{\mathcal R}$, as $E\otimes X$ is contractible, the following isomorphism shows the suspension is isomorphic to $\mathbf{L}\Sigma$
\[\textstyle (E\otimes X) \coprod_{I\otimes X}0\cong (E\coprod_I 0)\otimes X \cong k[1]\otimes X = \Sigma X.\]
The derived functor satisfies $\mathbf{L}\Sigma\dashv \mathbf{R}\Omega.$ 
Hence $\mathbf{R}\Omega$ is the loop functor and $\mathcal Q^\mathcal R$ is stable.
\end{remark}

\begin{example}
For the category $\Delta$ those objects are $[n]$ for $n\geq 0$, equip it with the degree function $\deg([n]) = n.$
The class $\Delta^+$ consists of injections $[m]\hookrightarrow [n]$ with $m\leq n,$ e.g. cofaces and $\Delta^-$ consists of surjection is $[n]\twoheadrightarrow [m]$, e.g., codegeneracies.
The category of cosimplicial simplicial sets inherits the Reedy model structure from
the Kan-Quillen model structure on $\mathbf{sSet}.$
The latching object at $[n]$ is $\mathrm{colim}_{\substack{[m]\hookrightarrow[n]\\m<n}}A^m$.
\end{example}

\begin{example}\label{example:Sconstruction}
Recall the category $[n]$ such that the map between $i$ and $j$ is $i\to j$ for all $0\leq i\leq j\leq n$. 
For each $n$ and the category $\mathcal I_n\coloneqq [n]\times [n]$, the degree function is defined as $d(i,j) = n-i+j \in \{0,\ldots,2n\}.$ 
With $\mathcal I_n^-$ consisting of $(i,j)\to(i',j)$ with $i\leq i'$ and $\mathcal I_n^+$ consisting of $(i,j)\to(i,j')$ with $j\leq j'$.
Then each map $(i,j)\to (i',j')$ factors into $(i,j)\to (i',j) \to (i',j'),$ where the first nonidentity map lowers degree and the second raises degree.
Hence $\mathcal I_n$ is a Reedy category.
We may apply Theorem~\ref{thm:ChkReedy} to $k[\mathcal I_n]$ and a $\mathrm{Ch}(k)$-model category $\underline{\mathcal Q}$ to get $\underline{\mathcal Q}^{\mathcal I_n}$.
If $\mathcal Q$ is a pointed model category,
then for the latching object $L_{(i,j)}X$ and the matching object $M_{(i,j)}X$ of the functor $X:k[\mathcal I_n]\to \underline{\mathcal Q}$, as the nonempty latching (resp. matching) category has the terminal object $(i,j-1)\to (i,j)$ (resp. initial object $(i,j)\to (i+1,j)$), we have 
\[L_{(i,j)}X = \begin{cases}X_{i,j-1}& j>0;\\
0& j=0,\end{cases}\text{ and }M_{(i,j)}X = \begin{cases}X_{i+1,j}& i<n;\\
0& i=n.\end{cases}\]
This also implies that $X$ is Reedy cofibrant if and only if all $X_{i0}$ are cofibrant and all $X_{i,j-1}\to X_{ij}$ are cofibrations, and that $X$ is Reedy fibrant if and only if all $X_{nj}$ are fibrant and all $X_{ij}\to X_{i+1,j}$ are fibrations.
\end{example}

\section{The cyclic nerve of a dg category}\label{appendix:cyclicnerve}
Suppose $\underline{\mathcal P}$ to be a small $k$-linear dg category with a strict zero object and finite biproduct.
We use cohomological grading, so that $|df|=|f|+1$, and write $gf=g\circ f$.
Then the dg cyclic nerve in degree $n$ is
\[\underline{\mathrm{CN}}_n(\underline{\mathcal P}) \coloneqq  \textstyle\bigoplus_{\substack{X_0,\ldots,X_n\text{ are}\\\text{objects in $\underline{\mathcal P}$}}}\underline{\mathcal P}(X_n,X_0)\otimes_k \underline{\mathcal P}(X_{n-1},X_n) \otimes_k\cdots \otimes_{k}\underline{\mathcal P}(X_{0},X_{1}).\]
For $n=0$, this means $\underline{\mathcal P}_0(\underline{\mathcal P})=\bigoplus_{X_0}\underline{\mathrm{CN}}(X_0,X_0)$.
We usually write $f_n\otimes\cdots\otimes f_0$ as the element, we have $ f_i\in\underline{\mathcal P}(X_i,X_{i+1})$ for $0\leq i<n,$ and $f_n\in\underline{\mathcal P}(X_n,X_0).$

\begin{lemma}\label{lemma:dgcyclicnerve}
The internal differential is written as
\begin{align}\textstyle d_{\mathrm{int}}(f_n\otimes \cdots \otimes f_0) = \sum_{j=0}^{n}(-1)^{|f_n|+\cdots+|f_{n-(j-1)}|}f_n\otimes\cdots\otimes df_{n-j}\otimes \cdots\otimes f_0.\label{formula:internaldifferential}\end{align}
The empty sum in the exponent is treated to be zero.
On $f_n\otimes \ldots \otimes f_0,$ with 
\begin{itemize}[noitemsep,leftmargin=*,nosep]
    \item degeneracies $\bar s_j^\mathrm H$ for $j=0,\ldots,n$ given by 
    \[\bar s_j^\mathrm H(f_n\otimes \cdots \otimes f_0) = f_n\otimes \cdots \otimes f_{j}\otimes \id_{X_j}\otimes f_{j-1} \otimes \cdots \otimes f_{0};\]
    \item the cyclic operator $t_n$ as below
    \[t_n^\mathrm H(f_n\otimes\cdots\otimes f_0) = (-1)^{|f_0|S}f_0\otimes f_n \otimes\cdots\otimes f_{1}.\]
    with $S=|f_1|+\cdots+|f_{n}|.$
    \item faces $\bar d_j^\mathrm H$ for $j=0,\ldots,n$ by composition
    \[\bar d_j^\mathrm H(f_n\otimes \cdots \otimes f_{0}) =\begin{cases}f_n\otimes\cdots\otimes (f_{j}f_{j-1})\otimes\cdots\otimes f_{0} & j> 1;\\
    (-1)^{|f_0|S}(f_0f_n)\otimes f_{n-1}\otimes\cdots\otimes f_1, & j = 0,\end{cases}\]
\end{itemize}
we have that $\underline{\mathrm{CN}}_\bullet\underline{\mathcal P}$ is a cyclic object in $\mathrm{Ch}(k)$.
\end{lemma}
At the endpoints, those factors with indices outside the indicated range are omitted. 
In particular, $\bar s^\mathrm H_0(f_n\otimes\cdots\otimes f_0)=f_n\otimes\cdots\otimes f_0\otimes\id_{X_0}$, and $t_0=\id$.
For the current cyclic operator, we have $\bar d_0^\mathrm H = \bar d_n^\mathrm H t_n^\mathrm H$ and $(t_n^\mathrm H)^{-1}$ is the cyclic operator in the usual convention such that $\bar d_0^\mathrm H (t_n^\mathrm H)^{-1}= \bar d_n^\mathrm H.$

\begin{proof}
The degeneracies $\bar s_j^H$ are morphisms because the differential of the identities vanishes.
We need $d_\mathrm{int} t_n^\mathrm H=t_n^\mathrm Hd_\mathrm{int}$. For this, we analyze the summand containing $df_j$. 
Put $A_j=\sum_{i=j+1}^{n}|f_i|,$ $S=\sum_{i=1}^{n}|f_i|.$
If $j>0$, then the left-hand side has sign $(-1)^{|f_0|S+|f_0|+A_j}$. 
On the right-hand side, as $df_j$ has degree $|f_j|+1$, the sign is $(-1)^{A_j+|f_0|(S+1)},$ and the two signs match. 
If $j=0$, then the left-hand side has sign $(-1)^{|f_0|S}$, whereas the right-hand side has sign $(-1)^S(-1)^{(|f_0|+1)S}=(-1)^{|f_0|S}.$
Hence the two signs also match in this case.

We then show that $d_\mathrm{int}\bar d^\mathrm H_j=\bar d^\mathrm H_jd_\mathrm{int}$. For $1\leq j\leq n$, this follows from the Leibniz rule for the composition of maps: $d(f_jf_{j-1})=(df_j)f_{j-1}+(-1)^{|f_j|}f_j(df_{j-1}).$
For $j=0$, as $\bar d^\mathrm H_0=\bar d^\mathrm H_n t_n^\mathrm H$, $\bar d^\mathrm H_0$ is also a cochain map.

The modified operators satisfy the simplicial identities. 
The current cyclic identities are
\[(t_{n}^\mathrm{H})^{n+1} = \id,\quad \bar d_j^\mathrm{H}t_n = \begin{cases}t_{n-1}^\mathrm{H}\bar d_{j+1}^\mathrm{H}& j<n;\\\bar d^\mathrm{H}_0\quad & j = n,
\end{cases}\quad  \bar s_j^\mathrm{H}t_n^\mathrm{H} = \begin{cases}t_{n+1}^\mathrm{H}\bar s_{j+1} & j<n;\\
(t_{n+1}^\mathrm{H})^2\bar s_0^\mathrm{H} & j=n.\end{cases}\]
For $(t_n^\mathrm H)^{n+1},$ each pair of homogeneous factors crosses twice, so its total sign is $1$. 
The proof, for example, we want (For explicitness, write $\bar d_i^{(n)}$ for $\bar d_i^\mathrm{H}$, $d_i^{(n)}$ for $d_i^\mathrm{H}$ and similar notations for degeneracies)
\begin{equation*}\bar d_i^{(n+1)}\bar s_j^{n} = \begin{cases}\bar s_{j-1}^{(n-1)}\bar d_i^{(n)}& i<j;\\
\id & \substack{i=j \text{ or }\\
i = j+1};\\
\bar s_j^{(n-1)}\bar d_{i-1}^{(n)} & i>j+1.
\end{cases}\end{equation*}
Set $a \coloneqq n+1-i$ and $b \coloneqq n-j$. 
The following equations hold 
\begin{equation*}d_a^{(n+1)}s_b^{(n)}
=
\begin{cases}
s_{b-1}^{(n-1)}d_a^{(n)}
    & a<b,\\
\operatorname{id}
    & \substack{a=b \text{ or }\\ a=b+1},\\
s_b^{(n-1)}d_{a-1}^{(n)}
    & a>b+1.
\end{cases}\end{equation*}
If $i<j$, then $a>b+1$ and we have 
\[\bar d_i^{(n+1)}\bar s_j^{(n)}=d_a^{(n+1)}s_b^{(n)} = s_b^{(n-1)}d_{a-1}^{(n)}=s_{n-j}^{(n-1)}d_{n-i}^{(n)} = \bar s_{j-1}^{(n-1)}\bar d_i^{(n)}.\]
The other desired equations could be checked similarly.
For the remaining identities, the Koszul signs agree because the two sides perform the same graded permutation. 
A composition has the degree to be the sum of maps, and an identity has degree zero. 
Thus all structure maps are cochain maps, finishing the proof.
\end{proof}

The Hochschild differential is \[\bar b \coloneqq \sum_{j=0}^n(-1)^j\bar d_j^\mathrm H:\underline{\mathrm{CN}}_n(\underline{\mathcal P})\to \underline{\mathrm{CN}}_{n-1}(\underline{\mathcal P})\text{ for }n\geq 1\text{ and }\bar b_0 = 0\] and, by this lemma, commutes with $d_{\mathrm{int}}$.
We have the following additivity following \cite[Proposition~2.4.9]{McCarthy}.
\begin{proposition}\label{prop:additivity}Retain the above notations. Let $\underline{\mathcal Q}$ be a dg category satisfy the same hypothesis as $\underline{\mathcal P}$.
With the above notation, the natural
morphism of cyclic objects in $\mathrm{Ch}(k)$
\begin{align*}
& p:\underline{\mathrm{CN}}_\bullet(\underline{\mathcal P}\times\underline{\mathcal Q}) \to \underline{\mathrm{CN}}_\bullet(\underline{\mathcal P}) \oplus \underline{\mathrm{CN}}_\bullet(\underline{\mathcal Q});\\
& p((f_n,g_n)\otimes \cdots \otimes (f_0,g_0)) = (f_n\otimes \cdots \otimes f_0,g_n\otimes \cdots\otimes g_0).
\end{align*}
is a dg special homotopy equivalence.
\end{proposition}
In particular, this dg special homotopy equivalence induces a homotopy equivalence of the Hochschild total complexes.
For this proof, we need a dg version of McCarthy's special homotopy criterion \cite[Proposition~2.3.4]{McCarthy}.
\begin{lemma}\label{lemma:specialhomotopy}
Resume the notations above.
For each object $X$ of $\underline{\mathcal P}$, there are finitely many objects $D_X(i)$ indexed by $i\in I_X$ and closed degree-zero morphisms
\[\alpha_i(X):X\to D_X(i);\quad \beta_i(X):D_X(i)\to X\]
satisfying $\sum_{i\in I_X} \beta_i(X)\alpha_i(X) = \id_X.$
With $f_i: X_i\to X_{i+1}$ for $r<n$ and $f_n:X_n\to X_0$, then the map $E_\bullet = E_\bullet(\alpha,\beta)$ defined below is dg special-homotopic to the identity on $\underline{\mathrm{CN}}_{\bullet}(\underline{\mathcal P})$
\[\textstyle E_n(f_n\otimes \cdots \otimes f_0) = \sum_{i_0\in I_{X_0},\ldots,i_n\in I_{X_n}}(\alpha_{i_0}(X_0)f_n\beta_{i_n}(X_n))\otimes \cdots\otimes (\alpha_{i_{1}}(X_1)f_0\beta_{i_{0}}(X_{0})).\]
\end{lemma}

\begin{proof}
We follow the construction in \cite[Theorem 2.3.4]{McCarthy}, using the indices of Lemma~\ref{lemma:dgcyclicnerve}. 
Recall that the maps $\alpha_i(X)$ and $\beta_i(X)$ are closed of degree zero and satisfy
\[\sum_{i\in I_X}\beta_i(X)\alpha_i(X)=\id_X.\]
Let $\mathcal{S}$ be the category with two objects $0,1$ and four morphisms $\id_0,$ $\id_1,$ $0\to 1,$ and $1\to 0.$ 
We regard $k\mathcal{S}$ as a dg category concentrated in degree zero.
For morphisms $f:X\to Y$ and $g$ in $\mathcal S$, define
\[
\mu_{i,j}(f;g)=
\begin{cases}
f,&g=\id_0,\quad i=j=1,\\
\alpha_i(Y)f,&g:0\to1,\quad j=1,\\
f\beta_j(X),&g:1\to0,\quad i=1,\\
\alpha_i(Y)f\beta_j(X),&g=\id_1,\\
0,&\text{otherwise},
\end{cases}
\]
where $i$ indexes the target and $j$ the source. 

For a cycle with vertices $(X_0,\ldots,X_n)$, write $f_r:X_r\to X_{r+1}$ for $r<n$ and $f_n:X_n\to X_0$. Similarly, write $g_r:a_r\to a_{r+1}$ and $g_n:a_n\to a_0$, where $a_r\in\{0,1\}$. Put
\[
\begin{aligned}
&h_n\bigl((f_n\otimes\cdots\otimes f_0)\otimes(g_n\otimes\cdots\otimes g_0)\bigr)\\
&\qquad=\sum_{i_0,\ldots,i_n}
\mu_{i_0,i_n}(f_n;g_n)\otimes\cdots\otimes
\mu_{i_1,i_0}(f_0;g_0),
\end{aligned}
\]
with $i_r=1$ if $a_r=0$, and $i_r\in I_{X_r}$ if $a_r=1$. This defines
\[
h:\underline{\mathrm{CN}}_\bullet(\underline{\mathcal P})\otimes_k\underline{\mathrm{CN}}_\bullet(k\mathcal{S})
\longrightarrow\underline{\mathrm{CN}}_\bullet(\underline{\mathcal P}),
\]
where the tensor product is degreewise in the cyclic direction.

The relation $\sum_i\beta_i\alpha_i=\id$ shows that $h$ commutes with the faces. 
Since each nonzero $\mu_{i,j}(f;g)$ has degree $|f|$, the cyclic operators on both sides produce the same Koszul sign. 
Thus $h$ is semi-cyclic. 
Its restrictions to the constant $0$- and $1$-cycles are $\id$ and $E_\bullet$, respectively.

Finally, closedness and degree zero of the $\alpha_i,\beta_i$ give $d\mu_{i,j}(f;g)=\mu_{i,j}(df;g).$
Hence $d_\mathrm{int} h=hd_\mathrm{int}$ and $d_\mathrm{int} E=Ed_\mathrm{int}$. This proves that $h$ is a dg special homotopy. We now spell out the associated Hochschild homotopy with our indexing.

For $0\leq j\leq n$, let $\bar\lambda_j^{(n)}:[n+1]\to[1]$ send $0,\ldots,n-j$ to $0$ and $n-j+1,\ldots,n+1$ to $1$. 
This map may written as a tuple and corresponds to an element $\Lambda_j^{(n)}$ in $\underline{\mathrm{CN}}_{n+1}(k\mathcal S)$
\[\bar\lambda_j^{(n)} = (\underbrace{0,\ldots,0}_{n-j+1},\underbrace{1,\ldots,1}_{j+1})\text{ and }\Lambda_j^{(n)} = (0\to 1)\otimes \id_0^{\otimes (n-j)} \otimes (1\to 0) \otimes \id_{1}^{\otimes j}.\]
Put
\[
\bar h_j^{(n)}:\underline{\mathrm{CN}}_n(\underline{\mathcal P})\rightarrow\underline{\mathrm{CN}}_{n+1}(\underline{\mathcal P});
\quad
\bar h_j^{(n)}(x)
=h_{n+1}(\bar s^\mathrm H_jx\otimes \Lambda_j^{(n)}).
\]
For any $x$, the compatibility of $h$ with the faces gives $\bar d^\mathrm H_i\bar h_j^{(n)}(x) = h_{n}(\bar d^\mathrm H_i\bar s^\mathrm H_jx\otimes \bar d^\mathrm H_i\Lambda_j^{(n)})$.
The simplicial identities and compatibility of $h$ with the faces implies
\begin{align*}
&\bar d^\mathrm H_i\bar h_j^{(n)}=\begin{cases}\bar h_{j-1}^{(n-1)}\bar d^\mathrm H_i & i<j;\\\bar h_j^{(n-1)}\bar d^\mathrm H_{i-1}&i>j+1,\end{cases}\,\,\bar d^\mathrm H_i\bar h_i^{(n)}=\bar d^\mathrm H_i\bar h_{i-1}^{(n)}\text{ for }1\leq i\leq n,\\
&\text{and }\bar d^\mathrm H_0\bar h_0^{(n)}=\id,\,\,\bar d^\mathrm H_{n+1}\bar h_n^{(n)}=E_n.
\end{align*}
Put $\textstyle H_n:=-\sum_{j=0}^{n}(-1)^j\bar h_j^{(n)}\text{ and }H_{-1}=0.$ 
Then $\bar b_{n+1}H_n+H_{n-1}\bar b_n=E_n-\id.$

Since $h$ and the degeneracies commute with $d_\mathrm{int}$, so does $H$. On the sum total complex, an element $x\in\underline{\mathrm{CN}}_n(\underline{\mathcal P})^q$ has total degree $q-n$. With $D(x)=d_\mathrm{int} x+(-1)^q\bar b x$ and $K(x)=(-1)^qH(x)$, we obtain $DK+KD=E-\id$. 
This is the desired homotopy on the Hochschild total complex.
\end{proof}

\begin{proof}[Proof of Proposition~\text{\rm{\ref{prop:additivity}}}]
The proof idea is to construct a quasi-inverse $\iota:\underline{\mathrm{CN}}_\bullet(\underline{\mathcal P})\oplus \underline{\mathrm{CN}}_\bullet(\underline{\mathcal Q}) \to \underline{\mathrm{CN}}_\bullet(\underline{\mathcal P}\otimes \underline{\mathcal Q})$
so that $p\iota = \id$ and $\iota p = E_\bullet$ for some degree zero closed homogeneous map $\alpha_i$, $\beta_i$ and the Lemma~\ref{lemma:specialhomotopy} implies the proposition.

Two dg functors $\underline{\mathcal P}\to \underline{\mathcal P}\times \underline{\mathcal Q};\,\,X \mapsto (X,0)\text{ and }\underline{\mathcal Q}\to \underline{\mathcal P}\times \underline{\mathcal Q};\,\,Y \mapsto (0,Y)$ induce the desired $\iota$: 
\begin{equation}\begin{split}
&\iota:\underline{\mathrm{CN}}_\bullet(\underline{\mathcal P})\oplus \underline{\mathrm{CN}}_\bullet(\underline{\mathcal Q})\to \underline{\mathrm{CN}}_\bullet(\underline{\mathcal P}\times\underline{\mathcal Q});\\
&(f_n\otimes\cdots\otimes f_0,g_n\otimes \cdots\otimes g_0)\mapsto (f_n,0)\otimes\cdots\otimes (f_0,0) + (0,g_n)\otimes\cdots\otimes (0,g_0).
\end{split}\label{formula:quasiinverseforadditivity}
\end{equation}
The desired data for Lemma~\ref{lemma:specialhomotopy} are $D_{(X,Y)}(1) = (X,0)$ and $D_{(X,Y)} = (0,Y)$ with the canonical projections and inclusions 
\[\alpha_1:(X,Y)\to (X,0),\,\,\beta_1:(X,0)\to(X,Y),\,\,\alpha_2:(X,Y)\to(0,Y),\text{ and }\beta_2:(0,Y)\to (X,Y)\]
that is closed of degree zero and satisfies $\beta_1\alpha_1 + \beta_2\alpha_2 = \id_{(X,Y)}.$
For the map $E_n$ defined in Lemma~\ref{lemma:specialhomotopy}, we have
all summands equal to zero except for $i_0 = \cdots = i_n = 1$ or $i_0 = \cdots = i_n = 2$ and hence 
\[E_n((f_n,g_n)\otimes\cdots\otimes(f_0,g_0))=(f_n,0)\otimes\cdots\otimes (f_0,0) + (0,g_n)\otimes\cdots\otimes (0,g_0)\]
Thus $E_\bullet = \iota p,$ as desired.
\end{proof}

\end{document}